\documentclass[11pt,a4paper]{article}
\usepackage{graphicx}
\usepackage{enumerate}
\usepackage{url} 
\RequirePackage[OT1]{fontenc}
\usepackage{bbm}
\usepackage[numbers]{natbib}
\usepackage{amsthm,amsmath}
\usepackage{graphicx,ifthen}
\usepackage{latexsym}
\usepackage{mathrsfs}
\usepackage{amssymb}
\usepackage{hyperref}
\usepackage{xcolor}
\usepackage{caption, subcaption}
\usepackage{stmaryrd}
\usepackage{circledsteps}
\usepackage{dsfont}

\newtheoremstyle{example}{\topsep}{\topsep}%
     {}
     {}
     {\rmfamily}
     {}
     {\newline}
     {\thmname{#1}\thmnumber{ #2}\thmnote{ #3}}

   \theoremstyle{example}

\numberwithin{equation}{section}
\theoremstyle{plain}
\newtheorem{thm}{Theorem}[section]

\newtheorem{exam}{Example}[subsection]
\newtheorem{prop}{Proposition}[section]

\newtheorem{rem}{Remark}[section]

\newtheorem{cor}{Corollary}[section]

\newcommand{\Lower}[2]{\smash{\lower #1 \hbox{#2}}}
\newcommand{\ben}{\begin{enumerate}}
\newcommand{\een}{\end{enumerate}}
\newcommand{\bi}{\begin{itemize}}
\newcommand{\ei}{\end{itemize}}

\def\[#1\]{\begin{equation}\begin{aligned}#1\end{aligned}\end{equation}}
\newcommand{\given}{{\hspace{0.08em}|\hspace{0.08em}}}
\newcommand{\tC}{\widetilde{C}}
\newcommand{\tN}{\widetilde{N}}
\newcommand{\tY}{\widetilde{Y}}
\newcommand{\bigcdot}{{\boldsymbol{\cdot}}}
\newcommand{\bakappa}{\kappa_{\boldsymbol{\cdot}}}
\newcommand{\bX}{\mathbf{X}}
\newcommand{\bZ}{\mathbf{Z}}
\newcommand{\bn}{\mathbf{n}}
\newcommand{\bbP}{\mathbb{P}}
\newcommand{\1}[1]{\mathds{1}_{\{#1\}}}
\newcommand{\calN}{\mathcal{N}}

\newcommand{\stirlingfirst}[2]{\genfrac{[}{]}{0pt}{}{#1}{#2}}

\begin{document}

\def\spacingset#1{\renewcommand{\baselinestretch}%
{#1}\small\normalsize} \spacingset{1}



\title{\bf Bayesian Analysis of Generalized Hierarchical Indian Buffet Processes for Within and Across Group Sharing of Latent Features}
\author{Lancelot F. James\thanks{
\textit{Supported in part by grants RGC-GRF 16301521, 16300217, and T31-604/18-N of the HKSAR.}}\hspace{.2cm}\\
Department of ISOM, HKUST\\
Juho Lee
\\
The Graduate School of AI, KAIST\\
Abhinav Pandey \\
Department of ISOM, HKUST\\
}
\maketitle

\begin{abstract}
Bayesian nonparametric hierarchical priors are highly effective in providing flexible models for latent data structures exhibiting sharing of information within and across groups. In this work, we focus on latent feature allocation models, where the data structures correspond to multi-sets or unbounded sparse matrices, which we refer to as generalized hierarchical Indian Buffet processes (HIBP). These are based on hierarchical versions of generalized spike and slab Indian Buffet processes (IBP), where the fundamental development in this regard is the Bernoulli-based HIBP, devised by Thibaux and Jordan (2007), as a hierarchical extension of the IBP devised by Griffiths and Ghahramani (2005). With a focus on Bayesian inference, we provide novel explicit descriptions of the joint, marginal, and posterior distributions of the HIBP, significantly advancing our understanding of these processes. Our results allow for exact sampling for the otherwise complex joint marginal distributions. We provide a general characterization of their posterior distributions as well as highlight bottlenecks for practical implementation. Our main focus then shifts to specific tractable results for the remarkable case of Poisson HIBP, which correspond to generalizations of mixed Poisson random count models arising in genetics, imaging, topic modeling, random occupancy, and species sampling models. We show they also have important relations to Bayesian nonparametric latent class models appearing in the literature. Furthermore, we show that all general HIBP may be coupled to Poisson HIBP, allowing for further analysis of such processes.
\end{abstract}

\noindent%
\vfill

\newpage


\section{Introduction}
Bayesian nonparametric hierarchical priors are highly effective in providing flexible models for latent data structures exhibiting sharing of information between and across groups. 
Most prominent is the Hierarchical Dirichlet Process (HDP) \citep{HDP}, and its subsequent variants, which model latent clustering/class between and across groups.  Note in these models individuals share at most one feature within and across groups. 

In this work, we focus on tractable Bayesian analysis for hierarchical notions of latent feature allocation models, where the data structures correspond to multi-sets or unbounded sparse matrices. The fundamental development in this regard is the Hierarchical Indian Buffet Process (HIBP), devised in \citep{Thibaux}. The HIBP utilizes a hierarchy of Beta processes to model the sharing of latent features across $J$ distinct groups, where each individual group $j$ (for $j = 1, 2,\ldots,J$) consists of $M_{j}$ customers. Specifically, each group $j$ generates a latent binary random matrix, reflecting the within-group sharing of features, according to a beta-Bernoulli Indian Buffet Process (IBP) prior first devised in \citep{GriffithsZ}. Similar to the HDP, this information is merged/shared across groups according to another common Beta process which is an example of a Completely Random measure (CRM).

Here we extend this to hierarchical notions of the general spike and slab IBP referred to as generalized HIBP, which include models already appearing in the literature. Most notably our focus will be on Poisson based HIBP of which, single customer,  Poisson-Gamma-Gamma or equivalently Negative Binomial-Gamma random count models treated in \citep{ZhouCarin2015, ZhouPadilla2016}, are notable special cases that otherwise have  rich modelling capabilities, and properties.  In general, the Poisson HIBP  generates complex random count matrices based on multivariate mixed Poisson processes and correspond to generalizations of models arising in genetics, imaging, topic modeling, random occupancy, and species sampling models. We show they also have important relations to Bayesian nonparametric latent class models appearing in the literature. A key feature of our exposition is that these connections allow for distributional recovery of individual latent counts that would otherwise be observables in the non-hierarchical setting.  Furthermore, we show that all general HIBP may be coupled to Poisson HIBP, allowing for further analysis of such processes. 

The basic Indian Buffet process (IBP) first devised in \citep{GriffithsZ}, constitutes a class of flexible, and tractable, Bayesian non-parametric priors for a wide range of statistical applications involving latent features. See for example \citep{MasoeroBiometrika}, for a recent interesting application to problems in genetics. See ~\citep{Broderick1} for more on the IBP sequential scheme for allocation of features to individual customers, originally described in~\citep{GriffithsZ},  and comparisons with general Chinese restaurant processes for latent class models. The IBP has inspired the development and investigation of generalizations which induce non-binary entries. Notably,~\citep{Titsias} pointed out that images, which play the analogous roles of customers, exhibit multiple occurrences of a feature within a single image, and hence in analogy to the classic gamma-Poisson Bayesian parametric model, proposed the use of Poisson-type IBP with a gamma process prior in that setting. Other possibilities for distributions have been treated in the literature. In particular \citep{James2017} provides a unified framework for the posterior analysis of  general \textit{spike and slab} IBP type models, and important multivariate extensions, which produce sparse random matrices with arbitrary valued entries based on random variables having mass at zero, but otherwise producing non-zero entries corresponding to any random variable. In addition, see~\citep{Basbug} for potential ideas in regards to matrices with general spike and slab entries. In this regard it is notable that  models with general compound Poisson entries, and also complex random count matrices, shall arise naturally within the context of our work. Also the work of \citep{CamerFavScaled}, describe modifications of the Bernoulli IBP utilizing, instead of CRM,  scaled subordinators~\citep{JOT} for better modelling of unseen-features/species models arising for instance in genetics. This modification may be easily applied to our present exposition although we refrain from doing so here. In particular,~\citep{BerahaFavscaled} discuss this modification for general IBP and hierarchical latent feature models proposed by \citep{Masoerotrait} and also the HIBP processes we analyze here. Their view of the HIBP is based on less developed properties which we aim to rectify here. 
Hereafter, for any integer $n,$ we will use the notation $[n]:= \{1, 2, ..., n\}$ to represent the set of integers from $1$ to $n$.
\subsection{Outline and contributions}
As already discussed in for instance the works of~\citep{Masoerotrait,Thibaux, ZhouPadilla2016} the extension of IBP type processes to processes sharing multiple features over groups has a vast number of possible applications. What remains is a detailed analysis, such as in the case of the Dirichlet process~\cite{Ferg1973,Pit96} and all subsequently widely used BNP models, that will allow for a better understanding of their structural/distributional properties hence allowing practical modelling for a wide range of applications and aid in the development of efficient computational methods for inference. In addition, it is important to establish connections to existing models that may arise in other contexts. Our aim is to do this for the HIBP models we describe here. Specifically in section 2 we begin with more details of the definitions of IBP and HIBP and their building blocks based on completely random measures. We first revisit the key results from \citep{James2017} in the general IBP setting, in~Proposition~\ref{IBPpost} highlighting a more accessible perspective in terms of a method of decompositions. We then use this along with basic properties of completely random measures to provide initial representations of the HIBP, and then proceed to provide analysis of a key, otherwise hidden, multivariate Poisson process dictating the allocation of features across groups. Notably this process can be seen as a multivariate extension of processes appearing in~\cite{Pit97}.  Section 3 then establishes a result to exactly sample from the marginal distribution of any HIBP process. Section 3.1 provides specific details for Poisson HIBP and importantly show any HIBP can be coupled to the Poisson framework. More variables, with interpretations gained in~\cite{Pit97}, are introduced and will form the basis for a method of distributional recovery of key quantities for posterior analysis. Section 4 presents a general description of the posterior distribution for all HIBP, highlighting general stumbling blocks in terms of recovery of latent quantities that would normally be observed in the non-hierarchical setting. Rather than attempting a general discussion of how to address this issue, the remainder of our work focuses on achieving full tractability in the case of general mixed Poisson processes, that is the Poisson HIBP. We also note that the fact that we can couple general HIBP to this framework allows one to obtain further details in those cases. Briefly, focusing on the Poisson HIBP, in section 5 we obtain tractable descriptions of the posterior and predictive descriptions and a description of the otherwise complex marginal distribution which one can see, due to our representations, as sums over Gibbs (product type) structure in relation to permanental processes,~see for instance~\citep{CraneEwens}, appearing in the study of multivariate negative binomial and more general processes, which are otherwise difficult to sample from. We also describe a highy tractable simplification for HIBP which otherwise exhibits a rich structure for modelling with perhaps surprising connections to Pitman-Yor and finite dimensional versions of such processes as arising in recent work of~\cite{LijoiPYmulti,LijoiPYmultiJASA}, where we provide new ways to exactly sample from the corresponding Pitman-Yor multinomial process, among other properties. 

 We note further the mixed Poisson framework has firm connections to the species sampling framework as it relates to latent class/clustering models. In particular, Pitman   
~\citep{PitmanPoissonMix}, in an unpublished manuscript, shows that the mixed Poisson process has an interpretation, already within the framework of Fisher/McCloskey, in terms of the number of customers/animals sampled/trapped at a given time and the corresponding distribution of the time when this occurs, see\citep{HJL}[Section 2.2.1. p. 321] for some more related details. In particular, relevant to our exposition, he points out connections between the work of~\cite{FoFZhou} on what is equivalently in distribution a single Poisson IBP, that is a mixed Poisson process,  to the use of the gamma-trick in~\cite{JLP2}. In fact the variables, corresponding to the gamma trick, can be directly interpreted as the corresponding arrival times of the mixed Poisson process. This plays an important role leading to one of our key results~in Proposition~\ref{prop:GibbsEPPF}. Our work, which involves $J$ groups, may be used to extend these notions with broader implications that go beyond the scope of this exposition.   

\section{Generalized HIBP }
We now discuss the formal definition of generalized HIBP, by first recalling the general $\mathrm{IBP}$ as 
discussed in~\cite{James2017}. We start by recalling definitions of completely random measures 
highlighting some properties that are of relevance to our exposition when the base measure $B_{0}$ 
defined below is discrete.  Let, for each $j\in[J],$ $\mu_{j}$ be a completely random measure defined as 
$\mu_{j} = \sum_{k=1}^{\infty} s_{j,k} \delta_{w_{j,k}} \sim \mathrm{CRM}(\rho_{j},B_{0})$, where 
$\rho_{j}(s)dsB_{0}(d\omega)$ is the mean measure of a Poisson random measure, see~\cite{James2017} 
for more on this relation, and specifies the distribution of $(s_{j,k},w_{j,k})_{\{k\ge 1\}}$ as follows: Here, $B_{0}(dw)$ is a finite measure on a Polish space $\Omega$ such 
that $(w_{j,k})_{\{k\ge 1\}} \overset{iid}\sim \bar{B}_{0}(dw) := B_{0}(dw)/B_{0}(\Omega)$, and 
$\rho_{j}(s)$ is the L\'evy density on $(0,\infty)$ satisfying 
$\int_{0}^{\infty}\min(s,1)\rho_{j}(s)ds<\infty$. The $\rho_{j}$ determines the distribution of the jumps of $\mu_{j},$ represented by 
$(s_{j,k})_{\{k\ge 1\}}$, and $\mu_{j}(\Omega):=\sum_{k=1}^{\infty}s_{j,k}$ is a positive 
infinitely divisible continuous random variable with Laplace transform 
$\mathbb{E}[e^{-t\mu_{j}(\Omega)}] = e^{-B_{0}(\Omega)\psi_{j}(t)}$, and Laplace exponent
$\psi_{j}(t) := \int_{0}^{\infty}(1-e^{-ts})\rho_{j}(s)ds.$ Furthermore, 
for a countable number of disjoint sets $(Q_{j,l})_{l\ge 1}$, $\mu_{j}(Q_{j,l})$ are independent with 
respective laws indicated by the Laplace exponents $B_{0}(Q_{j,l})\psi_{j}(t)$ for $l=1,2,\ldots.$

\begin{rem}\label{mudecomp}
The latter statement of independence over sets, which is a defining property of $\mathrm{CRM}$ applies to both continuous and discrete cases of $B_{0},$ and will
make more transparent some of our initial representations and derivations in the latter case. Consequently, for each $j$, the jumps $\{s_{j,k} : w_{j,k} \in Q_{j,l}, k\ge 1\}$, for $l\in\{1,2,\ldots\}$ countable number of disjoint sets, are independent and specified by Lévy densities $B_{0}(Q_{j,l})\rho_{j}$, and the atoms are selected iid proportional to $\mathbb{I}_{\{\omega\in Q_{j,l}\}}B_{0}(d\omega)$ $l=1,2,\ldots$. 
\end{rem}

Now, following~\cite{James2017}, we define the generalized spike and slab IBP for each group $j$ of $i\in[M_j]$ customers as: $Z^{(i)}_j = \sum_{k=1}^{\infty}A^{(i)}_{j,k}\delta_{w_{j,k}}|\mu_j \overset{iid}{\sim} \mathrm{IBP}_{A_j}(\mu_j),$ where $\mu_j \sim \mathrm{CRM}(\rho_j,B_0)$. Here, $\mathrm{IBP}_{A_j}$ denotes the processes where for each $j\in[J],$ $A^{(i)}_{j,k}|s_{j,k}=s$ are iid equal in distribution to a spike and slab variable $A_j$ with distribution: $\mathbb{P}_{A_j}(da|s) = [1-\pi_j(s)]\delta_0(da) + \mathbb{I}_{a \neq 0}\mathbb{P}_{A_j}(da|s),$ with $\mathbb{P}(A_j=0|s) = 1-\pi_j(s)$ such that $\int_0^{\infty}(1-[1-\pi_j(s)]^{M_j})\rho_j(s)ds$ is finite. Hence $\mathbb{P}((A^{(i)}=0,i\in[M_{j}]|s)={[1-\pi_j(s)]}^{M_j}.$ Note $A_{j} \overset{d}{=} b_j\hat{A}_j$, with $b_j \sim \text{Bernoulli}(\pi_j(s))$ and $\hat{A}_{j}|s$ has zero-truncated distribution $\mathbb{I}_{a \neq 0}\mathbb{P}_{A_j}(a|s)/\pi_{j}(s).$ As special cases $A_j \overset{d}{=} b_j \sim \text{Bernoulli}(s)$ yields $Z^{(i)}_j|\mu_j \overset{iid}{\sim} \text{BerP}(\mu_j),$ denoting a Bernoulli process. In this case $\hat{A}_{j}=1.$ If $A_{j}\sim \mathrm{Poisson}(r_{j}s),$ then  $Z^{(i)}_j|\mu_{j}\sim \mathrm{PoisP}(r_j\mu_{j})$ denoting  a Poisson porcess as in ~\cite{Titsias, FoFZhou},  see~\cite{James2017}[Section 4] for more details on these and other cases.  

As with the classic Bernoulli case~\cite{Thibaux} and Poisson-Gamma-Gamma HIBP~\cite{ZhouCarin2015,ZhouPadilla2016} cases discussed in the literature, we create general HIBP models by specifying $B_0$ to be a discrete random measure  fostering sharing across the $J$ groups.  Here specifically: $B_0=\sum_{l=1}^{\infty}\lambda_l\delta_{Y_{l}} \sim \mathrm{CRM}(\tau_0,F_0)$ where $F_0$ is a proper probability distribution, assumed to be non-atomic, on $\Omega$ and $\mathbb{E}[e^{-tB_0(\Omega)}]=e^{-\Psi_{0}(t)}$, where $\Psi_{0}(t)=\int_{0}^{\infty}(1-e^{-t\lambda})\tau_{0}(\lambda)d\lambda$. This gives the general HIBP as, for $i\in[M_{j}], j\in[J]$:
$Z^{(i)}_j|\mu_j \overset{iid}{\sim} \mathrm{IBP}_{A_j}(\mu_j), \mu_{j}|B_{0}\sim \mathrm{CRM}(\rho_{j},B_{0}), B_{0}\sim \mathrm{CRM}(\tau_{0},F_{0}).$ We note the Negative Binomial-Gamma process in ~\cite{ZhouPadilla2016}, is the Poisson HIBP case where $Z^{(1)}_j|\mu_{j}\sim \mathrm{PoisP}(r_j\mu_{j})$, $\mu_{j}|B_{0}$ are Gamma $\mathrm{CRM}$ and $B_{0}$ is another Gamma CRM. However in that setting $(M_{j}=1, j\in[J]).$ Interest is in the process $({Z}^{(1)}_{j},j\in[J])$ and growing a random matrix from $J$ to $J+1$ etc. via the predictive process $Z_{J+1}|(\hat{Z}^{(i)}_{j},j\in [J]).$ Similarly, although not obvious and perhaps not known, if $Z^{(1)}_{j}\sim~\mathrm{BerP}(\mu_{j}),$ then $(Z^{(1)}_{j},\in[J])$ equates in distribution to processes constructed in~\cite{Caron2012}, where interest again is growth in $J=1,2,\ldots.$ We will show, both these special cases have exact implementable solutions for their posterior and predictive distributions under general choices of $\mu_{j}.$ We note throughout that we use $\mathrm{Multi}$ to denote a Multinomial distribution and $\mathrm{tP}$ to denote a zero-truncated Poisson or multivariate Poisson distribution, furthermore let \( \mathbf{0} \) represent the vector of all zeros, where the dimension may vary depending on the context.

\subsection{Properties}
We now proceed to provide descriptions of the relevant marginal, posterior and predictive distributions. We first begin with an interpretation of some of the relevant results in \cite{James2017}[Section 3] for non-hierarchical setting which will also help us introduce some key quantities and distributions, where $B_{0}$ is considered to be given. Our description below based on decomposing the relevant space into a set where the $M_{j}$ customers do not sample from and its complement, where $\mu_{j}$ has a compound Poisson component, which identifies the appropriate marginals, joint and conditional distribution of variables. We believe this leads to a more clear description of the results and underlying mechanisms in~\cite{James2017}, which can be applied more broadly.

\begin{prop}\label{IBPpost}\cite{James2017}[Section 3] For each $j,$ consider $Z^{(i)}_j |\mu_j \overset{iid}{\sim} \mathrm{IBP}_{A_j}(\mu_j)$, for $\mu_{j}\sim \mathrm{CRM}(\rho_{j},B_{0}).$ Set 
$\Phi_{j}(M_{j})=\int_{0}^{\infty}(1-{[1-\pi_{j}(s)]}^{M_{j}})\rho_{j}(s)ds.$
The event $\{k: A^{(i)}_{j,k}=0, i\in[M_{j}]\}$, with $\mathbb{P}(A^{(1)}_{j}=0,\ldots,A^{(M_{j})}_{j}=0|s)=[1-\pi_{j}(s)]^{M_{j}}$
induces a decomposition of the  L\'evy density $\rho_{j}$ of $\mu_{j}$ as:
\begin{enumerate}
    \item $\rho_{j}(s)=\rho^{(M_{j})}_{j}(s)+\Phi_{j}(M_{j})f_{S}(s)$
    where $\rho^{(M_{j})}_{j}(s)=[1-\pi_{j}(s)]^{M_{j}}\rho_{j}(s)$
    \item $f_{S_{j}}(s)=(1-{[1-\pi_{j}(s)]}^{M_{j}})\rho_{j}(s)/\Phi_{j}(M_{j})$ is proper a probability density function (pdf).
    \item Hence there is a decomposition, 
\begin{equation}
\label{mudecomp}
\mu_{j}\overset{d}{=}\mu_{j,M_{j}}+\sum_{\ell=1}^{\xi_{j}}S_{j,\ell}\delta_{\omega_{j,\ell}}
\end{equation}
where $\mu_{j,M_{j}}\sim \mathrm{CRM}(\rho^{(M_{j})}_{j},B_{0})$,
 $\xi_{j}\sim \mathrm{Poisson}(\Phi_{j}(M_{j})B_{0}(\Omega))$,
the $(\omega_{j,\ell})$ are iid $\bar{B}_{0}$ and the $(S_{j,\ell})$ are iid with pdf $f_{S_{j}}(s).$
\item For each $\ell\in[\xi_{j}]$  there are independent vectors $(\tilde{A}^{(1)}_{j,\ell},\ldots, \tilde{A}^{(M_{j})}_{j,\ell})|S_{j,\ell}=s$ with distribution $\mathbb{I}_{(\mathbf{a}_{j,\ell}\neq \mathbf{0})}\prod_{i=1}^{M_{j}}\mathbb{P}_{A_{j}}(da^{(i)}_{j,\ell}|s)/[1-{(1-\pi_{j}(s))}^{M_{j}}].$  There is the joint distribution of $(\tilde{A}^{(1)}_{j,\ell},\ldots, \tilde{A}^{(M_{j})}_{j,\ell}),S_{j,\ell}$
\begin{equation}
\mathbb{I}_{(\mathbf{a}_{j,\ell}\neq \mathbf{0})}\prod_{i=1}^{M_{j}}\mathbb{P}_{A_{j}}(da^{(i)}_{j,\ell}|s)\rho_{j}(s)/\Phi_{j}(M_{j})
\label{joint1}
\end{equation}
\item Hence integrating over $\rho_{j}$ in~\eqref{joint1} leads to the marginal distribution, with distribution denoted  $\mathrm{IBP}(A^{(i)}_{j},i\in[M_{j}],\rho_{j},B_{0}),$
\begin{equation}
\label{item:key sum}
(Z^{(i)}_{j}, i\in[M_{j}])\overset{d}=(\sum_{\ell=1}^{\xi_{j}}\tilde{A}^{(i)}_{j,\ell}\delta_{\omega_{j,\ell}},i\in[M_{j}]),
\end{equation}
where $Z^{(i)}_{j}(\Omega)\overset{d}=\sum_{\ell=1}^{\xi_{j}}\tilde{A}^{(i)}_{j,\ell},$ component-wise and jointly.
\item The posterior distribution of $\mu_{j}|(Z^{(i)}_{j}, i\in[M_{j}])$ is described by the decomposition in~\eqref{mudecomp}, fixing the values of $\xi_{j}$ and $(\omega_{j,\ell},\ell\in [\xi_{j}])$ and for each  $\ell$ applying Bayes rule to
\eqref{joint1} to obtain the conditional density of $S_{j,\ell}|(\tilde{A}^{(1)}_{j,\ell},\ldots, \tilde{A}^{(M_{j})}_{j,\ell}).$ The law of $\mu_{j,M_{j}}\sim \mathrm{CRM}(\rho^{(M_{j})}_{j},B_{0})$ remains unchanged as no points are selected from that process.
\end{enumerate}
\end{prop}

\subsection{Descriptions for HIBP}
As an interesting technical matter, the result itself relies on the unicity of the jumps $S_{j,\ell}$ and not the concomittants $\omega_{j,\ell}$ drawn from $B_{0}.$ So the descriptions above hold whether $B_{0}$ is discrete or continuous. Now, in our HIBP setting,  when $B_{0}=\sum_{l=1}^{\infty}\lambda_{l}\delta_{Y_{l}}\sim \mathrm{CRM}(\tau_{0},F_{0}),$ and hence both discrete and random, we can use Proposition~\ref{IBPpost} and Remark~\ref{mudecomp}, setting for each $j,$ $Q_{j,l}:=\{\omega :\omega=Y_{l}\},$ with $B_{0}(Q_{j,l})=\lambda_{l},$ to obtain more refined and key important representations. That is: 1. $\mu_{j}|B_{0}\overset{d}=\sum_{l=1}^{\infty}\sigma_{j,l}(\lambda_{l})\delta_{Y_{l}}$ where, 2. $\sigma_{j,l}(\lambda)$ has distribution given by the Laplace exponent $\lambda \psi_{j}(t),$ with density denoted as $\eta_{j}(t_{j}|\lambda),$ furthermore $\sigma_{j,l}(\lambda)=\sum_{k=1}^{\infty}s_{j,k,l},$ where $(s_{j,k,l})_{\{k\ge 1\}}$ are the jumps of $\mu_{j}$ in an interval of length $\lambda$ with L\'evy density $\lambda\rho_{j}(s).$ 3. $\xi_{j}\overset{d}=\sum_{l=1}^{\infty}\xi_{j,l},$ for $\xi_{j,l}\sim \mathrm{Poisson}(\Phi_{j}(M_{j})\lambda_{l})$ 4. $(Z^{(i)}_{j},i\in[M_{j}])\overset{d}=(\sum_{l=1}^{\infty}[\sum_{k=1}^{\xi_{j,l}}\tilde{A}^{(i)}_{j,k,l}]\delta_{Y_{l}},i\in[M_{j}])$ 5. Where $((\tilde{A}^{(i)}_{j,k,l}, i\in[M_{j}]),S_{j,k,l})$ has the same distribution in~\eqref{joint1}. Use~\eqref{item:key sum}.

The sum representations, in particular the summing of the jumps of $(\mu_{j},j\in[J])$ are analogous to what happens to the probability masses of random probability measures in the HDP type settings, indicative of hidden random clustering mechanisms. We shall see the infinite sum representations do not pose any significant difficulties just as in the IBP settings. Furthermore, we note that the entries $(\sum_{k=1}^{\xi_{j,l}}\tilde{A}^{(i)}_{j,k,l},i\in[M_{j}])$ are compound Poisson process vectors and hence infinitely divisible. Key to our exposition, and otherwise hidden, is the multivariate Poisson process $\mathbf{X}:=(\sum_{l=1}^{\infty}\xi_{j,l}\delta_{Y_{l}},j\in[J])$ which dictates the features and their amounts distributed across each group $j\in[J]$.
In the next result, we describe properties of this multivariate Poisson $\mathrm{IBP}$ process, which, albeit not obvious, plays a crucial role throughout, as it connects as a multivariate extension of a discretization process in~\cite{Pit97}, see also~\cite{James2017}[Remark 4.2]. Our approach is again to use a description in terms of a method of decompositions, but otherwise relies on the results in~\cite{James2017}[Section 5] for multivariate IBP.

\begin{prop}\label{prop2}Consider $\mathbf{X}\overset{d}=(\sum_{l=1}^{\infty}\xi_{j,l}\delta_{Y_{l}}, j\in[J]).$ Then the event $\{l: \xi_{j,l}=0, j\in[J]\}$, with $\mathbb{P}((\xi_{j,l}=0,j\in[J])|\lambda)={\mbox e}^{-\lambda\sum_{j=1}^{J}\Phi_{j}(M_{j})}$
induces a decomposition of the  L\'evy density $\tau_{0}$ of $B_{0}$ as: 
\begin{enumerate}
    \item $\tau_{0}(\lambda)={\mbox e}^{-\lambda\sum_{j=1}^{J}\Phi_{j}(M_{j})}\tau_{0}(\lambda)+\Psi_{0}(\sum_{j=1}^{J}\Phi_{j}(M_{j}))f_{H}(\lambda)$.
    \item Where $f_{H}(\lambda)=\frac{(1-{\mbox e}^{-\lambda\sum_{j=1}^{J}\Phi_{j}(M_{j})})\tau_{0}(\lambda)}{\Psi_{0}(\sum_{j=1}^{J}\Phi_{j}(M_{j}))}$ is a proper density.
    \item Hence setting $\tau_{0,J}(\lambda)={\mbox e}^{-\lambda\sum_{j=1}^{J}\Phi_{j}(M_{j})}\tau_{0}(\lambda),$ this leads to a decomposition of $B_{0}:$
  \begin{equation}
    \label{postdisint}
B_{0}\overset{d}=B_{0,J} + \sum_{l=1}^{\varphi} H_{l}\delta_{\tilde{Y}_{l}}
    \end{equation}   
where $B_{0,J}\sim \mathrm{CRM}(\tau_{0,J},F_{0}),$  $(\tilde{Y}_{l})$ are iid $F_{0}$ and represent the $\varphi \sim \mathrm{Poisson}(\Psi_{0}(\sum_{j=1}^{J}\Phi_{j}(M_{j})))$ number of distinct features to be selected by the $\sum_{j=1}^{J}M_{j}$ customers. The $(H_{l},l\in[\varphi])\overset{iid}\sim f_{H}$ are the corresponding jumps of $B_{0}.$
\item Sampling from $\mathbf{X}$ leads to the selection of $\varphi$ iid vectors, where for each $l\in[\varphi],$ they correspond to,
$(X_{j,l}=x_{j,l}, j\in[J])|H_{l}=\lambda$ a multivariate zero-truncated Poisson distribution denoted as $\mathrm{tP}(\lambda (\Phi_{1}(M_{1}),\ldots,\Phi_{J}(M_{J}))$. Where $x_{l}:=\sum_{j=1}^{J}x_{j,l}=1,2,\ldots$
\item Hence $\tilde{X}_{l}:=\sum_{j=1}^{J}X_{j,l}|H_{l}=\lambda\sim \mathrm{tP}(\lambda \sum_{j=1}^{J}\Phi_{j}(M_{j}))$ with marginal distribution denoted as $\mathrm{MtP}(\sum_{j=1}^{J}\Phi_{j}(M_{j}),\tau_{0})$
\item $(X_{j,l},j\in[J])|\tilde{X}_{l}=x_{l}\sim \mathrm{Multi}(x_{l}; q_{1},\ldots, q_{J})$, for $q_{j}=\frac{\Phi_{j}(M_{j})}{\sum_{v=1}^{J}\Phi_{v}(M_{v})}.$
\end{enumerate}
\end{prop}
We next describe the posterior distribution of $B_{0}|\mathbf{X}$ which will play several important roles in the sequel. We see it's posterior has the same form as the univariate Poisson IBP in~\cite{James2017}[Section 4.2], as it depends only on 
$(\sum_{l=1}^{\infty}[\sum_{j=1}^{J}\xi_{j,l}]\delta_{Y_{l}})|B_{0}\sim \mathrm{PoiP}([\sum_{j=1}^{J}\Phi_{j}(M_{j})]B_0).$
\begin{prop}\label{propN} The posterior distribution of $B_{0}|\mathbf{X}$ only depends on the $(\tilde{X}_{l}=x_{l},\tilde{Y}_{l},l\in[\varphi],\varphi=r)$ 
 and follows the decomposition of $B_{0}$ in \eqref{postdisint} with given quantities as indicated, the distribution of $B_{0,J}\sim \mathrm{CRM}(\tau_{0,J},F_{0})$ unchanged, and $H_{l}|\tilde{X}_{l}=x_{l}$ having density proportional to $\lambda^{x_{l}}e^{-\lambda\sum_{j=1}^{J}\Phi_{j}(M_{j})}\tau_{0}(\lambda)$, conditionally independent for $l\in[r].$
\end{prop}

\section{Exact sampling of the HIBP marginal process}
We now use the results above to provide a description of the marginal process for the generalized HIBP, practically implementable in many situations by sampling the basic variables above which appear in the non-hierarchical setting. 
\begin{thm}\label{jointmarginal}
With reference to variables in Propositions~\ref{IBPpost} and \ref{prop2}, 
the marginal distribution of the HIBP $(Z^{(i)}_{j}, i\in[M_{j}], j\in[J])$ is equal to the distribution of the process:
    \begin{equation}
    \label{hibpmarginal}
    \left(\sum_{l=1}^{\varphi}\left[\sum_{k=1}^{X_{j,l}}\tilde{A}^{(i)}_{j,k,l}\right]\delta_{\tilde{Y}_{l}}, i\in [M_{j}], j\in [J]\right)
    \end{equation}
\end{thm}
Note despite the relative simplicity of sampling according to~\eqref{hibpmarginal}, these lead to exact sampling of quite complex multi-dimensional distributions. Furthermore for comparison with existing results, the random count matrices in~\cite{ZhouPadilla2016} corresponding to Poisson-Gamma-Gamma HIBP with $(M_{j}=1, j\in[J]),$ produce matrices with entries in $\mathbb{Z}^{J\times \varphi}$ where  $\mathbb{Z}=\{0,1,2,\ldots\},$ $J$ is the number of rows and $\varphi$ columns corresponding to observed features $(\tilde{Y}_{l},l\in[\varphi]).$ In our general mixed Poisson HIBP setting random count matrices take values in $\mathbb{Z}^{\sum_{j=1}^{J}M_{j}\times \varphi},$ changing with $(M_{j},j\in[J])$ for much more general choices of $\mu_{j},B_{0}$. See section~\ref{marginal} for precise descriptions of their apparently complex joint distributions. We show in the next section more details for the Poisson cases and also establish that all general HIBP can be sampled as processes coupled to the Poisson HIBP framework. 
\subsection{Exact sampling of every HIBP marginal distribution via the Poisson case}
While we can use Theorem \ref{jointmarginal} for $A_{j}$ directly, we show that one can sample from the marginal distribution of every HIBP by coupling it with a Poisson HIBP through the transformation $\mathbb{P}(A_{j}\neq 0|v)=\pi_{j}(v)=(1-e^{-s})$, where $\pi^{-1}_{j}$ denotes its inverse. Note we may couple many such process for each $j$ simultaneously by using different distributions for $A_{j}$.

We first provide some general important details for sampling the Poisson HIBP, as a special case of Theorem~\ref{jointmarginal}, which we shall use throughout the remainder of this work. Hereafter, to distinguish the Poisson cases, we will use $C_{j}:=A_{j}\sim \mathrm{Poisson}(s)$, setting $\tilde{C}^{(i)}_{j,k,l}=\tilde{A}^{(i)}_{j,k,l}$. Furthermore, we use $B_{j}$ in place of $\mu_{j},$ and set $\rho_{j}=\tau_{j},$ with $\tau^{(M_{j})}_{j}(s)=e^{-sM_{j}}\tau_{j}(s)$ and we set $\Phi_{j}(M_{j})=\psi_{j}(M_{j})=\int_{0}^{\infty}(1-e^{-sM_{j}})\tau_{j}(s)ds$. Here, $f_{S_{j}}(s)=(1-\mathrm{e}^{-M_{j}s})\tau_{j}(s)/\psi_{j}(M_{j})$ which specifies the distribution of the corresponding distribution of jumps $S_{j,k,l}$, where $(\tilde{C}^{(i)}_{j,k,l}, i\in[M_{j}])|S_{j,k,l}=s\sim\mathrm{tP}(s(M_{j},\ldots,M_{j}))$. Further details can be read from~\cite{James2017}[Section 4.2]. 

In particular, important for our general exposition,  for $\tilde{C}_{j,k,l}:=\sum_{i=1}^{M_{j}}\tilde{C}^{(i)}_{j,k,l},$ the joint distribution of $S_{j,k,l}, \tilde{C}_{j,k,l}$ may be expressed as 
\begin{equation}
\label{MtPsimple}
\frac{s^{c_{j,k,l}}e^{-sM_{j}}\rho_{j}(s)}{\psi^{(c_{j,k,l})}_{j}(M_{j})}\times 
\frac{M_{j}^{c_{j,k,l}}\psi^{(c_{j,k,l})}_{j}(M_{j})}{\psi_{j}(M_j)c_{j,k,l}!}
\end{equation}
where $\psi^{(c_{j,k,l})}_{j}(M_{j})=\int_{0}^{\infty}s^{c_{j,k,l}}e^{-sM_{j}}\rho_{j}(s)ds$ with interpretations of~\eqref{MtPsimple}, in these Poisson settings, given precisely in~\cite{Pit97}. Additionally, also important for our general exposition, work of \cite{Pit97}[Proposition 12, Corollary 13], as noted in ~\cite{James2017}[Remark 4.1], also gives the variables $(H_{l},\tilde{X}_{l})$ more precise meaning in the Poisson setting, in terms of decomposing the space. From that work the joint distributions $(H_{l},\tilde{X}_{l})$ can be expressed in terms of $H_{l}|\tilde{X}_{l}=x_{l}$ and $\tilde{X}_{l},$ as,  
\begin{equation}
\label{HXdecomp}
\frac{\lambda^{x_{l}}e^{-\lambda\sum_{j=1}^{J}\psi_{j}(M_{j})}\tau_{0}(\lambda)}{\Psi^{(x_{l})}_{0}(\sum_{j=1}^{J}\psi_{j}(M_{j}))}\times \frac{{(\sum_{j=1}^{J}\psi_{j}(M_{j}))}^{x_{l}}\Psi^{(x_{l})}_{0}(\sum_{j=1}^{J}\psi_{j}(M_{j}))}
{x_{l}!\Psi_{0}(\sum_{j=1}^{J}\psi_{j}(M_{j}))}
\end{equation}
for $\Psi^{(x_{l})}_{0}(\sum_{j=1}^{J}\psi_{j}(M_{j}))=\int_{0}^{\infty}
\lambda^{x_{l}}e^{-\sum_{j=1}^{J}\psi_{j}(M_{j})}\tau_{0}(\lambda)d\lambda.$
That is $\tilde{X}_{l}\sim \mathrm{MtP}(\sum_{j=1}^{J}\psi_{j}(M_{j}),\tau_{0})$

\begin{prop}\label{sumsampleprop}
Let for each $j$, $(Z^{(i)}_{j}, i\in [M_{j}]) | B_{j} \overset{iid}\sim \mathrm{PoiP}(B_{j})$, where $B_{j} | B_{0} \sim \mathrm{CRM}(\tau_{j}, B_{0})$ conditionally independent across $j\in[J]$. Then consider the sum process $(\sum_{i=1}^{M_{j}} Z^{(i)}_{j} , j\in[J])$, where $\sum_{i=1}^{M_{j}} Z^{(i)}_{j} | \mu_{j} \sim \mathrm{PoiP}(M_{j}B_{j})$.
    \begin{enumerate}
    \item The joint marginal distribution of the sum process $(\sum_{i=1}^{M_{j}} Z^{(i)}_{j}, i\in[M_{j}],j\in [J])$ is equal in distribution to
    \begin{equation}\label{prop:sumequation}
(\sum_{l=1}^{\varphi} N_{j,l}\delta_{\tilde{Y}_{l}}, j\in [J])\overset{d}=((\sum_{l=1}^{\varphi} \sum_{k=1}^{X_{j,l}}\tilde{C}_{j,k,l}\delta_{\tilde{Y}_{l}}, j\in [J])
    \end{equation}
\item where $(\tilde{C}_{j,k,l}=\sum_{i=1}^{M_{j}} \tilde{C}^{(i)}_{j,k,l}, k\in [X_{j,l}]) \overset{iid}\sim \mathrm{MtP}(M_{j}, \tau_{j})$ with the $X_{j,l}$ components independent across $l$.
\item $\tilde{X}_{l}\sim \mathrm{MtP}(\sum_{j=1}^{J}\psi_{j}(M_{j}),\tau_{0})$
\item $(X_{j,l},j\in[J])|\tilde{X}_{l}=x_{l}\sim \mathrm{Multi}(x_{l}; q_{1},\ldots, q_{J})$, for $q_{j}=\frac{\psi_{j}(M_{j})}{\sum_{v=1}^{J}\psi_{v}(M_{v})}.$    \item  $(\sum_{k=1}^{X_{j,l}} \tilde{C}^{(i)}_{j,k,l}, i\in [M_{j}]) |N_{j,l}:=\sum_{k=1}^{X_{j,l}} \tilde{C}_{j,k,l} = n_{j,l}$ is $\mathrm{Multi}(n_{j,l}; 1/M_{j}, \ldots, 1/M_{j})$, for each $j,l$,
    \end{enumerate}
\end{prop}

Proposition~\ref{sumsampleprop} indicates that one can sample from the marginal distribution of $(Z^{(i)}_{j},i\in[M_{j}], j\in[J]),$ by first sampling $((N_{j,l}, j\in [J]),\tilde{Y}_{l},l\in[\varphi],\varphi)$ from the sum process as in~\eqref{prop:sumequation} and then applying Multinomial sampling as in item 5. The next results shows how to couple and sample any general HIBP with a corresponding Poisson HIBP based on transformation of the jumps.

\begin{prop}\label{coupled}
Consider the Poisson HIBP setting with variables as defined in Proposition~\ref{sumsampleprop}.
For $A_{j}|v\sim \mathbb{P}_{A_{j}}(da|v)$ with 
$\mathbb{P}(A_{j}\neq 0|v)=\pi_{j}(v)$, specify a L\'evy density $\rho_{j}(v)$, such that $\psi_{j}(M_{j})=\int_{0}^{\infty}[(1-(1-\pi_{j}(v))^{M_{j}}]\rho_{j}(v)dv=\int_{0}^{\infty}(1-\mathrm{e}^{-sM_{j}})\tau_{j}(s)ds$.
Then consider variables $(\hat{A}^{(i)}_{j,k,l}, i\in[M_{j}])|S_{j,k,l}=s$ with iid zero-truncated distributions proportional to $\mathbb{I}_{a\neq 0}\mathbb{P}_{A_{j}}(da|\pi^{-1}_{j}(1-\mathrm{e}^{-s}))$, for $S_{j,k,l}$ iid with common density $f_{S_{j}}(s)=(1-\mathrm{e}^{-M_{j}s})\tau_{j}(s)/\psi_{j}(M_{j})$
\begin{enumerate}
\item Then the following equates to the marginal distribution of a coupled Poisson HIBP,
$$
\left(\sum_{l=1}^{\varphi}\left[\sum_{k=1}^{X_{j,l}}\mathbb{I}_{\{\tilde{C}^{(i)}_{j,k,l}>0\}}\hat{A}^{(i)}_{j,k,l}\right]\delta_{\tilde{Y}_{l}}, \sum_{l=1}^{\varphi}\left[\sum_{k=1}^{X_{j,l}}\tilde{C}^{(i)}_{j,k,l}\right]\delta_{\tilde{Y}_{l}}, i\in[M_{j}], j\in [J]\right)
$$
\item Where $(\sum_{l=1}^{\varphi}\left[\sum_{k=1}^{X_{j,l}}\mathbb{I}_{\{\tilde{C}^{(i)}_{j,k,l}>0\}}\hat{A}^{(i)}_{j,k,l}\right]\delta_{\tilde{Y}_{l}}, i\in[M_{j}], j\in [J])$ has the marginal distribution of a general HIBP with $A_{j},
\rho_{j},\tau_{0}$ as specified in \eqref{hibpmarginal} 
\item $(\sum_{l=1}^{\varphi}\left[\sum_{k=1}^{X_{j,l}}\mathbb{I}_{\{\tilde{C}^{(i)}_{j,k,l}>0\}}\right]\delta_{\tilde{Y}_{l}}, i\in[M_{j}], j\in [J])$ corresponds to a Bernoulli HIBP, where $Z^{(i)}_{j}|\mu_{j}\sim \mathrm{BerP}(\mu_{j})$ and $\mu_{j}|B_{0}\sim \mathrm{CRM}(\rho_{j},B_{0})$.
\end{enumerate}
\end{prop}

\section{Posterior distributions for $B_{0}$ and $(\mu_{j}, j\in [J])$}
While our results show that we can sample exactly the processes $((Z^{(i)}_{j}, i\in[M_{j}]), j\in[J])$ via the variables
$(X_{j,l},({A}^{(i)}_{j,k,l},k\in [X_{j,l}],\tilde{Y}_{l}, l\in[r],\varphi=r, i\in [M_{j}], j\in J),$ if we assume we are only given the information in the end result, then this information is given as 
$(\sum_{k=1}^{X_{j,l}}\tilde{A}^{(i)}_{j,k,l},\tilde{Y}_{l}, l\in[r],\varphi=r, i\in [M_{j}], j\in J).$ We now provide a description of the posterior distributions of $B_{0}$ and $(\mu_{j},j\in [J]).$ 

\begin{thm}\label{postBgivenZq}
A description of the posterior distributions of $(\mu_{j}\in [J],B_{0})|\mathbf{Z}_{J}$ for 
$\mathbf{Z}_{J}:=(\sum_{k=1}^{X_{j,l}}\tilde{A}^{(i)}_{j,k,l},\tilde{Y}_{l}, l\in[r],\varphi=r, i\in [M_{j}], j\in J)$ is as follows,
\begin{enumerate}
\item The joint distribution of $(\mu_{j},j\in[J]),B_{0}|\mathbf{Z}_{J}$ is such that component-wise and jointly, is equal in distribution to,
  \begin{equation}
  \label{postsumrep}
  (\tilde{\mu}_{j,M_{j}}+\sum_{l=1}^{r}\left[\tilde{\sigma}_{j,l}(H_{l})+\sum_{k=1}^{X_{j,l}}S_{j,k,l}\right]\delta_{\tilde{Y}_{l}}, j\in[J]),B_{0,J}+\sum_{l=1}^{r}H_{l}\delta_{\tilde{Y}_{l}} 
  \end{equation}
for $B_{0,J}\sim \mathrm{CRM}(\tau_{0,J},F_{0}),$ $\tilde{\mu}_{j,M_{j}}\sim\mathrm{CRM}(\rho^{(M_{j})}_{j},B_{0,J})$, and $(\tilde{\sigma}_{j,l})\overset{iid}\sim\mathrm{CRM}(\rho^{(M_{j})}_{j},\mathbb{L}),$ where $\mathbb{L}(dt)=dt$ on $(0,\infty),$  independent of $(H_{l})$ and $(S_{j,k,l}).$

\item The conditional distribution of $S_{j,k,l}|(\tilde{A}^{(1)}_{j,k,l},\ldots,\tilde{A}^{(M_{j})}_{j,k,l})$ is determined by $\eqref{joint1}$ and $H_{l}|\tilde{X}_{l}=x_{l}$ has density $\propto \lambda^{x_{i}}e^{-\lambda\sum_{j=1}^{J}\Phi_j(M_j)}\tau_{0}(\lambda)$,
\item \label{latentlaws} The posterior distribution is concluded by obtaining a 
description of $(X_{j,l},(\tilde{A}^{(i)}_{j,k,l},k\in [X_{j,l}], l\in[r], i\in [M_{j}], j\in J),$ given  $\mathbf{Z}_{J}$
\end{enumerate}
\end{thm}

Theorem~\ref{postBgivenZq} relies on inferring the conditional distribution of the latent variables $(X_{j,l},({A}^{(i)}_{j,k,l},k\in [X_{j,l}]), l\in[r], i\in [M_{j}], j\in J),$ based on the observed sums $(\sum_{k=1}^{X_{j,l}}\tilde{A}^{(i)}_{j,k,l},\tilde{Y}_{l}, l\in[r],\varphi=r, i\in [M_{j}], j\in J).$ While there is potential for developing computational methods to leverage this concept, the general approach remains unclear. Moreover, employing the results in \cite{James2017}[section 5] expressing the processes as multivariate IBP offers a means to describe the  posterior distributions solely based on the sums without necessitating knowledge of the conditional distribution of the individual latent variables. Nevertheless, this strategy is practically viable only  when ample information about the sum distribution is available. 

The Bernoulli HIBP case with $(M_{j}=1, j\in[J]),$ corresponds in distribution to the processes $(\sum_{l=1}^{\infty}\xi_{j,l}\delta_{Y_{l}}, j\in[J]),$ for $\xi_{j,l}\sim \mathrm{Poisson}(\Phi_{j}(1)\lambda_{l})$ suggesting a model where one is interested in only increasing $J$ to $J+1$ etc.. This, in fact corresponds to the bipartite graphs model in ~\citep{Caron2012}, with $J$ readers, $(Y_{l})$ books and $(\lambda_{l})$ their respective popularity parameters, set further $\gamma_{j}=\Phi_{j}(1).$ Theorem~\ref{postBgivenZq} gives a complete tractable description, which is otherwise unavailable and holds for arbitrary choice of $\gamma_{j}$ in~\citep{Caron2012}. In particular, it follows that $\tilde{A}^{(1)}_{j,k,l}=1$ and hence, in this case,  $X_{j,l}=\sum_{k=1}^{X_{j,l}}\tilde{A}^{(1)}_{j,k,l}$ is observed and $S_{j,k,l}$ has density proportional to $s \rho_{j}(s),$ for $s$ in $(0,1),$ and $j=1,2,\ldots.$ 
The Bernoulli case for $M_{j}>1$ is more complex. In contrast, as we shall show next a complete and tractable description of the posterior distribution in the Poisson HIBP case can be accomplished in a general fashion, based on novel usage of existing results in ~\cite{JLP2,Kolchin,Pit97,Pit02,FoFZhou}. This in turn, due to Proposition~\ref{coupled}, through coupling, can be used to provide more information for other general HIBP cases. 

\section{Posterior distribution for Poisson HIBP}
We now address the remarkable case of the Poisson HIBP. 
Where the given structure is  $\mathbf{Z}_{J}:=(\sum_{k=1}^{X_{j,l}}\tilde{C}^{(i)}_{j,k,l},\tilde{Y}_{l}, l\in[r],\varphi=r, i\in [M_{j}], j\in J).$ However in view of Proposition~\ref{sumsampleprop}, it suffices to work with the information in the sum process $(N_{j,l},\tilde{Y}_{l}, l\in[r],\varphi=r, j\in J),$ where again $N_{j,l}=\sum_{k=1}^{X_{j,l}}\tilde{C}_{j,k,l},$ for descriptions of the posterior distribution. The task then becomes to describe distributions of $(\tilde{C}_{j,k,l}, k\in [X_{j,l}], X_{j,l})$ given this information. See section~\ref{recovery} for distributional recovery of individual counts $\tilde{C}^{(i)}_{j,k,l}.$

One of our innovations here is to first note that conditioning on $(N_{j,l}, j\in[J])$, and $H_{l}=\lambda$ we can apply results in \cite{Kolchin,Pit97}, and also \cite{JLP2,FoFZhou}, to directly deduce the following key results based on finite Gibbs $\mathrm{EPPF},$ with these $\mathrm{EPPF}$ (Exchangeable Partition Probability Function) otherwise appearing  specifically in \cite{JLP2}{Proposition 4] and \cite{FoFZhou}[Sections 2.3,3]. See~\cite{Pit96,Pit02} for more on general EPPF's for species sampling models.

\begin{prop}\label{prop:GibbsEPPF}
Consider the Poisson HIBP setting in Proposition~\ref{sumsampleprop}. Then for each $j\in [J],$ $(\tilde{C}_{j,k,l}, k\in [X_{j,l}], X_{j,l})$ given $N_{j,l}=n_{j,l}$ and $H_{l}=\lambda$,  follows the distribution of a random partition of the integers $[n_{j,l}].$  Such that $\sum_{k=1}^{X_{j,l}}\tilde{C}_{j,k,l}=n_{j,l},$ where $\tilde{C}_{j,k,l}=c_{j,k,l}$ dictates the size of each block and $X_{j,l}=x_{j,l}$ corresponds to the random number of blocks, according to the finite Gibbs EPPF arising in~\cite{JLP2,FoFZhou} as follows:

\begin{enumerate}
\item The finite Gibbs partition distribution specifying these conditional distributions is, for
$\mathbf{c}_{j,l}=(c_{j,k,l},k\in[x_{j,l}]),$ and  $\psi^{(c)}_{j}(M_{j})=\int_0^{\infty} s^{c}e^{-sM_j}\tau_{j}(s)\,ds,$
:
\begin{equation}
\label{GibbsEPPF}
p^{[n_{j,l}]}(\mathbf{c}_{j,l}|\lambda\tau_{j}, M_j) = \frac{\lambda^{x_{j,l}}\prod_{k=1}^{x_{j,l}}\psi^{(c_{j,k,l})}_{j}(M_{j})}{\Xi^{[n_{j,l}]}_j(\lambda,M_j)}
\end{equation}
\item where, there is the identity:
$
\int_0^{\infty} t^{n_{j,l}}e^{-tM_j}\eta_j(t|\lambda)\,dt = e^{-\lambda\psi_{j}(M_{j})}\Xi^{[n_{j,l}]}_j(\lambda,M_{j})
$
for $\sum_{*}$ denoting sum over all positive integers $\mathbf{c}_{j,l}$ such that $\sum_{k=1}^{x_{j,l}}c_{j,k,l}=n_{j,l}$:
\begin{equation}
\Xi^{[n_{j,l}]}_j(\lambda,M_j)=\sum_{x_{j,l}=1}^{n_{j,l}}\frac{n_{j,l}!\lambda^{x_{j,l}}}{x_{j,l}!}\sum_{*}\frac{\prod_{k=1}^{x_{j,l}}\psi^{(c_{j,k,l})}_{j}(M_{j})}{\prod_{k=1}^{x_{j,l}} c_{j,k,l}!}
\end{equation}
\item The distribution of the number of blocks $X_{j,l}=x_{j,l}$ conditioned on $(N_{j,l}=n_{j,l}, j\in[J], H_{l}=\lambda$ is, $X_{j,l}=0$ when $n_{j,l}=0$, and otherwise:
\begin{equation}
\label{Kn}
p^{[n_{j,l}]}(x_{j,l}|\lambda\tau_{j},M_j) = \frac{n_{j,l}!\lambda^{x_{j,l}}}{x_{j,l}!}\sum_{*}\frac{\prod_{k=1}^{x_{j,l}}\psi^{(c_{j,k,l})}_{j}(M_{j})}{\Xi^{[n_{j,l}]}_j(\lambda,M_j)\prod_{k=1}^{x_{j,l}} c_{j,k,l}!}
\end{equation}
\item As in \cite{Kolchin, Pit97}, \eqref{Kn} and~\eqref{MtPsimple}  yields the identity, $p^{[n_{j,l}]}(x_{j,l}|\lambda\tau_{j},M_j)=$
\begin{equation}
\label{sumid}
\frac{n_{j,l}!\lambda^{x_{j,l}}({\psi_{j}(M_{j})})^{x_{j,l}}}{x_{j,l}!M^{n_{j,l}}_{j}\Xi^{[n_{j,l}]}_j(\lambda,M_j)}\mathbb{P}(\tilde{C}_{j,1,l}+\ldots+\tilde{C}_{j,x_{j,l},l}=n_{j,l})
\end{equation}
\end{enumerate}
\end{prop}

\begin{exam}[$B_{j}|B_{0}\sim \mathrm{G}(\theta_{j},\zeta;B_{0})$]
Now as important examples, see \cite{JLP2}[p. 86], when $\tau_{j}(s)=\theta_{j}s^{-1}e^{-s\zeta_{j}}$ corresponding to a gamma process,  with $\psi_{j}(M_{j})=\theta_{j}[\log(\zeta_{j}+M_{j})-\log(\zeta_{j})],$ $\psi^{(c)}_{j}(M_{j})=\theta_{j}\Gamma(c){(\zeta_{j}+M_{j})}^{-c},$
and, gamma density, $\eta_{j}(t_{j}|\lambda)=t^{\theta_{j}\lambda_{j}-1}_{j}{\mbox e}^{-t_{j}\zeta_{j}}\zeta^{\theta_{j}\lambda}/\Gamma(\theta_{j}\lambda)$ leading to
$\Xi^{[n_{j,l}]}_j(\lambda,M_{j})=\frac{\Gamma(\theta_{j}\lambda+n_{j,l})}{{(\zeta_{j}+M_{j})}^{n_{j,l}}\Gamma(\theta_{j}\lambda)}.$ It follows that,  
 $p^{[n_{j,l}]}(\mathbf{c}_{j,l}|\lambda\tau_{j}, M_j),$ equates with the well-known EPPF of a Dirichlet process with parameter $\theta_{j}\lambda:$
$p_{0,\theta_{j}\lambda}(\mathbf{c}_{j,l}) := \frac{\Gamma(\theta_{j}\lambda)\theta_{j}^{x_{j,l}} \lambda^{x_{j,l}}}{\Gamma(\theta_{j}\lambda+n_{j,l})}  \prod_{k=1}^{x_{j,l}} \Gamma(c_{j,k,l}),$ and does not depend on $(\zeta_{j},M_{j})$.
Now using this, it follows that in this case, $p^{[n_{j,l}]}(x_{j,l}|\lambda\tau_{j},M_j)$ is the same as, 
\begin{equation}
\mathbb{P}^{(n_{j,l})}_{0,\theta_{j} \lambda}(x_{j,l})=\frac{\Gamma(\theta_{j}\lambda)\theta_{j}^{x_{j,l}} \lambda^{x_{j,l}}}{\Gamma(\theta_{j}\lambda+n_{j,l})}\stirlingfirst{n_{j,l}}{x_{j,l}}
\label{DPKn}
\end{equation}
where on the left hand side is the distribution of the number of blocks of a partition of $[n_{j,l}]$ for a Dirichlet process, where $\stirlingfirst{n}{k}$ denotes the unsigned Stirling number of the first kind,as in~\cite{CraneEwens,TavareEwens}.
From~ \eqref{sumid} and \eqref{DPKn}, we arrive at the equality
\begin{equation}
\mathbb{P}(\tilde{C}_{j,1,l}+\ldots+\tilde{C}_{j,x_{j,l},l}=n_{j,l})=
\frac{x_{j,l}!}{n_{j,l}!}\frac{{(1-\frac{\
\zeta_{j}}{\zeta_{j}+M_{j}})}^{n_{j,l}}}{(\log(\zeta_{j
}+M_{j})-\log(\zeta_{j}))^{x_{j,l}}}\stirlingfirst{n_{j,l}}{x_{j,l}}
\end{equation}
indicating the probability mass function of the sum-log distribution of $\mathbb{P}(\tilde{C}_{j,1,l}+\ldots+\tilde{C}_{j,x_{j,l},l}=n_{j,l}),$
with $\mathbb{P}(\tilde{C}_{j,k,l}=c)=
\frac{{(\frac{M}{\zeta+M})}^{c}}
{c[\log(\zeta+M)-\log(\zeta)]},
$
which informs the strategy in ~\cite{ZhouCarin2015,ZhouPadilla2016}. 
\begin{rem}
The case of $\tau_{j}(s)=Cs^{\alpha-1},$ corresponding to $0<\alpha<1$ stable processes is treated in~\cite{Pit97}, see also~\cite{Pit06}[p. 73], and will be illustrated in the form of generalized gamma processes later. 
\end{rem}
\end{exam}

We now present an implementable description of the posterior distribution for Poisson HIBP using $p^{[n_{j,l}]}(\mathbf{c}_{j,l}|\tau_{j}, M_j)\Xi^{[n_{j,l}]}_j(1,M_j)=\prod_{k=1}^{x_{j,l}}\psi^{(c_{j,k,l})}_{j}(M_{j}).$

\begin{thm}\label{postPoissonHIBP}
A description of the posterior distributions of $(B_{j}\in [J],B_{0})|\mathbf{Z}_{J}$ is as follows,
\begin{enumerate}
\item The joint distribution of $(B_{j},j\in[J]),B_{0}|\mathbf{Z}_{J}$ is such that component-wise and jointly, is equivalent in distribution to,
  \begin{equation}
  \label{postsumrep}
  (B_{j,M_{j}}+\sum_{l=1}^{r}\left[\tilde{\sigma}_{j,l}(H_{l})+\sum_{k=1}^{X_{j,l}}S_{j,k,l}\right]\delta_{\tilde{Y}_{l}}, j\in[J]),B_{0,J}+\sum_{l=1}^{r}H_{l}\delta_{\tilde{Y}_{l}} 
  \end{equation}
for $B_{0,J}\sim \mathrm{CRM}(\tau_{0,J},F_{0}),$ $B_{j,M_{j}}\sim\mathrm{CRM}(\tau^{(M_{j})}_{j},B_{0,J})$, and 
\item $(\tilde{\sigma}_{j,l}(H_{l}),j\in[J])|H_{l}=\lambda$ are independent with density $\eta^{[0]}_{j}(t_{j}|\lambda,M_{j})=e^{-t_{j}M_{j}}\eta_{j}(t_{j}|\lambda)e^{\psi_{j}(M_{j})}$, independent of $(S_{j,k,l}).$ 
\item $S_{j,k,l}|\tilde{C}_{j,k,l}=c_{j,k,l}$ has density $\frac{s^{c_{j,k,l}}e^{-sM_{j}}\tau_{j}(s)}{\psi^{(c_{j,k,l})}(M_{j})}$
\item  $H_{l}|\tilde{X}_{l}=x_{l}$ has density $\frac{\lambda^{x_{l}}e^{-\lambda\sum_{j=1}^{J}\psi_{j}(M_{j})}\tau_{0}(\lambda)}{\Psi^{(x_{l})}_{0}(\sum_{j=1}^{J}\psi_{j}(M_{j}))}$, where $\tilde{X}_{l}=\sum_{j=1}^{J}X_{j,l}$
\item \label{latentlaws} $(X_{j,l},(\tilde{C}_{j,k,l},k\in [X_{j,l}], l\in[r], j\in J),$ given $(N_{j,l}, j\in[J]),l\in[r])$ are conditionally independent over $l\in[r],$ with distributions proportional to 
\begin{equation}
\label{jointrep1}
\frac{\Psi^{(x_{l})}_{0}(\sum_{j=1}^{J}\psi_{j}(M_{j}))}
{\Psi_{0}(\sum_{j=1}^{J}\psi_{j}(M_{j}))}\prod_{j=1}^{J}p^{[n_{j,l}]}(\mathbf{c}_{j,l}|\tau_{j}, M_j)\Xi^{[n_{j,l}]}_j(1,M_j)
\end{equation}
where $x_{l}=\sum_{j=1}^{J}x_{j,l},$ and noting $X_{j,l}=0$ if $N_{j,l}=0.$
\item In particular $\tilde{X}_{l}|(N_{j,l}=n_{j,l},j\in[J]),l\in[r]).$ for $x_{l}=\sum_{j=1}^{J}x_{j,l},$ with $x_{j,l}=0$ when $n_{j,l}=0,$ is determined by the joint distribution of $(X_{j,l}=x_{j,l}, j\in[J])$ proportional to 
$$
\frac{\Psi^{(x_{l})}_{0}(\sum_{j=1}^{J}\psi_{j}(M_{j}))}
{\Psi_{0}(\sum_{j=1}^{J}\psi_{j}(M_{j}))}\prod_{j=1}^{J}p^{[n_{j,l}]}(x_{j,l}|\tau_{j}, M_j)\Xi^{[n_{j,l}]}_j(1,M_j)
$$
\end{enumerate}
\end{thm}
The following proposition shows the posterior distribution in Theorem~\ref{postPoissonHIBP} only depends directly on $(X_{j,l}, j\in[J]).$ In other words knowledge of the sum counts $\tilde{C}_{j,k,l}$ is not necessary, although for computational purposes they may be incorporated. This result follows from the randomization of the posterior distribution of the total mass in \cite{JLP2}[Theorem 1], or equivalently \cite{James2002}[Corollary 5.1]. See also \cite{Pit02}. However see section~\ref{recovery} for interest in obtaining the counts $\tilde{C}_{j,k,l}.$
\begin{prop}\label{poissonequivalence}
Consider the Poisson HIBP setting in Theorem~\ref{postPoissonHIBP}. Then given $N_{j,l}=n_{j,l},H_{l}=\lambda,$ there is the distributional equivalence 
\begin{equation}
\hat{\sigma}_{j,l}(\lambda)\overset{d}=\tilde{\sigma}_{j,l}(\lambda)+\sum_{k=1}^{X_{j,l}}S_{j,k,l}
\end{equation}
where $\hat{\sigma}_{j,l}(\lambda)$ has density $\eta^{[n_{j,l}]}(t_{j}|\lambda,M_{j})=\frac{t_{j}^{n_{j,l}}\eta^{[0]}(t_{j}|\lambda,M_{j})}{\Xi^{[n_{j,l}]}_j(\lambda,M_{j})}$. In the case where $B_{j}\sim \mathrm{G}(\theta_{j},\zeta_{j};B_{0})$ it follows that $\hat{\sigma}_{j,l}(H_{l})\sim \mathrm{Gamma}(\theta H_{l}+n_{j,l},M_{j}+\zeta_{j}).$
\end{prop}

We have not discussed results for specific cases of $B_{0}$ so far. Here, we describe results for the case where $B_{0}\sim \mathrm{GG}(\alpha_{0},\theta_{0},\zeta_{0},F_{0})$, denoting a generalized gamma process as specified below, we will also look at the case where $B_{j}\sim \mathrm{GG}(\theta_{j},\zeta_{j},B_{0})$. Prior to that, for $0<\alpha_{j}<1,$ let $P_{\alpha_{j},0} \sim \mathscr{PY}(\alpha_{j},0,F_{0})$ denote the $(\alpha_{j},0)$ Pitman-Yor process, \cite{IJ2001,Pit96}, constructed from normalizing a $\alpha_{j}$ stable-subordinator. The corresponding EPPF  describing the distribution of a random partition of $[n_{j,l}]$ is then given by:
$$
p_{\alpha_{j},0}(\mathbf{c}_{j,l}) = \frac{\alpha^{x_{j,l}-1}_{j}\Gamma(x_{j,l})}{\Gamma(n_{j,l})}\prod_{k=1}^{x_{j,l}}\frac{\Gamma(c_{j,k,l}-\alpha_{j})}{\Gamma(1-\alpha_{j})}.
$$
From this, the probability of the number of blocks $\mathbb{P}_{\alpha_{j},0}^{(n_{j,l})}(k) =
\frac{{\alpha^{k-1}_{j}\Gamma(k)} S_{\alpha_{j}}(n_{j,l},k)}{ {\Gamma(n_{j,l})}}$, where $S_{\alpha_{j}}(n_{j,l},k)$ denotes the generalized Stirling number of the second kind, as in \cite{Pit06}[3.2.3, p. 65-66].

\begin{prop}\label{PropgenPGG}
Set $B_{0}\sim \mathrm{GG}(\alpha_{0},\theta_{0},\zeta_{0}, F_{0})$ denoting a generalized gamma process $\mathrm{CRM}$ with L\'evy density for $0<\alpha_{0}<1$:
$
\tau_{0}(\lambda)=\theta_{0}\alpha_{0}\lambda^{-\alpha_{0}-1}e^{-\zeta_{0}\lambda}/\Gamma(1-\alpha_{0})$ and $\Psi_{0}(t)=\theta_{0}[(\zeta_{0}+t)^{\alpha_{0}}-\zeta^{\alpha_{0}}_{0}].$ Then given $\mathbf{Z}_{J}$:
\begin{enumerate}
\item[1.] $B_{0,J}\sim \mathrm{GG}(\alpha_{0},\theta_{0},\zeta_{0}+\sum_{j=1}^{J}\psi_{j}(M_{j}))$
\item[2.] $H_{l}|\tilde{X}_{l}\sim \mathrm{Gamma}(\tilde{X}_{l}-\alpha_{0},\zeta_{0}+\sum_{j=1}^{J}\psi_{j}(M_{j}))$ and 
\item[3.] for $\tau_{j}(s) =\frac{ \alpha_{j}\theta_{j}e^{-s\zeta_{j}}s^{-\alpha_{j}-1}}{\Gamma(1-\alpha_{j})}$, for $j\in[J],$ $((\tilde{C}_{j,k,l}=c_{j,k,l},k\in [X_{j,l}],X_{j,l}=x_{j,l}) j\in J),$ given $(N_{j,l}, j\in[J])$
has, for $x_{l}=\sum_{j=1}^{J}x_{j,l},$ with $x_{j,l}=0$ when $n_{j,l}=0,$ and otherwise $x_{j,l}\in[n_{j,l}],$ distribution proportional to,
\begin{equation}
\frac{\Gamma(x_{l} - \alpha_{0})}{\prod_{j=1}^{J}\Gamma(x_{j,l})}\frac{\prod_{j=1}^{J}\alpha_{j}\theta^{x_{j,l}}_{j}{(\zeta_{j}+M_{j})}^{x_{j,l}\alpha_{j}}}{(\zeta_{0} + \sum_{j=1}^{J} \theta_{j} [(\zeta_{j} + M_j)^{\alpha_{j}} - \zeta_{j}^{\alpha_{j}}])^{x_{l} - \alpha_{0}}} \prod_{j=1}^{J}p_{\alpha_{j},0}(\mathbf{c}_{j,l}).
\label{specialstable}
\end{equation}
\end{enumerate}

\end{prop}
\subsection{A simplification for any $J,$ with Poisson-GG-GG examples}
\label{simple}
We have developed quite general results which allow for a great deal of flexibility. However, with respect to this particular exposition, we wish to provide simpler and more readable descriptions and implementations of models which are still sufficiently rich for practical use and have remarkable distributional properties. This in no way precludes the usage of the more complex generalized models, which allow for fitting of additional parameters, and mixed matching of processes.
The ability to sample over all groups $J$ simultaneously, via Proposition~\ref{sumsampleprop} and the prediction rule described in the forthcoming Proposition~\ref{Prediction}, makes the assumption of $M_{j}=M$ for all $j\in[J]$, where $M = \max_{j\in[J]} M_j$, a plausible simplification. Under this assumption, groups with customers less than $M$ can be taken to have $M - M_j$ additional "ghost" customers who selected no features, corresponding to submission of blank documents perhaps. Furthermore, we set $\tau_{j}(s)=\theta_{j}\tau(s)$ for all $j\in [J],$ so that $\psi_{j}(M)=\theta_{j}\psi(M)$ and hence $\sum_{j=1}^{J}\psi_{j}(M_{j})=\psi(M)\sum_{j=1}^{J}\theta_{j}$ 
and $q_{j}=\theta_{j}/\sum_{v=1}^{J}\theta_{v},$ we set $\mathbf{q}:=(q_{1},\ldots , q_{J}).$

\subsection{Poisson-GG-GG}\label{PGGGG}
For the simplified Poisson-Generalized Gamma-Generalized Gamma model, that is 
Poisson-GG-GG,  set in~\eqref{specialstable}, $\alpha_{0}=\beta/\alpha$ for $0<\beta<\alpha,$ $\tau_{j}(s)=\theta_{j}\frac{\alpha s^{-\alpha-1}e^{-\zeta s}}{\mathrm{\Gamma}(1-\alpha)}$
and $\tau_{0}(\lambda)=\frac{\beta\lambda^{-\frac{\beta}{\alpha}-1}e^{-\lambda \sum_{j=1}^{J}\theta_{j}\zeta^{\alpha}}}{\alpha\Gamma(1-\beta/\alpha)},$ and $M_{j}=M.$ For the next results we will need to recall the notion of Pitman-Yor finite dimensional distributions as in~\cite{Carlton}, this serves to also to connect our work to another recent line of investigations in~\cite{LijoiPYmulti,LijoiPYmultiJASA} and related references. Let $(Q_{1},\ldots,Q_{J})$ denote some disjoint partition of $\Omega$, such that $F_{0}(Q_{j}):=q_{j}.$ Then for $P_{\alpha,\theta}\sim\mathscr{PY}(\alpha,\theta,F_{0}),$ with $\theta>-\alpha,$  denoting a two parameter Pitman-Yor process with base-measure $F_{0}$~\cite{IJ2001,IJ2003,Pit96}, we set its finite-dimensional distribution $(P_{\alpha,\theta}(Q_{1}),\ldots,P_{\alpha,\theta}(Q_{J})):=(\tilde{P}_{\alpha,\theta}(q_{j}), j\in[J])\sim \mathrm{RS}(\mathbf{q},\alpha,\theta).$   We note in the purely stable case $(\tilde{P}_{\alpha,0}(q_{j}), j\in[J])\overset{d}=\left(\frac{\sigma_{j}q^{1/\alpha}_{j}}{\sum_{v=1}^{J}\sigma_{v}q^{1/\alpha}_{v}}, j\in[J]\right)$ for $(\sigma_{j}, j\in [J])$ iid stable $\alpha$ variables. 
\begin{rem}
Note the process $\sum_{j=1}^{J}\tilde{P}_{\alpha,\theta}(q_{j})\delta_{Y_{j}}$ with $q_{j}=1/J$ equates to Pitman-Yor multinomial processes with law $\mathrm{PY}(\alpha,\theta,J,F_{0})$ as in~\cite{LijoiPYmulti}[Proposition 1]. We can refer to the distribution in our case with general $q_{j}$ as $\mathrm{PY}(\alpha,\theta,\mathbf{q},F_{0})$.
\end{rem}
This leads to a highly tractable model with up to $J+2$ parameters $(\alpha,\beta,\zeta, q_{1},\ldots,q_{J}).$ 
In order to sample from the marginal distribution of the Poisson-GG-GG:
\begin{enumerate}
\item $\varphi\sim \mathrm{Poisson}({(\sum_{j=1}^{J}\theta_{j})}^{\beta/\alpha}[{(\zeta+M)}^{\beta}-\zeta^{\beta}]).$
\item $\psi_{j}(M)=\theta_{j}[{(\zeta+M)}^{\alpha}-\zeta^{\alpha}]$
\item $\tilde{C}_{j,k,l}\overset{iid}\sim \mathrm{MtP}(M,\tau_{j})$ with pmf
$
\frac{{(1-\frac{\zeta}{\zeta+M})}^{c}}
{(1-\frac{\zeta^{\alpha}}{{(\zeta+M)}^{\alpha}})}
\times
\frac{\alpha\Gamma(c-\alpha)}{c!\Gamma(1-\alpha)}
$
\item $\tilde{X}_{l}\sim \mathrm{MtP}([\sum_{j}^{J}\theta_{j}][{(\zeta+M)}^{\alpha}-\zeta^{\alpha}],\tau_{0})$ with pmf 
$
\frac{{(1-\frac{\zeta^{\alpha}}{{(\zeta+M)}^{\alpha}})}^{x_{l}}}
{(1-\frac{\zeta^{\beta}}{{(\zeta+M)}^{\beta}})} 
\times
\frac{\beta\Gamma(x_{l}-\frac{\beta}{\alpha})}{x_{l}!\alpha\Gamma(1-\frac{\beta}{\alpha})}
$
\item $X_{j,l}, j\in [J]|\tilde{X}_{l}=x_{l}\sim \mathrm{Multi}(x_{l},q_1,\ldots,q_{J}),$ for $q_{j}=\theta_{j}/\sum_{v=1}^{J}\theta_{v}$
\item $(\sum_{k=1}^{X_{j,l}}\tilde{C}_{j,k,l},j\in[J])\overset{d}=
(N_{j,l}, j\in[J])$ has joint pmf for $n_{l}=\sum_{j=1}^{J}n_{j,l}=1,2,\ldots$
\begin{equation}
\label{jointNformulastable}
\frac{{(1-\frac{\zeta}{{(\zeta+M)}})}^{n_{l}}}
{(1-\frac{\zeta^{\beta}}{{(\zeta+M)}^{\beta}})} \frac{\beta\Gamma(n_{l}-\beta)}{n_{l}!\Gamma(1-\beta)}
\times\frac{n_{l}!}{\prod_{j=1}^{J}n_{j,l}!}\mathbb{E}[\prod_{j=1}^{J}\tilde{P}^{n_{j,l}}_{\alpha,-\beta}(q_{j})]
\end{equation}
where $(\tilde{P}_{\alpha,-\beta}(q_{j}), j\in[J])\sim \mathrm{RS}(\mathbf{q},\alpha,-\beta),$
\end{enumerate}
Hence,  the joint distribution of $((\tilde{C}_{j,k,l}, k\in[X_{j,l}]), [X_{j,l}], j\in[J]), (N_{j,l}, j\in[J]))$ is, for $x_{l}=\sum_{j=1}^{J}x_{j,l},$ $\sum_{k=1}^{x_{j,l}}c_{j,k,l}=n_{j,l},$ with $x_{j,l}=0$ for $n_{j,l}=0,$:
\begin{equation}
\frac{\Gamma(x_{l}-\frac{\beta}{\alpha})\Gamma(1-\beta)\alpha^{x_{l}-1}}{\Gamma(1-\frac{\beta}{\alpha})\Gamma(n_{l}-\beta)}
\prod_{j=1}^{J}q_{j}^{x_{j,l}}\prod_{k=1}^{x_{j,l}}\frac{\Gamma(c_{j,k,l}-\alpha)}{\Gamma(1-\alpha)}
\times\frac{\beta\Gamma(n_{l}-\beta)}{n_{l}!\Gamma(1-\beta)}\frac{{(1-\frac{\zeta}{{(\zeta+M)}})}^{n_{l}}}
{(1-\frac{\zeta^{\beta}}{{(\zeta+M)}^{\beta}})}.
\label{specialstablestablederived1}
\end{equation}
For some relevant components of the posterior distribution given $\mathbf{Z}_{J}$
we have:
\begin{enumerate}
\item $\psi_{j}(M)=\theta_{j}[(\zeta+M)^{\alpha}-\zeta^{\alpha}],$ $\Psi_{0}(\sum_{j=1}^{J}\psi_{j}(M))={(\sum_{j=1}^{J}\theta_{j})}^{\beta/\alpha}[{(\zeta+M)}^{\beta}-\zeta^{\beta}].$
\item $B_{j,M}|B_{0,J}\sim \mathrm{GG}(\alpha,\theta_{j},\zeta+M;B_{0,J}),$ and $B_{0,J}\sim \mathrm{GG}(\beta/\alpha,1,\sum_{j=1}^{J}\theta_{j}{(\zeta+M)}^{\alpha};F_{0}).$
\item $\sum_{j=1}^{J}\theta_{j}{(\zeta+M)}^{\alpha}H_{l}\overset{d}=G_{\tilde{X}_{l}-\frac{\beta}{\alpha}}|\tilde{X}_{l}\sim \mathrm{Gamma}(\tilde{X}_{l}-\frac{\beta}{\alpha},1)$
\item The vector $(\hat{\sigma}_{j,l}(H_{l}), j\in[J])|(X_{j,l}, j\in[J])$ is equal in distribution to
\begin{equation}
\left(\frac{T^{(j,l)}_{\alpha}(q_{j}G_{\tilde{X}_{l}-\frac{\beta}{\alpha}}) + G_{n_{j,l} - X_{j,l}\alpha}}{\zeta + M}, j\in[J]\right).
\label{gengammasum2}
\end{equation}
where $G_{n_{j,l}-X_{j,l}\alpha}\overset{ind}\sim \mathrm{Gamma}(n_{j,l}-X_{j,l}\alpha)$ and $T^{(j,l)}_{\alpha}(y)$ are independent generalized gamma variables such that 
$\mathbb{E}[{\mbox e}^{-tT^{(j,l)}_{\alpha}(y)}]=e^{-y[{(1+t)}^{\alpha}-1]}.$ 
\item Equivalently the variables in~\eqref{gengammasum2} have the representation
\begin{equation}
\frac{G_{n_{l}-\beta}}{\zeta+M}\times \left(R_{X_{l}}\tilde{P}_{\alpha,\tilde{X}_{l}\alpha-\beta}(q_{j})+\tilde{D}_{j,l},j\in[J]\right)
\label{gengammasum3}
\end{equation}
where $G_{n_{l}-\beta}\sim \mathrm{Gamma}(n_{l}-\beta,1)$ independent of $(\tilde{P}_{\alpha,\tilde{X}_{l}\alpha-\beta}(q_{j}), j\in[J])\sim \mathrm{RS}(\mathbf{q},\alpha,\tilde{X}_{l}\alpha-\beta)$ and conditionally independent of this~$(R_{\tilde{X}_{l}},(\tilde{D}_{j,l},j\in[J]))\sim \mathrm{Dirichlet}(\tilde{X}_{l}\alpha-\beta,(n_{j,l}-X_{j,l}\alpha,j\in[J])).$ 
\item $(X_{j,l}, j\in[J])|(N_{j,l}=n_{j,l}, j\in[J])$ are further sampled according to \eqref{specialstablestablederived}.
\end{enumerate}
\begin{rem} The distributional correspondence between~\eqref{gengammasum2} and \eqref{gengammasum3} follows from applications of~\citep{PY97}[Proposition 21] and the usual beta-gamma algebra for gamma and Dirichlet variables. See~\cite{JamesStick} for more on representations of the form in~\eqref{gengammasum2}.
\end{rem}
\subsubsection{Alternative sampler}\label{sec:altsample}
Alternatively one may exactly sample $(N_{j,l},j\in [J])$ as follows. If we ignore the precise grouping in~\eqref{specialstablestablederived1}, then it follows that we can use for $\mathbf{c}_{l}:=(c_{k,l},k\in [x_{l}]),$ 

$$p_{\alpha,-\beta}(\mathbf{c}_{l})=\frac{\Gamma(x_{l}-\frac{\beta}{\alpha})\Gamma(1-\beta)}{\Gamma(1-\frac{\beta}{\alpha})\Gamma(n_{l}-\beta)}p_{\alpha,0}(\mathbf{c}_{l})
$$

denoting the $\mathrm{EPPF}$ of an $(\alpha,-\beta)$ partition of $[n_{l}],$ corresponding to sampling from $P_{\alpha,-\beta}\sim \mathscr{PY}(\alpha,-\beta,F_{0}),$ say
 $(C_{k,l}, k\in [\tilde{X}_{l}])|N_{l}=n_{l}$, where $\tilde{X
}_{l}|N_{l}=\sum_{k=1}^{\tilde{X}_{l}}C_{k,l}=n_{l}$ corresponds to $\mathbb{P}(\tilde{X}_{l}=x_{l}|n_{l})=\mathbb{P}^{(n_{l})}_{\alpha,-\beta}(x_{l})$. From this,
there is the distribution of 
$(C_{k,l},k\in[\tilde{X}_{l}]), (X_{j,l},j\in[J])|N_{l}=n_{l}$:

$$p_{\alpha,-\beta}(\mathbf{c}_{l})\frac{x_{l}!}{\prod_{j=1}^{J}x_{j,l}!}\prod_{j=1}^{J}q_{j}^{x_{j,l}}.
$$

More precisely, given $(C_{k,l}, k\in [\tilde{X}_{l}], \tilde{X}_{l}=x_{l})|N_{l}=n_{l},$
an $(\alpha,-\beta)$ partition of $[n_{l}],$
there are i.i.d. categorical variables $(\mathscr{L}_{1},\ldots,\mathscr{L}_{x_{l}})$ with $\mathbb{P}(\mathscr{L}_{i}=j)=q_{j}, j\in[J],$ with pairings 
$((C_{k,l},\mathscr{L}_{k}), k\in[x_{l}])$ where the sets  $((C_{k,l}:\mathscr{L}_{k}=j)\overset{d}=
(\tilde{C}_{j,k,l}, k\in[X_{j,l}])$ for $(X_{j,l}=\sum_{k=1}^{\tilde{X}_{l}}\mathbb{I}_{\{\mathscr{L}_{k}=j\}}, j\in[J])|\tilde{X}_{l}=x_{l}, N_{l}=n_{l}
\sim \text{Multi}(x_{l},q_{1},\ldots, q_{J}).$ Hence, given $N_{l}=n_{l},$

$$(\sum_{k=1}^{\tilde{X}_{l}}
C_{k,l}\mathbb{I}_{\{\mathscr{L}_{k}=j\}},j\in[J])\overset{d}=(\sum_{k=1}^{X_{j,l}}\tilde{C}_{j,k,l},j\in[J])
$$

following the distribution of $(N_{j,l}=n_{j,l},j\in [J])|N_{l}=n_{l},$ 
as indicated in \eqref{jointNformulastable} where $\mathbb{P}(N_{l}=n_{l})=
\frac{{(1-\frac{\zeta}{{(\zeta+M)}})}^{n_{l}}}
{(1-\frac{\zeta^{\beta}}{{(\zeta+M)}^{\beta}})} \frac{\beta\Gamma(n_{l}-\beta)}{n_{l}!\Gamma(1-\beta)}.$
\begin{rem}As can be seen from comparison with \eqref{specialstablestablederived1}, $\mathbb{E}[\prod_{j=1}^{J}\tilde{P}^{n_{j,l}}_{\alpha,-\beta}(q_{j})]$ is obtained as a special case of the general joint moment formula described in ~\cite{IJ2003}[Section 3.2, eq. (15)]. 
\end{rem}

\subsection{The marginal distribution of $(N_{j,l},j\in[J])$}\label{marginal}
We now provide descriptions of the general marginal process. First note that $(N_{j,l},j\in[J])|(\hat{\sigma}_{j,l}(H_{l})=t_{j}, j\in[J]),H_{l}$ is $\mathrm{tP}(t_{1}M_{1},\ldots,t_{J}M_{J}),$ and hence there is a joint distribution of  $(N_{j,l},\hat{\sigma}_{j,l}(H_{l}),j\in[J]),H_{l},$ with $n_{l}:=\sum_{j=1}^{J}n_{j,l}=1,2,\ldots,$ 
$$
\frac{\prod_{j=1}^{J}M^{n_{j,l}}_{j}t^{n_{j,l}}_{j}e^{-t_{j}M_{j}}}
{(1-{\mbox e}^{-\sum_{j=1}^{J}M_{j}t_{j}})\prod_{j=1}^{J}n_{j,l}!}
\frac{(1-{\mbox e}^{-\sum_{j=1}^{J}M_{j}t_{j}})\prod_{j=1}^{J}\eta_{j}(t_{j}|\lambda)}{(1-{\mbox e}^{-\lambda\sum_{j=1}^{J}\psi_{j}(M_{j})})}f_{H_{l}}(\lambda).
$$
For the next result define vectors, $\mathbf{k_J} = (k_1, k_2, \ldots, k_J) \in \mathbb{N}^J$, and $\mathbf{n}_l = (n_{1,l}, \ldots, n_{J,l})$. Where, $k_j(n_{j,l}) = 0$ if $n_{j,l} = 0$ and $k_j(n_{j,l}) \in [n_{j,l}]$ otherwise. Let the collection of all such vectors $\mathbf{k_J}$ be denoted by $\mathcal{K}_J({\mathbf{n}_l})$.

Now recall the identity,
$\int_0^{\infty} t^{n_{j,l}}e^{-tM_j}\eta_j(t|\lambda)\,dt = e^{-\lambda\psi_{j}(M_{j})}\Xi^{[n_{j,l}]}_j(\lambda,M_{j})$, for the results that follow.

\begin{prop}\label{post:marginalofN}
There are the following descriptions of the marginal and conditional distributions of $(N_{j,l}=n_{j,l}, j\in[J])$ defined for $n_{l}=\sum_{j=1}^{J}n_{j,l}=1,2,\ldots,,$ and each $l\in[\varphi]$,for $\varphi\sim \mathrm{Poisson}(\Psi_{0}(\sum_{j=1}^{J}\psi_{j}(M_{j}))).$
\begin{enumerate} 
    \item The conditional distribution of $(N_{j,l}=n_{j,l}, j\in[J])| H_{l}=\lambda$ can be expressed as,
    $$
    \frac{e^{-\lambda\sum_{j=1}^{J}\psi_{j}(M_{j})}\prod_{j=1}^{J}M^{n_{j,l}}_{j}\Xi^{[n_{j,l}]}_{j}(\lambda,M_{j})}{(1-e^{-\lambda\sum_{j=1}^{J}\psi_{j}(M_{j})})\prod_{j=1}^{J}n_{j,l}!}
    $$
\item The marginal of  $(N_{j,l}=n_{j,l}, j\in[J])$ has the following sum representation, with $x_{l}:=\sum_{j=1}^{J} x_{j,l},$ and $\mathbf{x}_{j,l}=(x_{j,l},j\in[J])$:
 $$
\frac{\sum_{\mathbf{x}_{j,l} \in \mathcal{K}_J(\mathbf{n}_l)}\Psi^{(x_{l})}_{0}(\sum_{j=1}^{J} \psi_j(M_j)) \prod_{j=1}^{J} p^{[n_{j,l}]}(k_j | \tau_{j}, M_j) \Xi^{[n_{j,l}]}_j(1, M_j) M^{n_{j,l}}_{j}}
{\Psi_0\left(\sum_{j=1}^{J} \psi_j(M_j)\right) \prod_{j=1}^{J}n_{j,l}!}
$$
where 
$\Psi^{(x_{l})}_{0}(\sum_{j=1}^{J} \psi_j(M_j))=\int_{0}^{\infty} \lambda^{x_{l}}e^{-\lambda \sum_{j=1}^{J} \psi_j(M_j)}\tau_0(\lambda)d\lambda$ as can be found in~\cite{Pit97}. Also, $p^{[n_{j,l}]}(\mathbf{c}_{j,l}|\tau_{j}, M_j)\Xi^{[n_{j,l}]}_j(1,M_j)=\prod_{k=1}^{x_{j,l}}\psi^{(c_{j,k,l})}_{j}(M_{j}).$
\end{enumerate}
\end{prop}
\subsection{Prediction rules}\label{predict}
We present the prediction rule for  $(Z^{(M_{j}+1)}_{j}, j\in[J])|\mathbf{Z}_{J},$ a vector of one customer from each group. Other descriptions follow by utilizing Theorems~\ref{postPoissonHIBP} and Proposition~\ref{sumsampleprop}.

\begin{prop}\label{Prediction} Let for each $j\in [J]$, $Z^{(M_{j}+1)}_{j}|B_{j}\sim \mathrm{PoiP}(B_{j}). B_{j}\sim \mathrm{CRM}(\tau_{j},B_{0}), B_{0}\sim \mathrm{CRM}(\tau_{0},F_{0})$. Then the predictive distribution of $(Z^{(M_{j}+1)}_{j}, j\in[J])|\mathbf{Z}_{J}$ is equal to that of 
$(\tilde{Z}^{(M_{j}+1)}_{j}+\sum_{l=1}^{r}\mathscr{P}^{(M_{j}+1)}_{j,l}(\sum_{k=1}^{X_{j,l}}S_{j,k,l})\delta_{\tilde{Y}_{l}}, j\in[J]),$
where $\mathscr{P}^{(M_{j}+1)}_{j,l}(\lambda)\overset{ind}\sim \mathrm{Poisson}(\lambda),$ $\tilde{Z}^{(M_{j}+1}_{j}|B_{0,J}+\sum_{l=1}^{r}H_{l}\delta_{\tilde{Y}_{l}}\sim \mathrm{IBP}(C_{j},\tau^{(M_{j})}_{j},B_{0,J}+\sum_{l=1}^{r}H_{l}\delta_{Y_{l}}),$ $C_{j}\sim \mathrm{Poisson}(s).$ Hence given $\mathbf{Z}_{J},$
 \begin{equation}
    \label{postdisint}
(\tilde{Z}^{(M_{j}+1)}_{j}, j\in[J])\overset{d}=    
    (\sum_{v=1}^{\varphi^{*}}[\sum_{k=1}^{X^{*}_{j,l}}\tilde{C}^{(*,1)}_{j,k,v}]\delta_{Y_{v}} + \sum_{l=1}^{r}\left[\sum_{k=1}^{\mathscr{P}^{*}_{j,l}(H_{l})}\tilde{C}^{(*,2)}_{j,k,l}\right]\delta_{\tilde{Y}_{l}}, j\in J)
    \end{equation}
where $(Y_{v})\overset{iid}\sim F_{0},$  $\varphi^{*}\sim \mathrm{Poisson}([\Psi_{0}(\sum_{j=1}^{J}\psi_{j}(M_{j}+1))-\Psi_{0}(\sum_{j=1}^{J}\psi_{j}(M_{j}))]),$  
\begin{enumerate}
\item $X^{*}_{l}:=\sum_{j=1}^{J}X^{*}_{j,l}\sim \mathrm{MtP}(\sum_{j=1}^{J}
[\psi_{j}(M_{j}+1)-\psi_{j}(M_{j})],\tau_{0,J})$ with probability mass function 
$\frac{{(\sum_{j=1}^{J}[\psi_{j}(M_{j}+1)-
\psi_{j}(M_{j})])}^{x_{l}}\Psi^{(x_{l})}_{0}(\sum_{j=1}^{J}\psi_{j}(M_{j}+1))}
{x_{l}![\Psi_{0}(\sum_{j=1}^{J}\psi_{j}(M_{j}+1))-\Psi_{0}(\sum_{j=1}^{J}\psi_{j}(M_{j}))]}$
\item $(X^{*}_{j,l},j\in[J])|X^{*}_{l}=x^{*}_{l}\sim \mathrm{Multi}(x^{*}_{l}; q_{1},\ldots, q_{J})$, for $q_{j}=\frac{[\psi_{j}(M_{j}+1)-\psi_{j}(M_{j})]}{\sum_{v=1}^{J}[
\psi_{v}(M_{v}+1)-\psi_{v}(M_{v})]}.$
\item For each $j\in[J],$ $(\tilde{C}^{(*,1)}_{j,k,v})$ and $(\tilde{C}^{(*,2)}_{j,k,l}))$ are iid collections of $\mathrm{MTP}(1,\tau^{(M_{j})}_{j})$  variables with pmf $\frac{\psi^{(c)}_{j}(M_{j}+1)}{[\psi_{j}(M_{j}+1)-\psi_{j}(M_{j})]c!}$
\item $\mathscr{P}^{*}_{l}(H_{l})=\sum_{j=1}^{J}\mathscr{P}^{*}_{j,l}(H_{l})\sim \mathrm{Poisson}(H_{l}\sum_{j=1}^{J}[\psi_{j}(M_{j}+1)-\psi_{j}(M_{j})])$
\item  $(\mathscr{P}^{*}_{j,l}(H_{l}), j\in[J])|\mathscr{P}^{*}_{l}(H_{l})=n^{*}_{l}\sim \mathrm{Multi}(n^{*}_{l}; q_{1},\ldots, q_{J})$
\item $H_{l},(X_{j,l},j\in[J])|\mathbf{Z}_{J}$ as in Theorem~\ref{postPoissonHIBP}.
\end{enumerate}
\end{prop}
If one is interested in adding groups, $J$ to $J+1,$  as in ~\cite{ZhouPadilla2016} the prediction rule is as follows
\begin{cor}Suppose that $Z^{(1)}_{J+1}|B_{J+1}\sim \mathrm{PoiP}(B_{J+1})$ and $B_{J+1}\sim \mathrm{CRM}(\tau_{J+1},B_{0}).$ Then $Z^{(1)}_{J+1}|\mathbf{Z}_{J}$ is equal in distribution to $\tilde{Z}^{(1)}_{J+1}$ expressed as a univariate version of \eqref{postdisint} with $M_{J+1}=0$ in place of $M_{j},$ $X^{*}_{J+1,l}\sim \mathrm{MtP}(\psi_{J+1}(1),\tau_{0,J})$ and $\mathscr{P}^{*}_{J+1,l}(H_{l})\sim \mathrm{Poisson}(H_{l}\psi_{J+1}(1)).$ Replace $\varphi^{*}$ with $\varphi^{*}_{J+1}\sim \mathrm{Poisson}(\Psi_{0}(\psi_{J+1}(1)+\sum_{j=1}^{J}\psi_{j}(M_{j}))-\Psi_{0}(\sum_{j=1}^{J}\psi_{j}(M_{j}))$  To obtain the simpler rule, corresponding to~\cite{ZhouPadilla2016}, where $\mathbf{Z}_{J}:=(Z^{(1)}_{j},j\in[J]),$ set $(M_{j}=1,j\in[J]).$ 
\end{cor}

\subsubsection{Prediction rule Poisson-GG-GG }\label{predictGG}
For the Poisson-GG-GG model in section~\ref{PGGGG}, where recall $M_{j}=M,$ and let $\mathrm{NB}(p,t)$ denote a Negative Binomial variable with parameters $(p,t),$ we have
\begin{enumerate}
\item $\varphi^{*}\sim\mathrm{Poisson}({(\sum_{j=1}^{J}\theta_{j})}^{\beta/\alpha}[{(\zeta+M+1)}^{\beta}-{(\zeta+M)}^{\beta}]),$ and for each $l\in[\varphi^{*}],$
\item
$\mathscr{P}^{*}_{l}(H_{l})\sim 
\mathrm{NB}(\frac{{(\zeta+M)}^{\alpha}}{{(\zeta+M+1)}^{\alpha}},\tilde{X}_{l}-\frac{\beta}{\alpha})$
\item $\mathscr{P}^{(M_{j}+1)}_{j,l}(\sum_{k=1}^{X_{j,l}}S_{j,k,l})\sim 
\mathrm{NB}(\frac{{(\zeta+M)}}{{(\zeta+M+1)}},n_{j,l}-X_{j,l}\alpha)$
\item $X^{*}_{l}\sim \mathrm{MtP}([\sum_{j=1}^{J}\theta_{j}[{(\zeta+M+1)}^{\alpha}-{(\zeta+M)}^{\alpha}],\tau_{0,J})$ with pmf 
$
\frac{{(1-\frac{{(\zeta+M)}^{\alpha}}{{(\zeta+M+1)}^{\alpha}})}^{x_{l}}}
{(1-\frac{{(\zeta+M)}^{\beta}}{{(\zeta+M+1)}^{\beta}})} 
\times
\frac{\beta\Gamma(x_{l}-\frac{\beta}{\alpha})}{x_{l}!\alpha\Gamma(1-\frac{\beta}{\alpha})}
$
\item $(\tilde{C}^{(*,1)}_{j,k,v})$ and $(\tilde{C}^{(*,2)}_{j,k,l}))$ are iid collections of $\mathrm{MtP}(1,\tau^{(M)}_{j})$  variables with pmf $
\frac{{(\zeta+M+1)}^{\alpha-c}}{c![{(\zeta+M+1)}^{\alpha}-{(\zeta+M)}^{\alpha}]}\frac{\alpha\Gamma(c-\alpha)}{\Gamma(1-\alpha)},$ not depending on $j\in[J].$
\item Sample $(X_{j,l},j\in[J])|(N_{j,l}=n_{j,l},j\in[J]),$ with reference to \eqref{specialstablestablederived1}, from distribution proportional to,
$
\frac{\Gamma(x_{l}-\frac{\beta}{\alpha})\alpha^{J-1}}{\Gamma(1-\frac{\beta}{\alpha})\Gamma(x_{l})}
\frac{\Gamma(x_{l})}{\prod_{j=1}^{J}\Gamma(x_{j,l})}
\prod_{j=1}^{J}q_{j}^{x_{j,l}}\mathbb{P}^{(n_{j,l})}_{\alpha,0}(x_{j,l})
$
\end{enumerate}

\subsection{Recovering latent variables in the Poisson HIBP case and coupling HIBP}\label{recovery}
While Proposition~\ref{poissonequivalence} shows that it is not necessary to recover the  individual variables 
 $(\tilde{C}^{(i)}_{j,k,l}),$ in Theorem~\ref{postPoissonHIBP}, as noted they have utility in terms of constructing coupled HIBP, leading to recovery of latent variables necessary for describing posterior distributions in those cases, and certainly any application where information recovery for individuals is of interest. We also note from a computational perspective it may be easier to include sampling from the $(\tilde{C}_{j,kl}).$ In particular one can employ the following scheme to recover these variables.

Conditional on $(N_{j,l}=n_{j,l}, j\in[J], H_{l}=\lambda)$ or $(N_{j,l}=n_{j,l}, j\in[J])$: 
\begin{enumerate}
\item[1.] Sample $(\tilde{C}_{j,k,l}, k\in [X_{j,l}], X_{j,l})$ from its conditional distribution as described in Propositions~\ref{GibbsEPPF} and \eqref{prop:GibbsEPPF} or as in Theorem~\ref{postPoissonHIBP}. 
\item [2.] Then, $(\tilde{C}^{(i)}_{j,k,l}, i\in[M{j}])|\tilde{C}_{j,k,l}=c_{j,k,l}) \sim \text{Multi}(c_{j,k,l}, \frac{1}{M_{j}}, \ldots, \frac{1}{M_{j}})$. This allocates the sampled $\tilde{C}_{j,k,l}$ into the individual $\tilde{C}^{(i)}_{j,k,l}$. 
\item[3.] Additionally, $S_{j,k,l}|\tilde{C}_{j,k,l}=c_{j,k,l}$ is sampled from a density proportional to $s^{c_{j,k,l}}e^{-sM_{j}}\tau_{j}(s)$. 
\item[4.] Posterior quantities for coupled Bernoulli HIBP can be constructed by $\sum_{k=1}^{X_{j,l}}\mathbb{I}_{\{\tilde{C}^{(i)}_{j,k,l}>0\}}$ and other HIBP using $\mathbb{P}{A_{j}}(da|\pi^{-1}(1-\exp(-S_{j,k,l})))$, as in Proposition~\ref{coupled}.   
\end{enumerate}

\section{Experiments with Simulated data}
We have shown that our simplified Poisson HIBP in sections \ref{PGGGG} and \ref{predictGG} has components corresponding to variants of the Pitman-Yor multinomial distributions arising in~\cite{LijoiPYmulti,LijoiPYmultiJASA}. Readers may consult those works for interesting real data applications that clearly demonstrate the practical broad applicability of our results, even in the simplest cases. As a by-product, we offer new ways to sample from such processes. This should also stimulate ideas about further applications involving our latent feature framework.

Here, in this section, we provide simulations for the more general form of the Poisson-GG-GG as described in Proposition~\ref{PropgenPGG} and present a sampling algorithm for posterior inference. Specifically, we present Poisson-GG-GG HIBP with the following parameterization.
\[
&\tau_0(\lambda) = \mathrm{GG}(\lambda\given\alpha_0, \theta_0/\alpha_0, \zeta_0) = \frac{\theta_0}{\Gamma(1-\alpha_0)}\lambda^{-\alpha_0-1}e^{-\zeta_0\lambda}, \quad B_0 \sim \mathrm{CRM}(\tau_0, F_0), \\
&\tau_j(s) = \mathrm{GG}(s\given \alpha_j, \theta_j/\alpha_j, \zeta_j)= \frac{\theta_j}{\Gamma(1-\alpha_j)}s^{-\alpha_j-1}e^{-\zeta_js}, \quad B_j \given B_0 \sim \mathrm{CRM}(\tau_j, B_0) 
\text{ for } j \in [J],\\
&Z_j^{(i)}\given B_j \sim \mathrm{PoiP}(B_j) \text{ for } i \in [M_j], j \in [J].\nonumber
\]
In the main text, we describe the simulation and synthetic document classification result based on it. In the Section C in the supplementary material, we describe all the details, see also Proposition~\ref{PropgenPGG}, including the derivations for marginal distribution, predictive distribution, and Gibbs sampling procedures for parameter estimation.

\subsection{Simulation}
We first visualize the word-count matrices generated from our model
to see how the parameters of the model affect the overall patterns in the generated count matrices. The generated word count matrices are visualized in Figure~\ref{fig_Z_vis}. For word count matrices generated with larger values of $\alpha_0$, which encourage a stronger power-law distribution in word count frequencies, there tends to be a greater total number of words (Figure~\ref{fig_Z_vis}, right). Additionally, these matrices exhibit more diverse patterns of word occupancy, with fewer overlaps between groups. 

\begin{figure}
\centering
\includegraphics[width=0.32\linewidth]{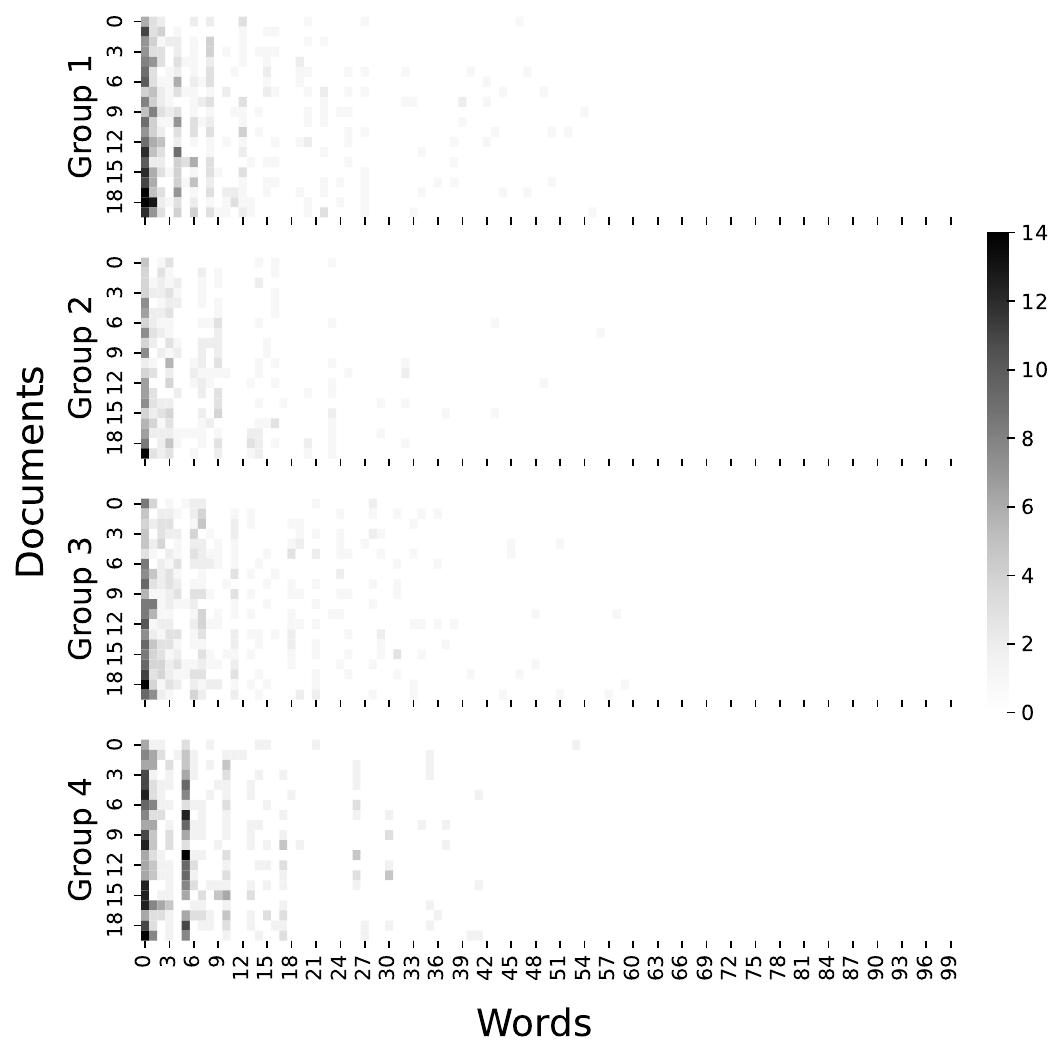}
\includegraphics[width=0.32\linewidth]{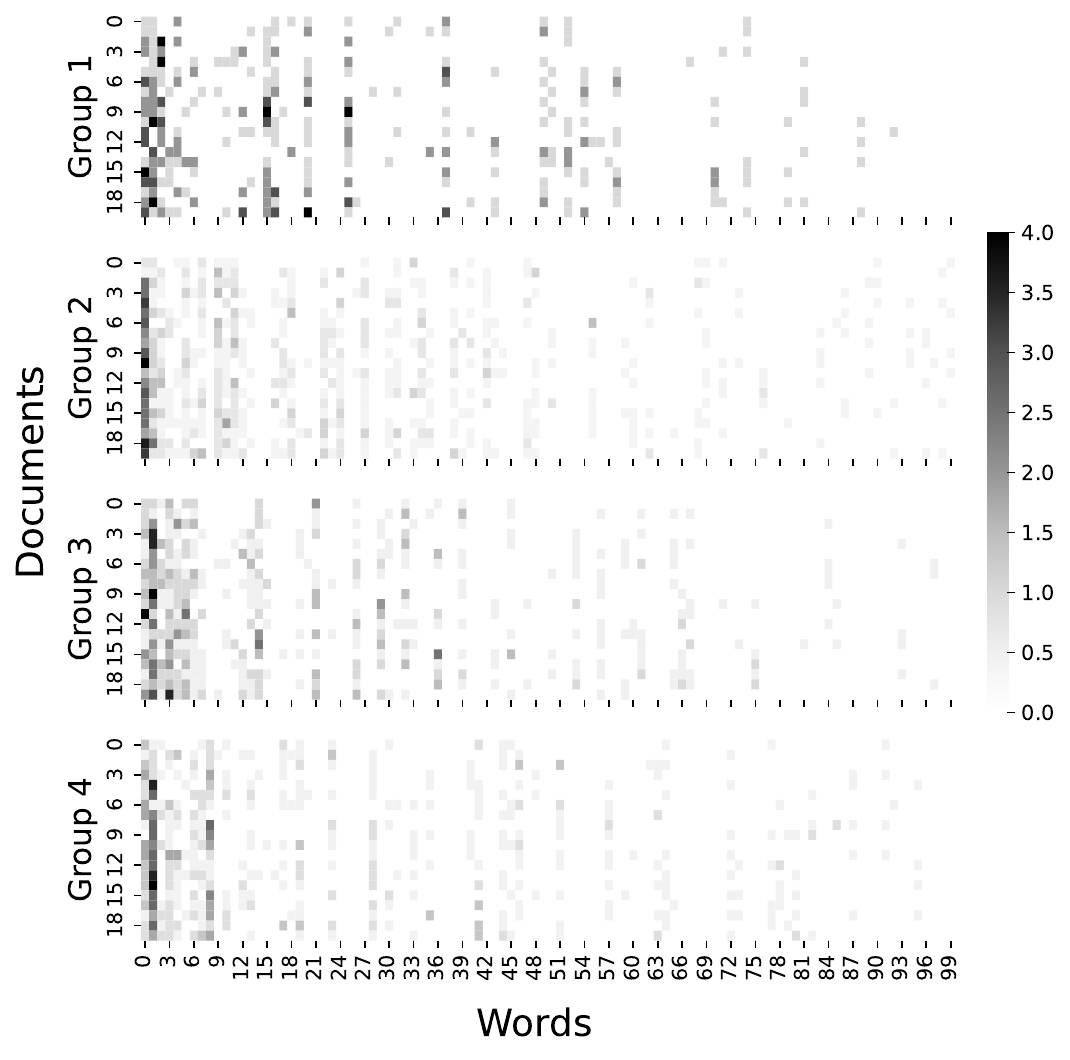}
\includegraphics[width=0.32\linewidth]{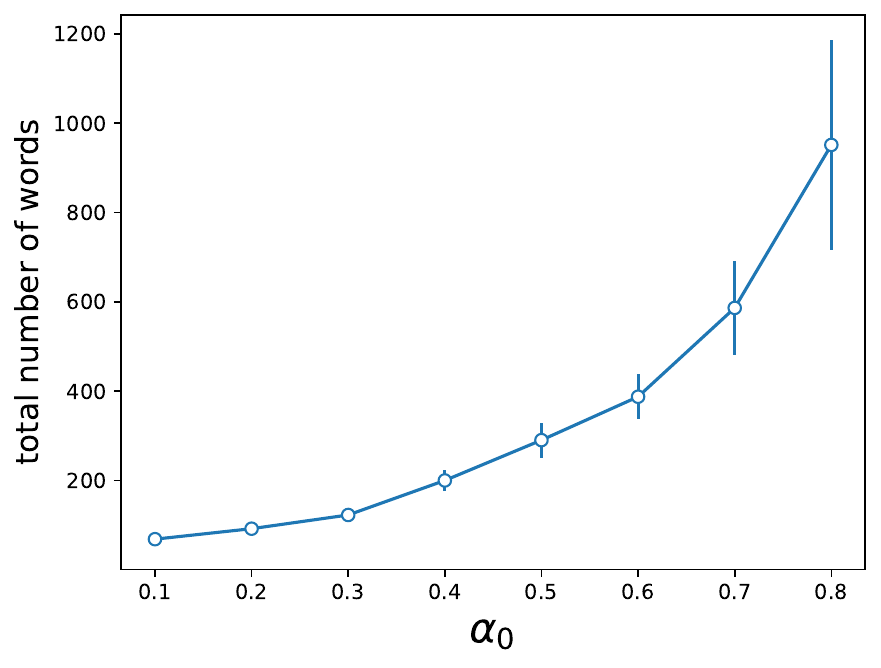}
\caption{(Left, middle) Word count matrices generated from Poisson-GG-GG HIBP model with parameters $J=4$, $\theta_0=10, \zeta=1$, $M_j=20$, $\theta_j \sim \mathrm{Unif}(2, 4)$, $\alpha_j \sim \mathrm{Unif}(0.2, 0.5)$, $\zeta_j=1$, and $\beta_j=1$ for $j\in [J]$.  The left one is generated from $\alpha_0=0.2$ and the middle one is generated from $\alpha_0=0.7$. Showing top 100 words only. (Right) Average total number of words across groups with respect to varying $\alpha_0$ value. We generate 10 matrices with each $\alpha_0$ value and compute the averages.}
\label{fig_Z_vis}
\end{figure}

\subsection{Posterior Inference for the parameters}
\label{gg_gg_poisson_hibp_posterior_inference}
Given a set of observed documents $\mathbf{Z}_J = ((Z_j^{(i)}, i\in [M_j]), j\in [J])$, one might be interested in inferring the parameters $\phi := (\alpha_0, \zeta_0, \theta_0, ((\alpha_j, \zeta_j, \theta_j, \beta_j), j \in [J]))$. In order to compute the marginal distribution required for the parameter estimation, we may need 
$\varphi=r$,
 $(X_{j,l}=x_{j,l}, j \in [J], l \in [r])$, and $(\tC^{(i)}_{j,k,l}=c_{j,k,l}^{(i)}, i\in [M_j], k\in[x_{j,l}], j\in [J], l\in[r])$, or the aggregated word counts 
 $\tN_{j,l}^{(i)}=\sum_{k=1}^{X_{j,l}} \tC_{j,k,l}^{(i)}$. In our experiments, we estimate the parameters using the marginal distribution of $(\varphi, N_{j,l}, X_{j,l})$ which can be derived from the results in the main text. The derivations for Poisson-GG-GG HIBP model are given in Section C, along with the Gibbs sampling algorithm for the parameter estimation. 
 Note that we should recover the counts $(X_{j,l})$, which can be done based on Theorem \ref{postPoissonHIBP}.

For demonstration, we generate documents $\mathbf{Z}_J$ with the true parameter values set as $\theta_0=10, \zeta=1, J=4, M_j=200, \theta_j\sim \mathrm{Unif}(2.0, 4.0), \alpha_j \sim \mathrm{Unif}(0.1, 0.6), \zeta_j=1, \beta_j=1$ for $j\in [J]$, and $\alpha_0 \in \{0.2, 0.6\}$. We fix $\zeta_0, (\zeta_j, j\in [J])$ to $1$, and inferred $\phi := (\theta_0, \alpha_0, ((\theta_j, \alpha_j), j\in [J]))$.  In the posterior inference, for all parameters, we place uninformative priors, such as $\mathrm{Gamma}(1, 1)$ on $\theta_0$ and $\mathrm{logit}\, \calN(0, 1)$ on $\alpha_0$. See Section C.II in the supplementary material for the detailed prior settings.

For each of the two datasets (one with $\alpha_0=0.2$ and another with $\alpha_0=0.6$), we run three independent MCMC chains with each chain run for 30,000 steps with 15,000 burn-in samples. Figure~\ref{fig:mcmc_posteriors} presents posterior samples for $\theta_0$ and $\alpha_0$. The results for the remaining parameters are shown in Section C.V in the supplementary material.
Overall, the proposed Gibbs sampler correctly identified the posteriors that are properly concentrated around the true values. 

\begin{figure}
\centering
\includegraphics[width=0.47\linewidth]{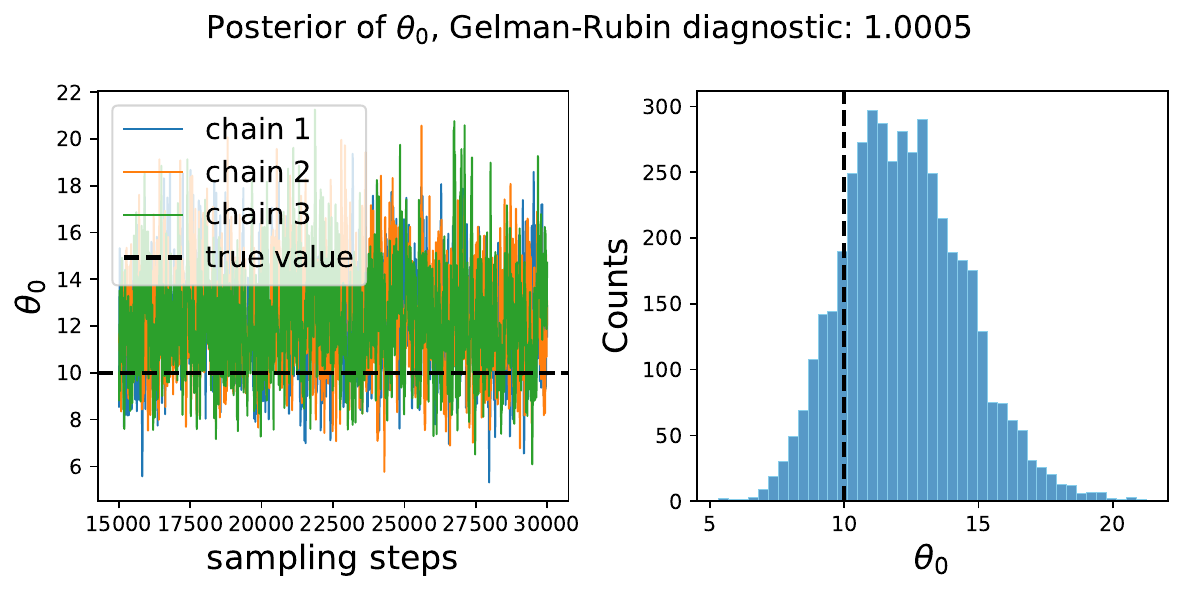}
\includegraphics[width=0.47\linewidth]{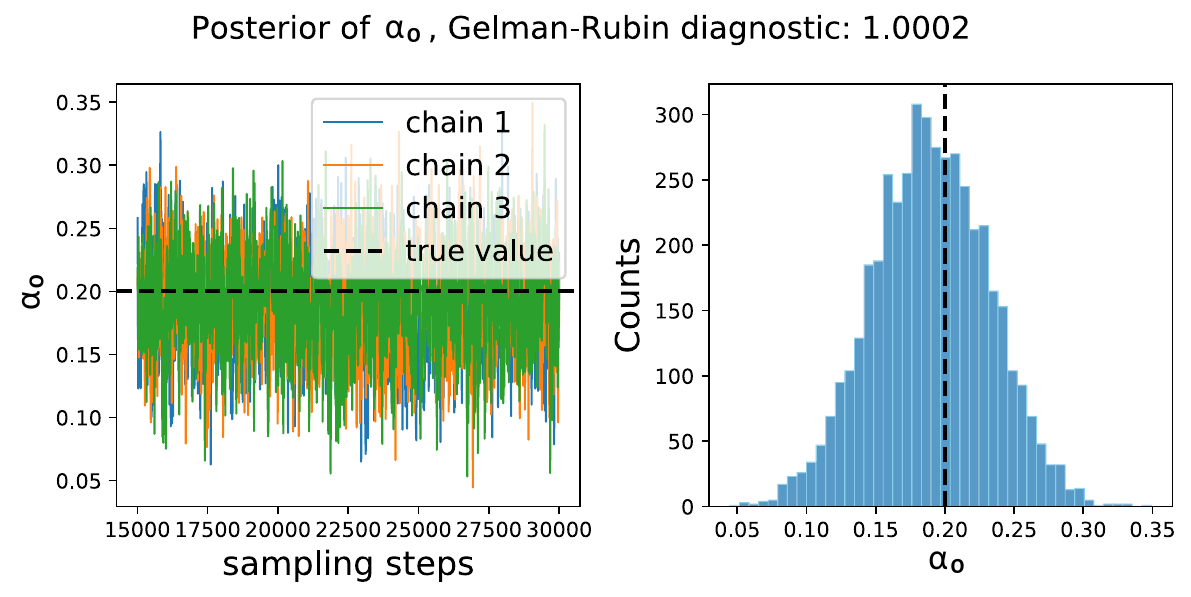}
\includegraphics[width=0.47\linewidth]{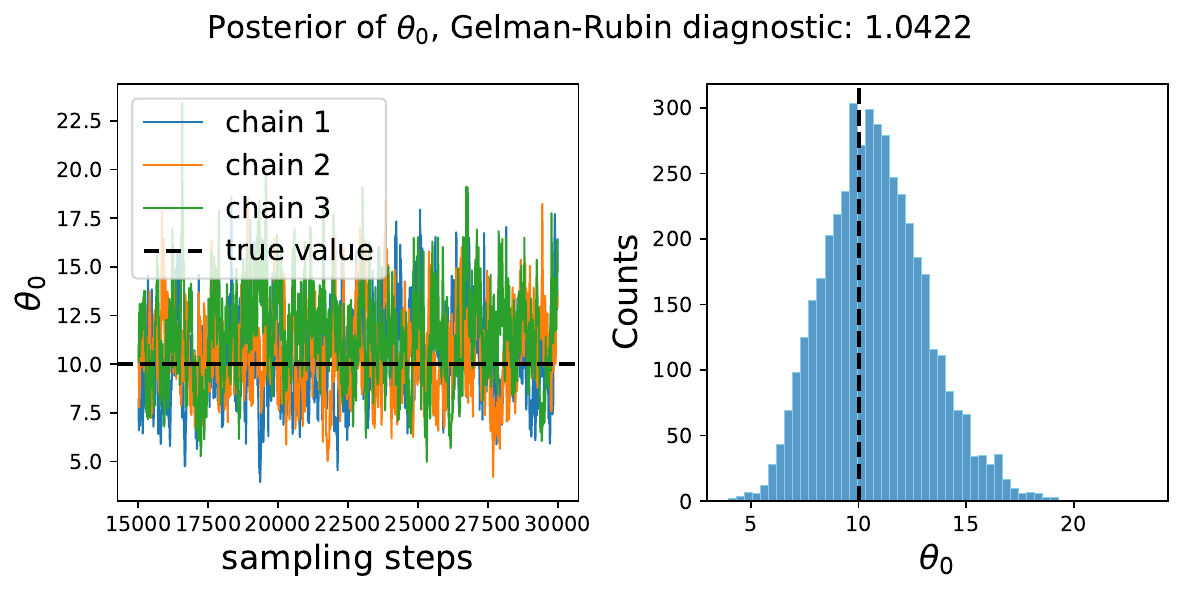}
\includegraphics[width=0.47\linewidth]{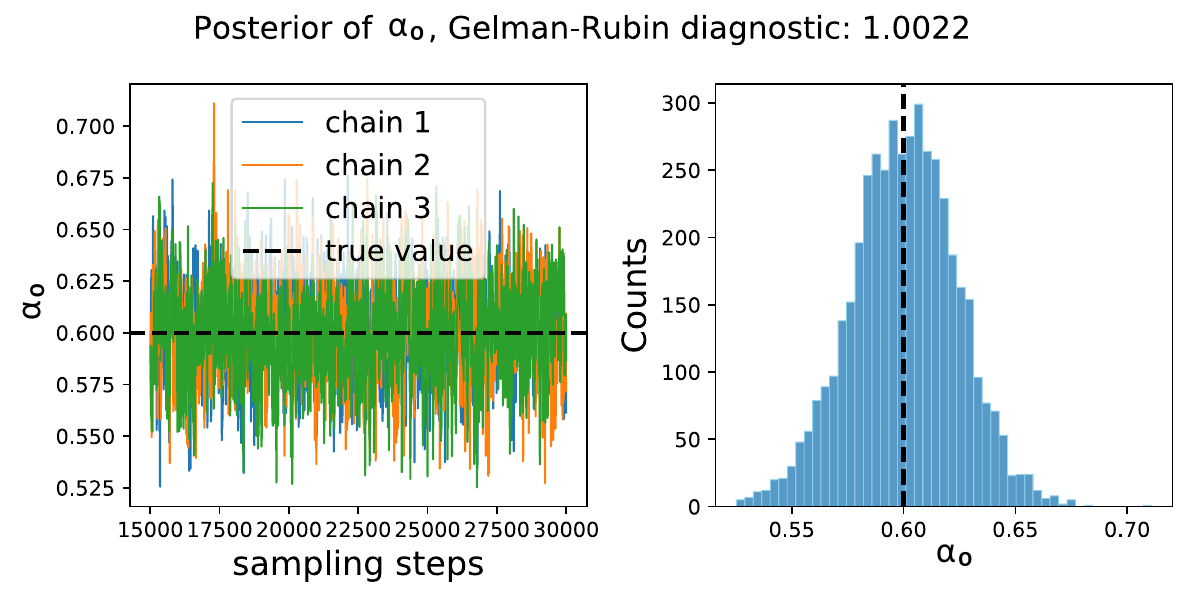}
\caption{Posterior samples for $\theta_0$ and $\alpha$. Top row shows the results for the data generated with $\alpha_0=0.2$ and bottom row shows the ones with $\alpha_0=0.6$. We present Gelman-Rubin convergence diagnostic~\cite{BrooksGeneral, GelmanInference}; Gelman-Rubin diagnostic values are clearly smaller than 1.2, indicating convergences.}
\label{fig:mcmc_posteriors}
\end{figure}

\subsection{Prediction and Classification}

\begin{figure}[h]
\centering
\includegraphics[width=0.4\linewidth]{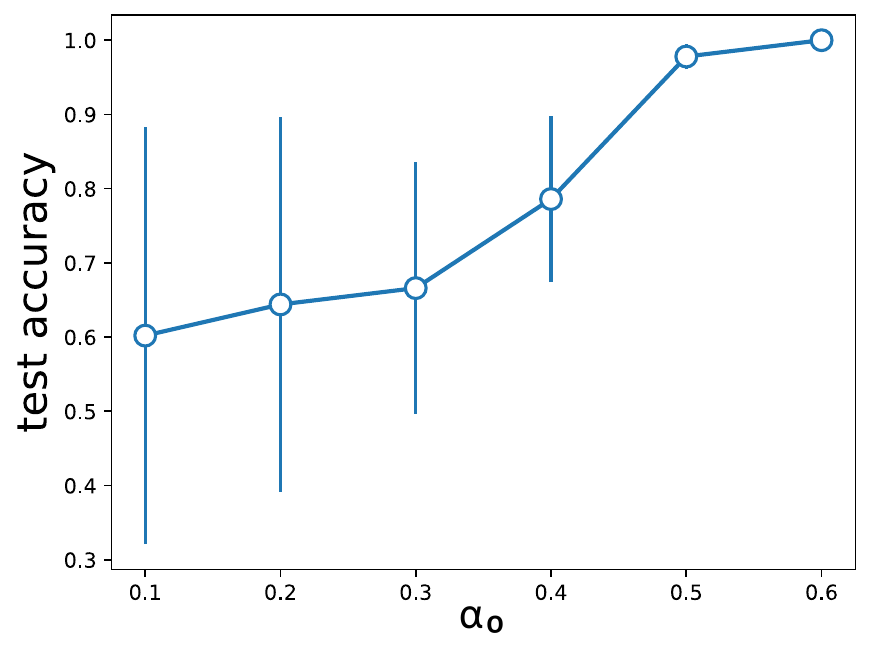}
\includegraphics[width=0.4\linewidth]{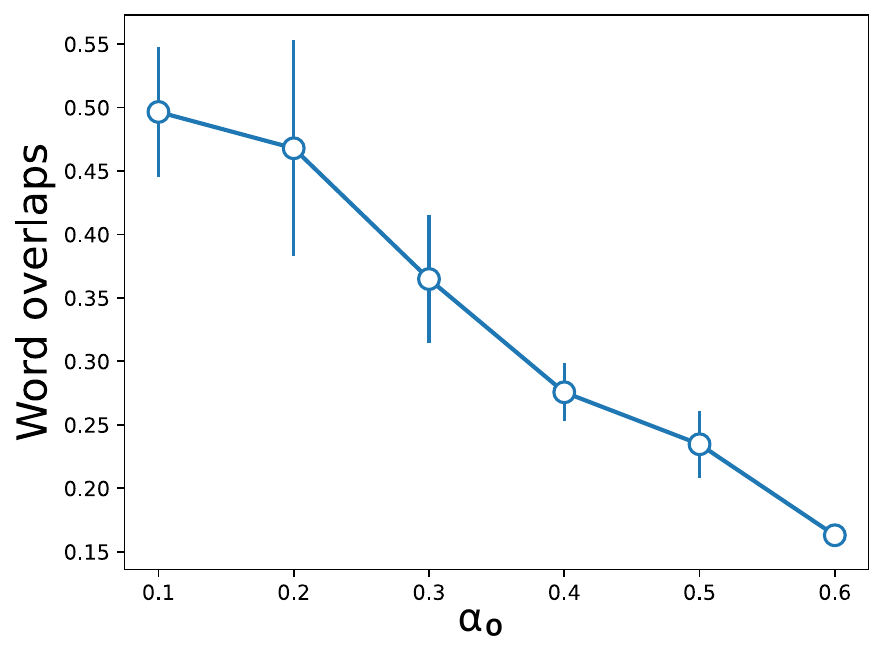}
\caption{Left: average group classification accuracies for test documents generated with varying $\alpha_0$ values. Right: average overlap of words (see the main text for definition) across groups for the data generated with varying $\alpha_0$ values.}
\label{fig:gg_gg_poisson_hibp_classify}
\end{figure}

We further demonstrate that we can use our model for classifying documents by their groups. Given a set of documents $\bZ_J$, one can simulate or compute predictive probabilities of future documents using Proposition \ref{Prediction}. Specifically, for a new document $Z_*$, for each $j\in [J]$, we compute the probability that $Z_*$ was generated from group $j$, $\mathbb{P}(Z_*=Z_{j}^{(M_j+1)} | \mathbf{Z}_J)$. Then we can classify $Z_*$ to be in group if $\arg\max_{j} \mathbb{P}(Z_*=Z_{j}^{(M_j+1)} | \mathbf{Z}_J)$. 

As an illustration, we generate training documents $\mathbf{Z}_J$ from the model we described in the previous section, and further generate $100$ test documents per each group, using the prediction rule described in Proposition \ref{Prediction}. See Section C.III for the detailed description for the prediction rules for Poisson-GG-GG HIBP models. We then run the posterior inference algorithm with the training data, collect the estimated parameters, and compute the predictive probabilities of each test document belonging to one of $J$ groups. Here, similar to training, computing the predictive probabilities also requires working with the aggregated counts. In Section C.IV, we present a Gibbs sampling algorithm to recover the latent variables required for computing predictive probabilities.

Setting the parameters as in Section \ref{gg_gg_poisson_hibp_posterior_inference} except for $\alpha_0$, we generate $\mathbf{Z}_J$ with $\alpha_0 \in \{0.1, 0.2, 0.3, 0.4, 0.5, 0.6\}$. For each $\alpha_0$ we generate five datasets, constituting 30 datasets in total. (Figure \ref{fig:gg_gg_poisson_hibp_classify}, left) summarizes the average test group classification accuracies w.r.t. $\alpha_0$ values used to generate data. We see that documents generated with larger $\alpha_0$ values are easier to classify; this is presumably because documents generated with larger $\alpha_0$ values are less likely to share words across groups and have their distinct words. To quantitatively see this, we measure the average overlap of words across groups,
\begin{align*}
\mathrm{Overlap}(\mathbf{Z}_J\,|\, \varphi=r, ((X_{j,l}=x_{j,l}))) = \frac{1}{\binom{J}{2}}\sum_{j=1}^{J-1}\sum_{j'=j+1}^{J} \sum_{l=1}^r \frac{\1{x_{j,l} > 0 \,\wedge\, x_{j',l}>0}}{r},
\end{align*}
that is, average portion of words co-occurring between all pairs of groups. (Figure \ref{fig:gg_gg_poisson_hibp_classify}, right) presents the $\mathrm{Overlap}(\mathbf{Z}_J)$ values averaged across five datasets per each $\alpha_0$. As $\alpha$ increases, the average overlap decreases, indicating that documents from different groups become easier to distinguish.

\clearpage

\appendix

\renewcommand\thesection{\Alph{section}}
\renewcommand\thesubsection{\thesection.\Roman{subsection}}

\section{Proofs and discussions}

\subsection{Results in sections 2 and 3}
Proposition 2.1, as mentioned in the main text, is just a rephrasing of the results in~\cite{James2017}[Section 3]. However, it is not difficult for the reader to check everything directly after the decomposition of the L\'evy densities $\rho_{j}(s)=[1-\pi_{j}(s)]^{M_{j}}\rho_{j}(s)+(1-[1-\pi_{j}(s)]^{M_{j}})\rho_{j}(s),$ where the condition  $\Phi_{j}(M_{j})=\int_{0}^{\infty}(1-{[1-\pi_{j}(s)]}^{M_{j}})\rho_{j}(s)ds<\infty,$ induces the decomposition of $\mu_{j},$ where the finite component induces a compound Poisson representation based on the components $(\xi_{j}, S_{j,\ell},\omega_{j,\ell},\ell\in[\xi_{j}])$ as described in Proposition 2.1, where the $M_{j}$ customers will sample from.  The rest, as explained is manipulation of the joint, conditional and marginal distributions which involves random vectors. 
\begin{rem}A point to note is that the results rely on the unicity of the jumps $(S_{j,\ell})$ and hence unique pairs $(S_{j,\ell},\omega_{j,\ell}),$ regardless of whether $(\omega_{j,\ell},\ell\in[\xi_{j}])$ has ties or not.
\end{rem}
 Propositions 2.2 and 2.3 follow again as a special case of~\cite{James2017}[Section 5, Propositions 5.1 and 5.2], for a $J$-length multivariate IBP process
say $Z^{(1)}_{0}=\mathbf{X}$ and  based on a single univariate $\mathrm{CRM}(\tau_{0},F_{0}),$ $\mu_{0}:=B_{0}$, where the notations on the left hand side are otherwise meant to match with respective quantities in\cite{James2017}[Section 5, Propositions 5.1 and 5.2]. Note $F_{0}$ plays the role of the variable denoted as $B_{0}$ in that text. The results are seen more transparently by following the decomposition of $\tau_{0},$ dictated by the zero set, and applying Bayes' rule appropriately. Here, one uses that $(\xi_{1},\ldots,\xi_{J})|\lambda$ are conditionally independent $\mathrm{Poisson}(\Phi_{j}(M_{j})\lambda)$ variables for $j\in[J]$, such that 
$$
\mathbb{P}(\xi_{1}=0,\ldots,\xi_{J}=0|\lambda) = e^{-\lambda \sum_{j=1}^{J} \Phi_{j}(M_{j})}
$$ 
denotes the spike distribution, which determines the region where features will or will not be sampled according to the process $\mathbf{X}$, leading to the decomposition of $\tau_{0}(\lambda)=e^{-\lambda \sum_{j=1}^{J} \Phi_{j}(M_{j})}\tau_{0}(\lambda)+(1-e^{-\lambda \sum_{j=1}^{J} \Phi_{j}(M_{j})})\tau_{0}(\lambda)$, which identifies the distribution of the variables $H_{l},$ $\varphi$ and $B_{0,J}.$ One then works with the multivariate zero-truncated Poisson distribution $\mathrm{tP}(\lambda (\Phi_{1}(M_{1}),\ldots,\Phi_{J}(M_{J})))$ such that for $l\in[\varphi]$ there is the joint distribution of $(X_{j,l}=x_{j,l}, j\in[J]), H_{l}$ with $x_{l}:=\sum_{j=1}^{J} x_{j,l}=1,2,\ldots,$ 
$$
\frac{\prod_{j=1}^{J} \Phi^{x_{j,l}}_{j}(M_{j}) \lambda^{x_{j,l}} e^{-\lambda \Phi_{j}(M_{j})}}{(1 - e^{-\lambda \sum_{j=1}^{J} \Phi_{j}(M_{j})}) \prod_{j=1}^{J} x_{j,l}!} f_{H_{l}}(\lambda).
$$
One then uses this to get the appropriate joint, marginal and posterior distributions, similar to the discussion in Proposition 2.1. Again a key observation for later usage is that the variables $(H_{l},\tilde{X}_{l})$ may be interpreted as in~\cite{Pit97}. 

\subsection{Proof of Theorem 3.1 }

Here we use the fact that given $B_{0}$ $(Z^{(i)}_j,i\in[M_{j}],j\in[J])$ is equal in distribution to $(\sum_{l=1}^{\infty}[\sum_{k=1}^{\xi_{j,l}}\tilde{A}^{(i)}_{j,k,l}]\delta_{Y_{l}},i\in[M_{j}],j\in[J])$ which can be viewed as a multivariate IBP whose zero set is in fact the same as the zero-set of the process $\mathbf{X}.$  Hence there is the same decomposition of $\tau_{0}$ as in Proposition 2.2. One can then work with, for each $l,$ the joint distributions of $(\sum_{k=1}^{\xi_{j,l}}\tilde{A}^{(i)}_{j,k,l},i\in[M_{j}],j\in[J])$ and $(\xi_{j,l},j\in[J])|\lambda$ similar to the case of Proposition 2.2. One then mixes over the distribution of $\mathbf{X}$ as in Proposition 2.2 to complete the result. The result in effect is from an application of \cite{James2017}[Proposition 5.2], however it is not necessary to have the precise conditional distributions of the sums $(\sum_{k=1}^{\xi_{j,l}}[\tilde{A}^{(i)}_{j,k,l}],i\in[M_{j}])$ which may be quite complex.
\qed

As to results in section 3.1, details for Proposition 3.1 follow from the non-hierarchical Poisson setting in~\cite{James2017}[Section 4.2], as well as the corresponding descriptions in~\cite{Pit97}, coupled with Proposition 2.2. Proposition 3.2 just follows from the appropriate transformation of the jumps of a CRM based on Levy density $\rho_{j}$ to that of $\tau_{j}$ via the map as indicated. 

\begin{rem}\label{remarkPost}In principle the setup in the proof of Theorem 3.1, again based on~\cite{James2017}[Section 5], is enough to establish an expression for the posterior distribution of $B_{0}|(Z^{(i)}_j,i\in[M_{j}],j\in[J])$ by working with the joint distribution of $(\sum_{k=1}^{X_{j,l}}\tilde{A}^{(i)}_{j,k,l},\tilde{Y}_{l}, l\in[r],\varphi=r, i\in [M_{j}], j\in J),$ and $(H_{l},l\in[r])$ and then using Bayes rule to obtain the posterior distribution of the jumps $H_{l},
l\in[r])$ given $(\sum_{k=1}^{X_{j,l}}\tilde{A}^{(i)}_{j,k,l}, l\in[r], i\in [M_{j}], j\in J),$ as otherwise the distribution of $B_{0,J}$ remains unchanged. However, as mentioned in the main text, we seek tractable descriptions of the posterior distribution that are conducive to practical sampling and help our understanding of its structural properties. While the resulting descriptions are correct, it does not meet that standard without more information about the distribution of the sums. Further refinements are also necessary, as we have demonstrated in the Poisson HIBP case.
\end{rem}

\subsection{Proof of Theorem 4.1}
As in Remark~\ref{remarkPost}, it is correct that one can use \cite{James2017}[Section 5] to obtain a description of \( B_{0} | \left( \sum_{k=1}^{X_{j,l}} \tilde{A}^{(i)}_{j,k,l}, \tilde{Y}_{l}, l \in [r], \varphi = r, i \in [M_{j}], j \in [J] \right) \), where, in any event, the posterior distribution of the \( (H_{l}) \) depends on the sums through the \( (X_{j,l}) \). Instead, we work with the enlarged space including the individual variables \( \left( X_{j,l}, (A^{(i)}_{j,k,l}, k \in [X_{j,l}], l \in [r], i \in [M_{j}], j \in [J]) \right) \). It is then evident that \( B_{0} \), given this information, has the posterior distribution given in Proposition~2.3. Given \( B_{0} \) and these variables, for each \( j \in [J] \), we may equate \( (S_{j,k,l}, \tilde{Y}_{l}, k \in [X_{j,l}], l \in [\varphi]) \) with \( ((S_{j,\ell}, \omega_{j,\ell}), \ell \in [\xi_{j}]) \) and substitute appropriately in Proposition 2.1. Noting further that \( \mu_{j,M_{j}} \), given \( B_{0} \) and this extra information, is $\mathrm{CRM}(\rho^{(M_{j})},B_{0})$ as in Proposition 2.1. We then apply the distribution of \( B_{0} \)  given \( (X_{j,l}, (A^{(i)}_{j,k,l}, k \in [X_{j,l}], l \in [r], i \in [M_{j}], j \in [J]) \) to show that the joint distribution of \( (\mu_{j,M_{j}}, j \in [J]) \) is component-wise \( \mathrm{CRM}(\rho^{(M_{j})},B_{0,J} + \sum_{l=1}^{r} H_{l} \delta_{\tilde{Y}_{l}}) \) and hence \( \mu_{j,M_{j}} \overset{d}{=} \tilde{\mu}_{j,M_{j}} + \sum_{l=1}^{r} \tilde{\sigma}_{j,l}(H_{l}) \delta_{\tilde{Y}_{l}} \). These arguments give the representations in statements 1 and 2 of Theorem 4.1. Statement 3 of that result then rounds out the appropriate descriptions. \qed

\subsection{Results for the Poisson HIBP in section 5}
We now establish the results in Section 5. Instead of using Theorem 4.1 we approach Theorem 5.1 and results by a more direct approach exploiting the sum properties of Poisson variables, and directly the infinite-divisibility of $\sigma_{j,l}(\lambda).$ We first summarize some properties/representations of the Poisson HIBP which showcase its mixed Poisson representations, and describe $(B_{j},j\in[J]),B_{0})$ as a multivariate CRM as in~\cite{James2017}. Using these representatons, as we shall show,  an alternative posterior analysis follows from \cite{James2017}[Section 5]. However tractability still requires the usage of augmentation. Here again we use the notation~$\mathscr{P}(\lambda)\sim \mathrm{Poisson}(\lambda)$

Facts for $Z^{(i)}_{j}|B_{j}\sim\mathrm{PoiP}(B_{j}), B_{j}|B_{0}\sim\mathrm{CRM}(\tau_{j},B_{0}),B_{0}\sim\mathrm{CRM}(\tau_{0},F_{0})$:
\begin{enumerate}\label{Poissonfact}
    \item Let $\eta_{j}(t_{j}|\lambda)$ denote the density of $\sigma_{j,l}(\lambda)=\sum_{k=1}^{\infty}s_{j,k,l}$ with corresponding Laplace exponent $\lambda\psi_{j}(t)$ and $(s_{j,k,l})_{\{k\ge 1\}}$ determined by $\lambda\tau_{j}(s).$
    \item $(B_{j}\in[J],B_{0})$ is a multivariate $\mathrm{CRM}$ in the sense of~\cite{James2017}[Section 5] with joint L\'evy density,$\prod_{j=1}^{J}\eta_{j}(t_{j}|\lambda)\tau_{0}(\lambda),$ and is equivalent in distribution to 
$
(\sum_{l=1}^{\infty}\sigma_{j,l}(\lambda_{l})\delta_{Y_{l}},j\in[J],\sum_{l=1}^{\infty}\lambda_{l}\delta_{Y_{l}})
$      
    \item $\mathbf{Z}_{J}:=(Z^{(i)}_{j}, i\in[M_{j}],j\in [J])$ is equal in distribution to 
    \begin{equation}
    \label{poissonrep}
   (\sum_{l=1}^{\infty}\mathscr{P}^{(i)}_{j,l}(\sigma_{j,l}(\lambda_{l}))\delta_{Y_{l}}, i\in [M_{j}],j\in[J]),
    \end{equation}
    where $\mathscr{P}^{(i)}_{j,l}(\sigma_{j,l}(\lambda_{l}))=\sum_{k=1}^{\infty}\mathscr{P}_{j,k,l}(s_{j,k,l})\sim \mathrm{Poisson}(\sigma_{j,l}(\lambda_{l}))$
 \item The posterior distribution of $(B_{j}\in[J],B_{0})|\mathbf{Z}$ only depends on the information in the sum process 
$$ 
(\sum_{l=1}^{\infty}\mathscr{P}_{j,l}(M_{j}\sigma_{j,l}(\lambda_{l}))\delta_{Y_{l}},j\in[J])
 $$
 where $\mathscr{P}_{j,l}(M_{j}\sigma_{j,l}(\lambda_{l}))=\sum_{i=1}^{M_{j}}\mathscr{P}^{(i)}_{j,l}(\sigma_{j,l}(\lambda_{l})).$
\end{enumerate}

We now describe the relevant components for the Poisson process by solely utilizing the observed sums of this process, once again utilizing the framework in~\cite{James2017}[Section 5], but aiming for a clearer explanation through a method of decompositions. The sum representations, in terms of $\sigma_{j,l},$  facilitate this approach. For each fixed $l,$ given $(\sigma_{j,l}(\lambda)=t_{j},j\in[J]),\lambda_{l}=\lambda$ the zero set below is given by 
$$
\mathbb{P}(\mathscr{P}_{1,l}(M_{1}t_{1})=0,\ldots,\mathscr{P}_{J,l}(M_{J}t_{J})=0)= {\mbox e}^{-\sum_{j=1}^{J}M_{j}t_{j}}
$$
 
The next results follow:

\begin{enumerate}
    \item The event \( \{ l : \mathscr{P}^{(i)}_{j,l}(\sigma_{j,l}(\lambda_{l})) = 0), i \in [M_{j}], j \in [J] \} \) induces a decomposition of the joint L\'evy density of \( (B_{j}, j \in [J], B_{0}) \) as:
    \[
    e^{-\sum_{j=1}^{J} M_{j} t_{j}} \prod_{j=1}^{J} \eta_{j}(t_{j} | \lambda) \tau(\lambda) + \Psi_{0}( \sum_{j=1}^{J} \psi_{j}(M_{j})) \varrho(t_{1}, \ldots, t_{j} | \lambda) f_{H}(\lambda).
    \]

    \item The joint density \( \varrho(\mathbf{t} | \lambda) f_{H}(\lambda) \) is proper, where:
    \begin{equation}
    \label{jointcompounddensity}
    \varrho(\mathbf{t} | \lambda) = \frac{(1 - e^{-\sum_{j=1}^{J} M_{j} t_{j}}) \prod_{j=1}^{J} \eta_{j}(t_{j} | \lambda)}{(1 - e^{-\lambda \sum_{j=1}^{J} \psi_{j}(M_{j})})}.
    \end{equation}

    \item The joint density denotes the joint density of variables \( (\hat{\sigma}_{j,l}(H_{l}), j \in [J], H_{l}) \) which are iid for \( l \in [\varphi] \), where \( (H_{l}) \) are iid with common density 
    \[
    f_{H}(\lambda) = \frac{(1 - e^{-\lambda \sum_{j=1}^{J} \psi_{j}(M_{j})}) \tau_{0}(\lambda)}{\Psi_{0} \left( \sum_{j=1}^{J} \psi_{j}(M_{j}) \right)}.
    \]

    \item There is a decomposition of \( (B_{j}, j \in [J], B_{0}) \) as:
    \begin{equation}
    \label{postdisint}
    \left( B_{j,M_{j}} + \sum_{l=1}^{\varphi} \hat{\sigma}_{j,l}(H_{l}) \delta_{\tilde{Y}_{l}}, j \in [J], B_{0,J} + \sum_{l=1}^{\varphi} H_{l} \delta_{\tilde{Y}_{l}} \right).
    \end{equation}

    \item Where \( B_{0,J} \sim \mathrm{CRM}(\tau_{0,J}, F_{0}) \), \( B_{j,M_{j}} \sim \mathrm{CRM}(\tau^{(M_{j})}_{j}, B_{0,J}) \), and have the same distributions as in Theorem 5.1. The \( (\tilde{Y}_{l}) \) are iid \( F_{0} \) and represent the \( \varphi \sim \mathrm{Poisson}(\Psi_{0}(\sum_{j=1}^{J} \psi_{j}(M_{j}))) \) number of distinct features to be selected by the \( \sum_{j=1}^{J} M_{j} \) customers.

    \item \( (N_{j,l}, j \in [J]) | (\hat{\sigma}_{j,l} = t_{j}, j \in [J]), (H_{l}) \sim \mathrm{tP}((t_{j} M_{j}, j \in [J]) \).
\end{enumerate} 
 
What remains to establish Theorem 5.1 is to determine the marginal and posterior quantities to describe the overall posterior distribution of $(B_{j}, j\in [J], B_{0})|\mathbf{Z}_{J}$ based on the joint distribution of $(N_{j,l},j\in[J]), (\hat{\sigma}_{j,l}=t_{j}, j\in[J]), H_{l}),$ for each $l\in[\varphi]$
expressed as 
\begin{equation}
\frac{\prod_{j=1}^{J}M^{n_{j,l}}_{j}t^{n_{j,l}}_{j}e^{-t_{j}M_{j}}}{(1-e^{-\sum_{j=1}^{J}M_{j}t_{j}})\prod_{j=1}^{J}n_{j,l}!}\times \varrho(\mathbf{t}|\lambda)f_{H}(\lambda)
\end{equation}
It follows that the joint distribution of $(N_{j,l},j\in[J]), (\hat{\sigma}_{j,l}(H_{l})=t_{j}, j\in[J])| H_{l}=\lambda,$  has the form
\begin{equation}
\frac{\prod_{j=1}^{J}M^{n_{j,l}}_{j}t^{n_{j,l}}_{j}e^{-t_{j}M_{j}}\eta_{j}(t_{j}|\lambda)}{\prod_{j=1}^{J}n_{j,l}!}
\end{equation}
which is equal to
\begin{equation}
\prod_{j=1}^{J}\frac{t_{j}^{n_{j,l}}\eta^{[0]}_{j}(t_{j}|\lambda,M_{j})}{\Xi^{[n_{j,l}]}_j(\lambda,M_{j})}\times\frac{e^{-\lambda\sum_{j=1}^{J}\psi_{j}(M_{j})}\prod_{j=1}^{J}M^{n_{j,l}}_{j}\Xi^{[n_{j,l}]}_{j}(\lambda,M_{j})}{(1-e^{-\lambda\sum_{j=1}^{J}\psi_{j}(M_{j})})\prod_{j=1}^{J}n_{j,l}!}
\label{sigmaNjoint}
\end{equation}
where $\eta^{[0]}_{j}(t_{j}|\lambda,M_{j}):=e^{-t_{j}M_{j}}\eta_{j}(t_{j}|\lambda)e^{\lambda\psi_{j}(M_{j})},$
identifying the densities of $\hat{\sigma}_{j,l}(\lambda)|N_{j,l}=n_{j,l},H_{l}=\lambda$ as $\eta^{[n_{j,l}]}_{j}(t_{j}|\lambda,M_{j})$ which agrees with the description in Proposition 5.2. 

\subsubsection{Comments on Proposition 5.1}
That is the distribution of  $\hat{\sigma_{j,l}}(\lambda)|N_{j,l}=n_{j,l},H_{l}=\lambda$ is equivalent to the distribution of 
$\sigma_{j,l}(\lambda)|\mathscr{P}_{j,l}(M_{j}\sigma_{j,l}(\lambda))=n_{j,l},$ and for $n_{j,l}>0,$ is the same as the distribution of  $\sigma_{j,l}(\lambda)|\frac{\sum_{v=1}^{n_{j,l}}\mathbf{e}_{v,l}}{\sigma_{j,l}(\lambda)}=M_{j},$ for $(\mathbf{e}_{v,l})_{v\ge 1}$ an iid collection of exponential (1) variables,
This sets up relations of the variables $\tilde{C}_{j,k,l},X_{j,l}$ to corresponding variables arising in~\cite{JLP2,FoFZhou}, which either by direct comparison or as informed in the works of~\cite{PitmanPoissonMix} and \cite{Kolchin,Pit97} leads to the  finite Gibbs partitions in Proposition 5.1. 

\begin{rem}
As an aside,  formally we may normalize the jumps $(s_{j,k,l})_{\{k\ge 1}\}$ of $\sigma_{j,l}(\lambda)=\sum_{k=1}^{\infty}s_{j,k,l}$ to form Poisson-Kingman classes as in ~\cite{JLP2,Pit02}, then condition these on $N_{j,l}=n_{j,l}.$ Finer details will be discussed elsewhere. 
\end{rem}
 
\subsubsection{establishing Propositions 5.2, 5.4 and the remainder of Theorem 5.1}

As mentioned in the main text, Proposition 5.2 then follows immediately otherwise from the existing quoted results in \citep{James2002,JamesStick,JLP2}, yielding
$$
\hat{\sigma}_{j,l}(\lambda)\overset{d}=\tilde{\sigma}_{j,l}(\lambda)+\sum_{k=1}^{X_{j,l}}S_{j,k,l}
$$
The other expression in~\eqref{sigmaNjoint} corresponds to the joint distribution of $(N_{j,l}, j\in[J])|H_{l}=\lambda$ as it appears in 1. of Proposition 5.4. The representations for $\Xi^{[n_{j,l}]}_{j}(\lambda,M_{j}),$ as can be read from for instance~\cite{James2002,JLP2}, are due to the cumulant expansion of the $n_{j,l}$ moments of the infinitely divisible variables with exponentially tilted densities $\eta^{[0]}(t_{j}|\lambda,M_{j}).$ Proposition 5.4 is completed by integrating with respect to $f_{H},$ keeping in mind the meaning and marginal distribution of the variables $X_{j,l},\tilde{X}_{l}, C_{j,k,l}$ etc. deduced within our exposition from Proposition 3.1. The joint distribution of $(X_{j,l},(\tilde{C}_{j,k,l},k\in [X_{j,l}], l\in[r], j\in J),$ given $(N_{j,l}, j\in[J]),l\in[r])$ appearing in statement 5. of Theorem 5.1 can also be verified directly by removing the sums in the following equivalent representation of the marginal distribution of $(N_{j,l}=n_{j,l}, j\in[J])$ with $x_{l}:=\sum_{j=1}^{J} x_{j,l},$ (which needs to be further divided by $\prod_{j=1}^{J}n_{j,l}!$)
\begin{equation}
\label{abstractmarginalN}
\sum_{\mathbf{x}_{j,l} \in \mathcal{K}_J(\mathbf{n}_l)}\mathbb{P}(\tilde{X}_{l}=x_{l})\mathbb{P}((X_{j,l}=x_{j,l},j\in[J])|\tilde{X}_{l}=x_{l})\prod_{j=1}^{J}
\mathbb{P}(\sum_{i=1}^{x_{j,l}}\tilde{C}_{j,i,l}=n_{j,l})
\end{equation}
where, as in Proposition~3.1, $\tilde{X}_{l}\sim \mathrm{MTP}(\sum_{j=1}^{J}\psi_{j}(M_{j}),\tau_{0}),$ $X_{j,l},j\in[J])|\tilde{X}_{l}-x_{l}\sim 
\mathrm{Multi}(x_{l}; q_{1},\ldots, q_{J})$, for $q_{j}=\psi_{j}(M_{j})/\sum_{v=1}^{J}\psi_{v}(M_{v})$ and $(\tilde{C}_{j,i,l}, i\in[x_{j,l})\overset{iid}\sim \mathrm{MTP}(M_{j}, \tau_{j}).$ \eqref{abstractmarginalN} is easily deduced by first using item 2. of Proposition 5.4 and then the identity deduced from item 4 of Proposition 5.1.

\subsubsection{Results for the general prediction rule in Proposition 5.5}
Once the posterior distribution of \( (B_{j}, j \in [J], B_{0}) | \mathbf{Z}_{J} \) has been established as in Theorem 5.1 or the discussion above, it is straightforward to obtain the prediction rule in Proposition 5.5. We provide a few details. 

\begin{enumerate}
\item From \eqref{postdisint} given $B_{j,M_{j}} + \sum_{l=1}^{\varphi} \hat{\sigma}_{j,l}(H_{l}) \delta_{\tilde{Y}_{l}}$ it follows that 
$Z^{(M_{j}+1)}_{j}\sim \mathrm{PoiP}(B_{j,M_{j}} + \sum_{l=1}^{\varphi} \hat{\sigma}_{j,l}(H_{l}) \delta_{\tilde{Y}_{l}}).$
\item Hence there is a decomposition of $(Z^{(M_{j}+1)}_{j},j\in[J])|\mathbf{Z}_{J}$ as a vector 
$$
(\hat{Z}_{j},j\in[J])+(\sum_{l=1}^{r}\mathscr{P}_{j,l}(\hat{\sigma}_{j,l}(H_{l}))\delta_{\tilde{Y}_{l}}, j\in[J])
$$
\item where further from Proposition 5.2 or directly from Theorem 5.1, 
$$
\mathscr{P}_{j,l}(\hat{\sigma}_{j,l}(H_{l}))\overset{d}=
\mathscr{P}^{(2)}_{j,l}(\tilde{\sigma}_{j,l}(\lambda))+\mathscr{P}_{j,l}(\sum_{k=1}^{X_{j,l}}S_{j,k,l})
$$
where $\mathscr{P}^{(2)}_{j,l}$ are independent Poisson variables with random intensities as indicated
\item $\hat{Z}_{j}|B_{j,M_{j}}\sim \mathrm{PoiP}(B_{j,M_{j}})$ and $B_{j,M_{j}}|B_{0,J}\sim \mathrm{CRM}(\tau^{(M_{j})}_{j},B_{0,J})$ with $B_{0,J}\sim CRM(\tau_{0,J},F_{0})$
\item hence the vector $(\hat{Z}_{j},j\in[J]),$ has a representation of the form ~\eqref{poissonrep} is sampled according to Proposition 3.1, based on the variable $C^{(*,1)}_{j,k,l},X^{*}_{j,l}$ as indicated in Proposition Proposition 5.5
    \item \( \sum_{l=1}^{r} \mathscr{P}^{(2)}_{j,l}(\tilde{\sigma}_{j,l}(\lambda)) \delta_{\tilde{Y}_{l}} \) has distribution \( \mathrm{PoiP}\left( \sum_{l=1}^{\varphi} \tilde{\sigma}_{j,l}(H_{l}) \delta_{\tilde{Y}_{l}} \right) \). Hence, the vector is also generated by applying Proposition 3.1 based on \( C^{(*,2)}_{j,k,l}, \mathscr{P}_{j,l}(H_{l}) \).
\item  Equivalently, one may first consider \( (\tilde{Z}^{(M_{j}+1)}_{j}, j \in [J]) \) given \( B_{0} \), where as an application of Proposition 2.1, \( \tilde{Z}^{(M_{j}+1)}_{j} | \mathbf{Z}_{J}, B_{0} \) has distribution \( \mathrm{IBP}(C_{j}, \tau^{(M_{j})}_{j}, B_{0}) \) as described in Proposition 2.1, item 5. Then one uses the posterior representation of \( B_{0} \) and applies Proposition 3.1.
\end{enumerate}

\section{Some calculations for the simplified Poisson-GG-GG}
Although the results concerning the simplified Poisson-GG-GG are supposed to be simpler in details the refined results are by no means obvious. We provide a derivation of the marginal distribution and related quantities. 

\subsection{Poisson GG-GG marginal/joint derivation}
For the simplified Poisson-Generalized Gamma-Generalized Gamma model, that is 
Poisson-GG-GG,  recall that we set in~Proposition 5.3, $\alpha_{0}=\beta/\alpha$ for $0<\beta<\alpha,$ $\tau_{j}(s)=\theta_{j}\frac{\alpha s^{-\alpha-1}e^{-\zeta s}}{\mathrm{\Gamma}(1-\alpha)}$
and $\tau_{0}(\lambda)=\frac{\beta\lambda^{-\frac{\beta}{\alpha}-1}e^{-\lambda \sum_{j=1}^{J}\theta_{j}\zeta^{\alpha}}}{\alpha\Gamma(1-\beta/\alpha)},$ and $M_{j}=M.$ 
Recall that the general form of the joint distribution of  $(N_{j,l},\hat{\sigma}_{j,l}(H_{l}),j\in[J]),H_{l},$ with $n_{l}:=\sum_{j=1}^{J}n_{j,l}=1,2,\ldots,$ is,
$$
\frac{\prod_{j=1}^{J}M^{n_{j,l}}_{j}t^{n_{j,l}}_{j}e^{-t_{j}M_{j}}}
{(1-{\mbox e}^{-\sum_{j=1}^{J}M_{j}t_{j}})\prod_{j=1}^{J}n_{j,l}!}
\frac{(1-{\mbox e}^{-\sum_{j=1}^{J}M_{j}t_{j}})\prod_{j=1}^{J}\eta_{j}(t_{j}|\lambda)}{(1-{\mbox e}^{-\lambda\sum_{j=1}^{J}\psi_{j}(M_{j})})}f_{H_{l}}(\lambda).
$$
Here in this case
$$
f_{H_{l}}(\lambda)=\frac{(1-e^{-\lambda\sum_{j=1}^{J}\theta_{j}[{(\zeta+M)}^{\alpha}-{\zeta}^{\alpha}]})\beta\lambda^{-\frac{\beta}{\alpha}-1}e^{-\lambda \sum_{j=1}^{J}\theta_{j}\zeta^{\alpha}}}{\alpha\Gamma(1-\frac{\beta}{\alpha}){(\sum_{j=1}^{J}\theta_{j})}^{\frac{\beta}{\alpha}}[{(\zeta+M)}^{\beta}-{\zeta}^{\beta}]}
$$
Also $(\hat{\sigma}_{j,l}(\lambda), j\in[J])$ has joint density expressed as
$$
\frac{(1-{\mbox e}^{-M\sum_{j=1}^{J}t_{j}})\prod_{j=1}^{J}\eta_{j}(t_{j}|\lambda)}{(1-e^{-\lambda\sum_{j=1}^{J}\theta_{j}[{(\zeta+M)}^{\alpha}-{\zeta}^{\alpha}]})}
$$
where
$$
\eta_{j}(t_{j}|\lambda)dt_{j}={\mbox e}^{-t_{j}\zeta}\mathbb{P}(\sigma_{j}(\theta_{j}\lambda)\in dt_{j})e^{\theta_{j}\lambda\zeta^{\alpha}}
$$
where $\sigma_{j}(\theta_{j}\lambda)\overset{d}=\sigma_{j}\theta^{\frac{1}{\alpha}}_{j}\lambda^{\frac{1}{\alpha}}$ where $\sigma_{j}$ is an $\alpha$-stable random variable with Laplace transform $\mathbb{E}[{\mbox e}^{-t\sigma_{j}}]={\mbox e}^{-t^{\alpha}}.$ It follows that we can always write 
$$
\frac{\prod_{j=1}^{J}M^{n_{j,l}}_{j}t^{n_{j,l}}_{j}e^{-t_{j}M_{j}}}
{\prod_{j=1}^{J}n_{j,l}!}=
\frac{M^{n_{l}}(\sum_{j=1}^{J}t_{j})^{n_{l}}{\mbox e}^{-M\sum_{j=1}^{J}t_{j}}}{n_{l}!}\times \frac{n_{l}!}{\prod_{j=1}^{J}n_{j,l}!}\prod_{j=1}^{J}p^{n_{j,l}}_{j}
$$
where $p_{j}=t_{j}/\sum_{v=1}^{J}t_{v}.$ Set $\sigma=\sum_{j=1}^{J}\sigma_{j}\theta^{\frac{1}{\alpha}}_{j}.$ Then the marginal distribution of $(N_{j,l},j\in [J])|H_{l}=\lambda$ can be expressed as the expectation over 
$$
\frac{M^{n_{l}}\sigma^{n_{l}}\lambda^{\frac{n_{l}}{\alpha}}{\mbox e}^{-(M+\zeta)\sigma \lambda^{1/\alpha}}}{n_{l}!(1-{\mbox e}^{-\lambda\sum_{j=1}^{J}\psi_{j}(M_{j})})}\times \frac{n_{l}!}{\prod_{j=1}^{J}n_{j,l}!}\prod_{j=1}^{J}P^{n_{j,l}}_{j}\times {\mbox e}^{\lambda\sum_{j=1}^{J}\theta_{j}\zeta^{\alpha}}
$$
where $P_{j}=\sigma_{j}\theta^{\frac{1}{\alpha}}_{j}/\sigma.$ However we may integrate over $f_{H_{l}}$ first which amounts to the integral
$$
\frac{\int_{0}^{\infty}\lambda^{\frac{n_{l}}{\alpha}}{\mbox e}^{-(M+\zeta)\sigma \lambda^{1/\alpha}}\beta\lambda^{-\frac{\beta}{\alpha}-1}d\lambda}{\alpha\Gamma(1-\frac{\beta}{\alpha}){(\sum_{j=1}^{J}\theta_{j})}^{\frac{\beta}{\alpha}}[{(\zeta+M)}^{\beta}-{\zeta}^{\beta}]}=
\frac{\beta\Gamma(n_{l}-\beta)\sigma^{-n_{l}+\beta}{(M+\zeta)}^{-n_{l}+\beta}}
{\Gamma(1-\frac{\beta}{\alpha}){(\sum_{j=1}^{J}\theta_{j})}^{\frac{\beta}{\alpha}}[{(\zeta+M)}^{\beta}-{\zeta}^{\beta}]}
$$
Note that 
$$
\mathbb{E}[\sigma^{\beta}]=(\sum_{j=1}^{J}\theta_{j})^{\beta/\alpha}\frac{\Gamma(1-\frac{\beta}{\alpha})}{\Gamma(1-\beta)}
$$
and $\tilde{\sigma}=\sigma/(\sum_{j=1}^{J}\theta_{j})^{\frac{1}{\alpha}}$ is an $\alpha$ stable variable with scale $1.$ Furthermore, we may write $P_{j}\overset{d}=\sigma_{j}
q^{\frac{1}{\alpha}}_{j}/\tilde{\sigma},$ where for clarity
$$
\mathbb{E}[\tilde{\sigma}^{\beta}]=\frac{\Gamma(1-\frac{\beta}{\alpha})}{\Gamma(1-\beta)}
$$
and write
$$
\frac{\beta\Gamma(n_{l}-\beta)\sigma^{-n_{l}+\beta}{(M+\zeta)}^{-n_{l}+\beta}}
{\Gamma(1-\frac{\beta}{\alpha}){(\sum_{j=1}^{J}\theta_{j})}^{\frac{\beta}{\alpha}}[{(\zeta+M)}^{\beta}-{\zeta}^{\beta}]}=
\frac{\beta\Gamma(n_{l}-\beta)\sigma^{-n_{l}}{(M+\zeta)}^{-n_{l}+\beta}}
{\Gamma(1-\beta)[{(\zeta+M)}^{\beta}-{\zeta}^{\beta}]}\times 
\frac{\tilde{\sigma}^{\beta}}
{\mathbb{E}[\tilde{\sigma}^{\beta}]}
$$
Hence it follows that the joint marginal distribution of $(N_{j,l}, j\in [J])$ can be expressed exactly as,
$$
\frac{{(1-\frac{\zeta}{{(\zeta+M)}})}^{n_{l}}}
{(1-\frac{\zeta^{\beta}}{{(\zeta+M)}^{\beta}})} \frac{\beta\Gamma(n_{l}-\beta)}{n_{l}!\Gamma(1-\beta)}\times \frac{n_{l}!}{\prod_{j=1}^{J}n_{j,l}!}
\mathbb{E}[\frac{\tilde{\sigma}^{\beta}}
{\mathbb{E}[\tilde{\sigma}^{\beta}]}\times \prod_{j=1}^{J}P^{n_{j,l}}_{j}]
$$
where, see for instance~\citep{PY97},
$$
\mathbb{E}\left[\frac{\tilde{\sigma}^{\beta}}
{\mathbb{E}[\tilde{\sigma}^{\beta}]}\times \prod_{j=1}^{J}P^{n_{j,l}}_{j}\right]=\mathbb{E}[\prod_{j=1}^{J}\tilde{P}^{n_{j,l}}_{\alpha,-\beta}(q_{j})]
$$
Hence,  the joint distribution of $((\tilde{C}_{j,k,l}, k\in[X_{j,l}]), [X_{j,l}], j\in[J]), (N_{j,l}, j\in[J]))$ is, for $x_{l}=\sum_{j=1}^{J}x_{j,l},$ $\sum_{k=1}^{x_{j,l}}c_{j,k,l}=n_{j,l},$ with $x_{j,l}=0$ for $n_{j,l}=0,$:
\begin{equation}
\frac{\Gamma(x_{l}-\frac{\beta}{\alpha})\Gamma(1-\beta)\alpha^{x_{l}-1}}{\Gamma(1-\frac{\beta}{\alpha})\Gamma(n_{l}-\beta)}
\prod_{j=1}^{J}q_{j}^{x_{j,l}}\prod_{k=1}^{x_{j,l}}\frac{\Gamma(c_{j,k,l}-\alpha)}{\Gamma(1-\alpha)}
\times\frac{\beta\Gamma(n_{l}-\beta)}{n_{l}!\Gamma(1-\beta)}\frac{{(1-\frac{\zeta}{{(\zeta+M)}})}^{n_{l}}}
{(1-\frac{\zeta^{\beta}}{{(\zeta+M)}^{\beta}})}.
\label{specialstablestablederived}
\end{equation}

\section{Details on the posterior inference algorithms for general Poisson-GG-GG HIBP}
\label{sec:post}

\subsection{Derivation of the marginal distribution}
In our simulation study, we consider the following Poisson-GG-GG-HIBP model with generic parameterizations,
\[
&\tau_0(\lambda) = \mathrm{GG}(\lambda\given\alpha_0, \theta_0/\alpha_0, \zeta_0) = \frac{\theta_0}{\Gamma(1-\alpha_0)}\lambda^{-\alpha_0-1}e^{-\zeta_0\lambda}, \quad B_0 \sim \mathrm{CRM}(\tau_0, F_0), \\
&\tau_j(s) = \mathrm{GG}(s\given \alpha_j, \theta_j/\alpha_j, \zeta_j), \quad B_j \given B_0 \sim \mathrm{CRM}(\tau_j, B_0) 
\text{ for } j \in [J],\\
&Z_j^{(i)}\given B_j \sim \mathrm{PoiP}(B_j) \text{ for } i \in [M_j], j \in [J].
\]
Before proceeding, let us calculate some quantities needed for the description,
\[
&\kappa_j := \psi_j(M_j) = \int_0^\infty (1-e^{-sM_j})\tau_j(s)ds = \frac{\theta_j}{\alpha_j}((\zeta_j+M_j)^{\alpha_j} - \zeta_j^{\alpha_j}), \\
&\bakappa := \sum_{j=1}^J \kappa_j, \quad q_j := \frac{\kappa_j}{\bakappa}, \quad
\Psi_0(\bakappa) = \frac{\theta_0}{\alpha_0}((\zeta_0 + \bakappa)^{\alpha_0} - \zeta_0^{\alpha_0}), 
\]
Let $Z_j := \sum_{i=1}^{M_j} Z_j^{(i)}$ and $\bZ_J := (Z_j, j \in [J])$. By Theorem 3.1, 
\[
\bZ_J \overset{d}{=} \left(
\sum_{l=1}^{\varphi} \left[
\sum_{k=1}^{X_{j,l}} \tC_{j,k,l} 
\right]\delta_{\tY_l}, i \in [M_j], j \in [J]
\right),
\]
where the corresponding marginal density can be computed as,
\[
\lefteqn{\bbP(\varphi=r, (X_{j,l}=x_{j,l}, \tC_{j,k,l}=c_{j,k,l}, j\in [J], k\in [x_{j,l}], l\in[r]))} \\
&= \bbP(\varphi=r) \prod_{l=1}^r \bbP(X_l=x_l) \prod_{j=1}^J \bbP(X_{j,l}=x_{j,l}\given X_l=x_l) 
\prod_{k=1}^{x_{j,l}} \bbP(C_{j,k,l}=c_{j,k,l}),
\]
where
\[
&\bbP(\varphi=r) = \frac{\Psi_0(\bakappa)^r e^{-\Psi_0(\bakappa)}}{r!},\\
&\bbP(X_l=x_l) = \frac{\1{x_l > 0} p_0^{x_l}}{1 - (1-p_0)^{\alpha_0}} \frac{\alpha_0\Gamma(x_l-\alpha_0)}{x_l!\Gamma(1-\alpha_0)}\text{ where } p_0 := \frac{\bakappa}{\zeta_0 + \bakappa}, \\
&\bbP((X_{j,l}=x_{j,l}, j\in[J])\given X_l=x) = \frac{x_{l}!}{\prod_{j=1}^J x_{j,l}!}\prod_{j=1}^J q_j^{x_{j,l}}, \\
&\bbP(\tC_{j,k,l}=c_{j,k,l}) = \frac{\1{c_{j,k,l}>0}p_j^{c_{j,k,l}}}{1 - (1-p_j)^{\alpha_j}} \frac{\alpha_j\Gamma(c_{j,k,l}-\alpha_j)}{c_{j,k,l}!\Gamma(1-\alpha_j)} \text{ where } p_j = \frac{M_j}{\zeta_j + M_j}.
\]
The marginal distribution is thus computed as,
\[
\lefteqn{\bbP(\varphi=r, (X_{j,l}=x_{j,l}, \tC_{j,k,l}=c_{j,k,l}, j\in [J], k\in [x_{j,l}], l\in[r]))} \\
&= \frac{\Psi_0(\bakappa)^r e^{-\Psi_0(\bakappa)}}{r!} \prod_{l=1}^r 
 \frac{\1{x_l>0}p_0^{x_l}}{1-(1-p_0)^{\alpha_0}} \frac{\alpha_0\Gamma(x_l-\alpha_0)}{x_l!\Gamma(1-\alpha_0)}
 \frac{x_{l}!}{\prod_{j=1}^J x_{j,l}!} \\
 &\times \prod_{j=1}^J \frac{ q_j^{x_{j,l}}\1{c_{j,k,l}>0}p_j^{c_{j,k,l}}}{1-(1-p_j)^{\alpha_j}} \frac{\alpha_j\Gamma(c_{j,k,l}-\alpha_j)}{c_{j,k,l}!\Gamma(1-\alpha_j)}.
\]
As  a note, recall from the main text, for $0<\alpha_{j}<1,$ let $P_{\alpha_{j},0} \sim \mathscr{PY}(\alpha_{j},0,F_{0})$ denote the $(\alpha_{j},0)$ Pitman-Yor process, \cite{IJ2001,Pit96}, constructed from normalizing a $\alpha_{j}$ stable-subordinator. The corresponding EPPF  describing the distribution of a random partition of $[n_{j,l}]$ is then given by:
$
p_{\alpha_{j},0}(\mathbf{c}_{j,l}) = \frac{\alpha^{x_{j,l}-1}_{j}\Gamma(x_{j,l})}{\Gamma(n_{j,l})}\prod_{k=1}^{x_{j,l}}\frac{\Gamma(c_{j,k,l}-\alpha_{j})}{\Gamma(1-\alpha_{j})}.
$
From this, the probability of the number of blocks $\mathbb{P}_{\alpha_{j},0}^{(n_{j,l})}(k) =
\frac{{\alpha^{k-1}_{j}\Gamma(k)} S_{\alpha_{j}}(n_{j,l},k)}{ {\Gamma(n_{j,l})}}$, where $S_{\alpha_{j}}(n_{j,l},k)$ denotes the generalized Stirling number of the second kind, as in \cite{Pit06}[3.2.3, p. 65-66].

In practice, as we only observe the aggregated counts $N_{j,l} = \sum_{k=1}^{X_{j,l}} \tC_{j,k,l}$, we have to work with the marginal distribution involving $(N_{j,l}, j\in [J], l \in [r])$.
From Proposition 5.4 and \cite{Pit97},
\[
\lefteqn{\bbP(\tN_{j,l}=n_{j,l}\given X_{j,l}=x_{j,l}) = \1{n_{j,l} \geq x_{j,l}}\sum_{*} \prod_{k=1}^{x_{j,l}} \bbP(\tC_{j,k,l}=c_{j,k,l})} \\
&= \1{n_{j,l} \geq x_{j,l}}\frac{\alpha_j^{x_{j,l}}p_j^{n_{j,l}}}{(1-(1-p_j)^{\alpha_j})^{x_{j,l}}}\sum_*  \prod_{k=1}^{x_{j,l}} \frac{\Gamma(c_{j,k,l}-\alpha_j)}{c_{j,k,l}!\Gamma(1-\alpha_j)} \\
&= \frac{\1{n_{j,l} \geq x_{j,l}}\alpha_j^{x_{j,l}} p_j^{n_{j,l}}}{(1 - (1-p_j)^{\alpha_j})^{x_{j,l}}}\frac{x_{j,l}! S_{\alpha_j}(n_{j,l}, x_{j,l})}{n_{j,l}!},
\]
where the summation $\sum_*$ is over all set of positive integers $(c_{j,1,l},\dots, c_{j,x_{j,l},l})$ partitioning $n_{j,l}$. The generalized Stirling number of the second kind $S_\alpha(n, x)$ is defined as follows:
\[
S_\alpha(n, x) = \frac{n!}{x!}\sum_{*} \prod_{k=1}^x \frac{\Gamma(c_k - \alpha)}{c_k!\Gamma(1-\alpha)},
\]
and can easily be computed via the following recursive formula,
\[
&S_{\alpha}(n, 1) = \frac{\Gamma(n-\alpha)}{\Gamma(1-\alpha)}, \quad S_{\alpha}(n, n) = 1, \\
&S_{\alpha}(n+1, x) = (n - \alpha x)S_\alpha(n, x) + S_\alpha(n, x-1).
\]
Based on this identity, the marginal distribution we will be working with is computed as,
\[\label{eq:poisson_gg_gg_hibp_marginal}
\lefteqn{\bbP(\varphi=r, (X_{j,l}=x_{j,l}, N_{j,l}=n_{j,l}, j\in [J], l\in [r]))} \\
&= \frac{\Psi_0(\bakappa)^r e^{-\Psi_0(\bakappa)}}{r!} \prod_{l=1}^r
\frac{\1{x_l>0} p_0^{x_l}}{1-(1-p_0)^{\alpha_0}} \frac{\alpha_0\Gamma(x_l-\alpha_0)}{\Gamma(1-\alpha_0)} \prod_{j=1}^J\frac{ q_j^{x_{j,l}}\1{n_{j,l}\geq x_{j,l}} \alpha_j^{x_{j,l}} p_j^{n_{j,l}}}{(1-(1-p_j)^{\alpha_j})^{x_{j,l}}} \frac{S_{\alpha_j}(n_{j,l}, x_{j,l})}{n_{j,l}!} \\
&= \frac{\theta_0^r e^{-\Psi_0(\bakappa)}}{r!(\zeta_0+\bakappa)^{x_\bigcdot - \alpha_0 r}}
\prod_{l=1}^r \frac{\1{x_l>0}\Gamma(x_l-\alpha_0)}{\Gamma(1-\alpha_0)}
\prod_{j=1}^J \frac{\1{n_{j,l}\geq x_{j,l}}\theta_j^{x_{j,l}}M_j^{n_{j,l}}}{(\zeta_j + M_j)^{n_{j,l}-\alpha_j x_{j,l}}} \frac{S_{\alpha_j}(n_{j,l}, x_{j,l})}{n_{j,l}!},
\]
where $x_\bigcdot := \sum_{l=1}^r x_l = \sum_{l=1}^r \sum_{j=1}^J x_{j,l}$.

\subsection{Detailed description of the sampler}
Let $\bn_J := (n_{j,l}, j\in[J], l \in [r])$ be the set of observed word counts. We aim to estimate the parameters $\phi := \{\alpha_0, \theta_0, \{\alpha_j, \theta_j\}_{j=1}^J\}$, with the parameters $\zeta_0$ and $\{\zeta_j\}_{j=1}^J$ fixed to $1$. We introduce vague priors on these parameters,
\[
\pi_0(\alpha_0) &:= \mathrm{logit}\, \calN(\alpha_0 \given 0, 1)
\propto \frac{1}{\alpha_0(1-\alpha_0)}\exp\bigg(-\frac{\Big(\log \frac{\alpha_0}{1-\alpha_0}\Big)^2}{2}\bigg),\\
\pi_0(\theta_0) &:= \mathrm{Gamma}(\theta_0\given 1, 1), \\
\pi_0(\alpha_j) &:= \mathrm{logit}\,\calN(\alpha_j\given 0, 1) \text{ for } j \in [J],\\
\pi_0(\theta_j) &:= \log \, \calN(\theta_j\given 0, 1) 
\propto \frac{1}{\theta_j}\exp\Big(-(\log\theta_j)^2\Big)
\text{ for } j \in [J], \\
\pi_0(\phi) &:= \pi_0(\alpha_0)\pi_0(\theta_0)\prod_{j=1}^J \pi_0(\alpha_j)\pi_0(\theta_j).
\]
The target posterior distribution to simulate is then defined as,
\[
\pi(\phi, \bX_J\given \bn_J)
\propto \bbP(\varphi=r, (X_{j,l}, \tN_{j,l}=n_{j,l}, j\in [J], l\in [r]))
\pi_0(\phi),
\]
where $\bX_J := (X_{j,l}, j\in [J], l \in [r])$. We proceed by doing Gibbs sampling the latent variables and the parameters $\phi$.

\paragraph{Sampling $\bX_J$.} The conditional distribution for each entry $X_{j,l}$ in $\bX_J$ is simply given as a discrete distribution,
\[
\bbP(X_{j,l}=x \given \text{others}) 
\propto \frac{\Gamma(x + x_l^{-j} - \alpha_0) \theta_j^x (\zeta_j+M_j)^{\alpha_j x} S_{\alpha_j}(n_{j,l}, x)}{(\zeta_0 + \bakappa)^{x}}, \quad x \in [n_{j,l}],
\]
where $x_l^{-j} := \sum_{j'\neq j} x_{j',l}$.

\paragraph{Sampling $\theta_0$.}
Given the gamma prior, the posterior of $\theta_0$ is also given as a Gamma distribution.
\[
\bbP(\theta_0 \in d\theta \given \mathrm{others}) = \mathrm{Gamma}\bigg(\theta \,\Big|\, r + a, b + \frac{(\zeta_0+\bakappa)^{\alpha_0}-\zeta_0^{\alpha_0}}{\alpha_0}\bigg).
\]

\paragraph{Sampling the other parameters.}
For the parameters other than $\theta_0$, unfortunately, the conditional distributions do not admit closed forms. Hence we sample them via random-walk metropolis-Hastings with the following proposal distributions.
\[
q(\alpha_0'\given \alpha_0) &= \mathrm{logit}\, \calN\bigg(\alpha_0' \,\Big|\,
\log\frac{\alpha_0}{1-\alpha_0}, (0.01)^2\bigg) \text{ for } j \in \{0\} \cup [J], \\
q(\theta_j'\given \theta_j) &= \mathrm{log}\, \calN(\theta_j'\given \log \theta_j, (0.01)^2) \text{ for } j \in [J].\]

\subsection{Prediction rules}

Let $Z_j^* := Z_j^{(M_j+1)}$ be a new document to be attached to the group $j$ in $\bZ_J$. 
By Proposition 5.5, with slight changes in the notations, 
\[
Z_j^* \given \bZ_J \overset{d}{=} \underbrace{\sum_{v=1}^{\varphi_j^*} \bigg[ \sum_{k=1}^{X_{j,v}^{(1)}}  \tC_{j,k,v}^{(1)}\bigg] \delta_{Y_{v}}}_{:=\,\, \Circled{1}} + \underbrace{\sum_{l=1}^r \bigg[
\sum_{k=1}^{X_{j,l}^{(2)}} \tC_{j,k,l}^{(2)} 
\bigg]\delta_{\tY_l}}_{:= \,\, \Circled{2}} + \underbrace{\sum_{l=1}^r \bigg[ \sum_{k=1}^{x_{j,l}} \tC^{(3)}_{j,k,l} \bigg] \delta_{\tY_l}.}_{:=\,\,\Circled{3} }
\]
Below, we describe each component in detail.

\paragraph{\Circled{1}.} Denote
\[
&\Delta_j := \psi_j(M_j+1) - \psi_j(M_j) = \frac{\theta_j}{\alpha_j}((\zeta_j + M_j+1)^{\alpha_j} - (\zeta_j+M_j)^{\alpha_j}).
\]
Then
\[
\bbP(\varphi_j^*=r_j^*) &\sim \mathrm{Poisson}\left(
\Psi_0(\bakappa+\Delta_j) - \Psi_0(\bakappa)\right) \\
&= \frac{(\Psi_0(\bakappa+\Delta_j)-\Psi_0(\bakappa))^{r_j^*} e^{-\Psi_0(\bakappa+\Delta_j)+\Psi_0(\bakappa)}}{(r_j^*)!}, \\
\bbP(X_{j,v}^{(1)}=x_{j,v}^{(1)}) &=
\int_0^\infty \frac{\theta_0\lambda^{-\alpha_0-1}e^{-(\zeta_0+\bakappa)\lambda} (1 - e^{-\Delta_j\lambda}) }{\Gamma(1-\alpha_0) (\Psi_0(\bakappa+\Delta_j) - \Psi_0(\bakappa))}
\frac{\1{x_{j,v}^{(1)} > 0} (\Delta_j\lambda)^{x_{j,v}^{(1)}}e^{-\Delta_j\lambda} }{(x_{j,v}^{(1)})!(1 - e^{-\Delta_j\lambda})} d \lambda \\
&= \frac{\1{x_{j,v}^{(1)}>0} (p_{0,j}^*)^{x_{j,v}^{(1)}} }{1 - (1-p_{0,j}^*)^{\alpha_0}}\frac{\alpha_0\Gamma(x_{j,v}^{(1)}-\alpha_0)}{(x_{j,v}^{(1)})!\Gamma(1-\alpha_0)} \text{ where } p_{0,j}^* := \frac{\Delta_j}{\zeta_0 + \bakappa + \Delta_j}, \\
\bbP(\tC_{j,k,v}^{(1)}=c_{j,k,v}^{(1)}) &= 
\int_0^\infty \frac{\theta_j s^{-\alpha_j-1}e^{-(\zeta_j+M_j)s}(1-e^{-s})}{\Gamma(1-\alpha_j)
\Delta_j} \frac{\1{c_{j,k,v}^{(1)}>0} s^{c_{j,k,v}^{(1)}}e^{-s} }{(1-e^{-s})} d s \\
&= \frac{\1{c_{j,k,v}^{(1)}>0} (p_j^*)^{c_{j,k,v}^{(1)}} }{1 - (1 - p_j^*)^{\alpha_j}}
\frac{\alpha_j\Gamma(c_{j,k,v}^{(1)}-\alpha_j)}{(c_{j,k,v}^{(1)})!\Gamma(1-\alpha_j)}
\text{ where } p_j^* := \frac{1}{\zeta_j+M_j+1}.
\]
For $\tN_{j,v}^{(1)} := \sum_{k=1}^{X_{j,k,v}^{(1)}} \tC_{j,k,v}^{(1)}$, we can compute the PMF as follows.
\[
\bbP(\tN_{j,v}^{(1)}=n_{j,v}^{(1)} \given X_{j,v}^{(1)}=x_{j,v}^{(1)}) = 
\frac{\1{n_{j,v}^{(1)}\geq x_{j,v}^{(1)}} \alpha_j^{x_{j,v}^{(1)}} (p_j^*)^{n_{j,v}^{(1)}} }{(1-(1-p_j^*)^{\alpha_j})^{x_{j,v}^{(1)}}} \frac{(x_{j,v}^{(1)})!S_{\alpha_j}(n_{j,v}^{(1)}, x_{j,v}^{(1)})}{n_{j,v}^{(1)}!}.
\]

\paragraph{\Circled{2}.}
From Proposition 5.5, we have
\[
\bbP(X_{j,l}^{(2)}=x_{j,l}^{(2)}) &= \int_0^\infty \mathrm{Gamma}(\lambda\given x_l - \alpha_0, \zeta_0 + \bakappa) \mathrm{Poisson}(x_{j,l}^{(2)}\given \Delta_j\lambda) d \lambda \\
&= \frac{\Gamma(x^{(2)}_{j,l}+x_l-\alpha_0)}{(x^{(2)}_{j,l})!\Gamma(x_l-\alpha_0)} (p_{0,j}^*)^{x_{j,l}^{(2)}} (1-p_{0,j}^*)^{x_l-\alpha_0}  \\
&= \mathrm{NB}(x_{j,l}^{(2)}\given x_l-\alpha_0, p_{0,j}^*), \\
\bbP(\tC_{j,k,l}^{(2)}=c_{j,k,l}^{(2)}) 
&= \frac{\1{c^{(2)}_{j,k,l}>0} (p_j^*)^{c_{j,k,l}^{(2)}}}{1-(1-p_j^*)^{\alpha_j}} 
\frac{\alpha_j\Gamma(c_{j,k,l}^{(2)}-\alpha_j)}{(c_{j,k,l}^{(2)})!\Gamma(1-\alpha_j)}, \\
\bbP(\tN_{j,l}^{(2)}=n_{j,l}^{(2)} \given X_{j,l}^{(2)}=x_{j,l}^{(2)}) &= 
\frac{\1{n_{j,l}^{(2)}\geq x_{j,l}^{(2)}} \alpha_j^{x_{j,l}^{(2)}} (p_j^*)^{n_{j,l}^{(2)}} }{(1-(1-p_j^*)^{\alpha_j})^{x_{j,l}^{(2)}}} \frac{(x_{j,l}^{(2)})!S_{\alpha_j}(n_{j,l}^{(2)}, x_{j,l}^{(2)})}{n_{j,l}^{(2)}!}.
\]

\paragraph{\Circled{3}.}
From Proposition 5.5, for $k \in [x_{j,l}]$,
\[
\bbP(\tC_{j,k,l}^{(3)}=c_{j,k,l}^{(3)}) = \mathrm{NB}(c_{j,k,l}^{(3)}\given 
c_{j,k,l}-\alpha_j, p_j^*)^{\1{x_{j,l}>0}}.
\]
Since a sum of negative binomial variables sharing a success probability is also a negative binomial,
\[
\bbP(\tN_{j,l}^{(3)}=n_{j,l}^{(3)}\given X_{j,l}=x_{j,l}) = \mathrm{NB}(n_{j,l}^{(3)}\given n_{j,l}-\alpha_jx_{j,l}, p_j^*)^{\1{x_{j,l}>0}}.
\]

To summarize, the predictive probability of a new document $Z_j^*$, computed for the aggregated counts, is given as
\[\label{eq:poisson_gg_gg_hibp_pred}
\lefteqn{\bbP( \varphi^*_j=r^*_j, (\tN_{j,v}^{(1)}=n_{j,v}^{(1)}, X_{j,v}^{(1)}=x_{j,v}^{(1)}, v \in [r^*_j]), (\tN_{j,l}^{(2)}=n_{j,l}^{(2)}, X_{j,l}^{(2)}=x_{j,l}^{(2)}, \tN^{(3)}_{j,l}=n_{j,l}^{(3)}, l \in [r]) \given \bZ_J)}  \\
&= \frac{(\Psi_0(\bakappa+\Delta_j)-\Psi_0(\bakappa))^{r_j^*} e^{-\Psi_0(\bakappa+\Delta_j)+\Phi_0(\bakappa)}}{(r_j^*)!}\\
&\times \prod_{v=1}^{r^*_j}
\frac{\1{x_{j,v}^{(1)}>0} (p_{0,j}^*)^{x_{j,v}^{(1)}} }{1 - (1-p_{0,j}^*)^{\alpha_0}}\frac{\alpha_0\Gamma(x_{j,v}^{(1)}-\alpha_0)}{(x_{j,v}^{(1)})!\Gamma(1-\alpha_0)} 
\frac{\1{n_{j,v}^{(1)}\geq x_{j,v}^{(1)}} \alpha_j^{x_{j,v}^{(1)}} (p_j^*)^{n_{j,v}^{(1)}} }{(1-(1-p_j^*)^{\alpha_j})^{x_{j,v}^{(1)}}} \frac{(x_{j,v}^{(1)})!S_{\alpha_j}(n_{j,v}^{(1)}, x_{j,v}^{(1)})}{n_{j,v}^{(1)}!} \\
&\times \prod_{l=1}^r \frac{\Gamma(x^{(2)}_{j,l}+x_l-\alpha_0)}{(x^{(2)}_{j,l})!\Gamma(x_l-\alpha_0)} (p_{0,j}^*)^{x_{j,l}^{(2)}} (1-p_{0,j}^*)^{x_l-\alpha_0}
\Bigg[\frac{\1{n_{j,l}^{(2)}\geq x_{j,l}^{(2)}} \alpha_j^{x_{j,l}^{(2)}} (p_j^*)^{n_{j,l}^{(2)}} }{(1-(1-p_j^*)^{\alpha_j})^{x_{j,l}^{(2)}}} \frac{(x_{j,l}^{(2)})!S_{\alpha_j}(n_{j,l}^{(2)}, x_{j,l}^{(2)})}{n_{j,l}^{(2)}!}\Bigg]^{\1{x_{j,l}^{(2)}>0}} \\
& \times \Bigg[\frac{\Gamma(n_{j,l}^{(3)}+n_{j,l}-\alpha_jx_{j,l})}{(n_{j,l}^{(3)})!\Gamma(n_{j,l}-\alpha_jx_{j,l})} (p_j^*)^{n_{j,l}^{(3)}}(1-p_j^*)^{n_{j,l}-\alpha_jx_{j,l}}\Bigg]^{\1{x_{j,l}>0}}.
\]

\subsection{Computing the predictive probability of a test document}
Now assume that we have observed a test document. In order to compute \eqref{eq:poisson_gg_gg_hibp_pred}, we need three types of counts $(\tN_{j,v}^{(1)}, \tN_{j,l}^{(2)}, \tN_{j,l}^{(3)})$ and latent variables $(X_{j,v}^{(1)}, X_{j,l}^{(2)}, X_{j,l}^{(3)})$, but from the document we can only get the counts
\[
\bn^* := (n_{j,v}^{(1)}, v\in [r_j^*]), (n_{j,l}^{(2,3)}, l \in [r]),
\]
where $n_{j,l}^{(2,3)} = n_{j,l}^{(2)} + n_{j,l}^{(3)}$. So in order to compute the predictive probability of $\bn^*$ given $\bZ_J$, we should first recover the latent variables and counts via the Gibbs sampling.

\paragraph{Sampling $X_{j,v}^{(1)}$.}
The conditional distribution of $X_{j,v}^{(1)}$ given others is a discrete distribution,
\[
\bbP(X_{j,v}^{(1)}=x\given\text{others})
\propto \frac{\Gamma(x-\alpha_0) \theta_j^x (\zeta_j+M_j+1)^{\alpha_j x}S_{\alpha_j}(n_{j,v}^{(1)},x)}{(\zeta_0+\bakappa+\Delta_j)^{x}}\text{ for } x \in [n_{j,v}^{(1)}].
\]

\paragraph{Sampling $X_{j,l}^{(2)}$, $\tN_{j,l}^{(2)}$ and $\tN_{j,l}^{(3)}$.}
Note that $\tN_{j,l}^{(2)} + \tN_{j,l}^{(3)} = n_{j,l}^{(2,3)}$. We have three possible cases:
\begin{enumerate}
    \item $x_{j,l}=0$ and $n_{j,l}^{(2,3)}=0$: in this case, the test document do not contain the $l$th word, so $X_{j,l}^{(2)}=\tN_{j,l}^{(2)}=\tN_{j,l}^{(3)}=0$ with probability 1.
    \item $x_{j,l}=0$ and $n_{j,l}^{(2,3)} > 0$: the $l$th word in the test document solely comes from the $\tN_{j,l}^{(2)}$, so $\tN_{j,l}^{(2)}=n_{j,l}^{(2,3)}$ and $\tN_{j,l}^{(3)}=0$ with probability 1. Given $\tN_{j,l}^{(2)}=n_{j,l}^{(2,3)}$, the conditional distribution of $X_{j,l}^{(2)}$ is given as a discrete distribution.
    \[
    \lefteqn{\bbP(X_{j,l}^{(2)}=x\given \tN_{j,l}^{(2)}=n_{j,l}^{(2,3)}, \text{others})}\\
    &\propto
    \frac{\Gamma(x + x_l - \alpha_0) \theta_j^x (\zeta_j + M_j+1)^{\alpha_jx}S_{\alpha_j}(n_{j,l}^{(2,3)}, x)}{(\zeta_0+\bakappa+\Delta_j)^x} \text{ for } x \in [n_{j,l}^{(2,3)}].
    \]
    \item $x_{j,l} > 0$ and $n_{j,l}^{(2,3)}>0$: the $l$th word in the test document either
    comes from both $\tN_{j,l}^{(2)}$ and $\tN_{j,l}^{(3)}$ or only from one of $\tN_{j,l}^{(2)}$ and $\tN_{j,l}^{(3)}$. The conditional distribution of $X_{j,l}^{(2)}$ is a discrete distribution,
    \[
    \lefteqn{\bbP(X_{j,l}^{(2)}=x\given \tN_{j,l}^{(2)}+\tN_{j,l}^{(3)}=n_{j,l}^{(2,3)}, \text{others})}\\
    &\propto
    \frac{\Gamma(x + x_l - \alpha_0) \theta_j^x (\zeta_j + M_j+1)^{\alpha_jx}
    (S_{\alpha_j}(n_{j,l}^{(2,3)}, x))^{\1{x>0}}
    }{(\zeta_0+\bakappa+\Delta_j)^x} \text{ for } x \in \{0\} \cup [n_{j,l}^{(2,3)}].
    \]
    If $x_{j,l}^{(2)}$ drawn from the above conditional is zero, we have
    $\tN_{j,l}^{(2)}=0$ and $\tN_{j,l}^{(3)}=n_{j,l}^{(2,3)}$ with probability 1. Otherwise,
    the conditional distribution for $\tN_{j,l}^{(2)}$ is given as,
    \[
    \lefteqn{\bbP(\tN_{j,l}^{(2)}=n \given \tN_{j,l}^{(2)}+\tN_{j,l}^{(3)}=n_{j,l}^{(2,3)}, X_{j,l}^{(2)}=x_{j,l}^{(2)}, \text{others})} \\
    &\propto \frac{S_{\alpha_j}(n, x_{j,l}^{(2)})}{n!}
    \frac{\Gamma(n_{j,l}^{(2,3)}-n + n_{j,l}-\alpha_j x_{j,l})}{(n_{j,l}^{(2,3)}-n)!}
    \text{ for } n \in [n_{j,l}^{(2,3)}].
    \]
    Given $n_{j,l}^{(2)}$ drawn from the above conditional, we have $n_{j,l}^{(3)}=n_{j,l}^{(2,3)}-n_{j,l}^{(2)}$ with probability 1.
\end{enumerate}

\subsection{Estimation results for the remaining parameters}
Here we present the parameter estimation results for the parameters $(\alpha_j, \theta_j, j \in [J])$.

\begin{figure}
\centering
\includegraphics[width=0.7\linewidth]{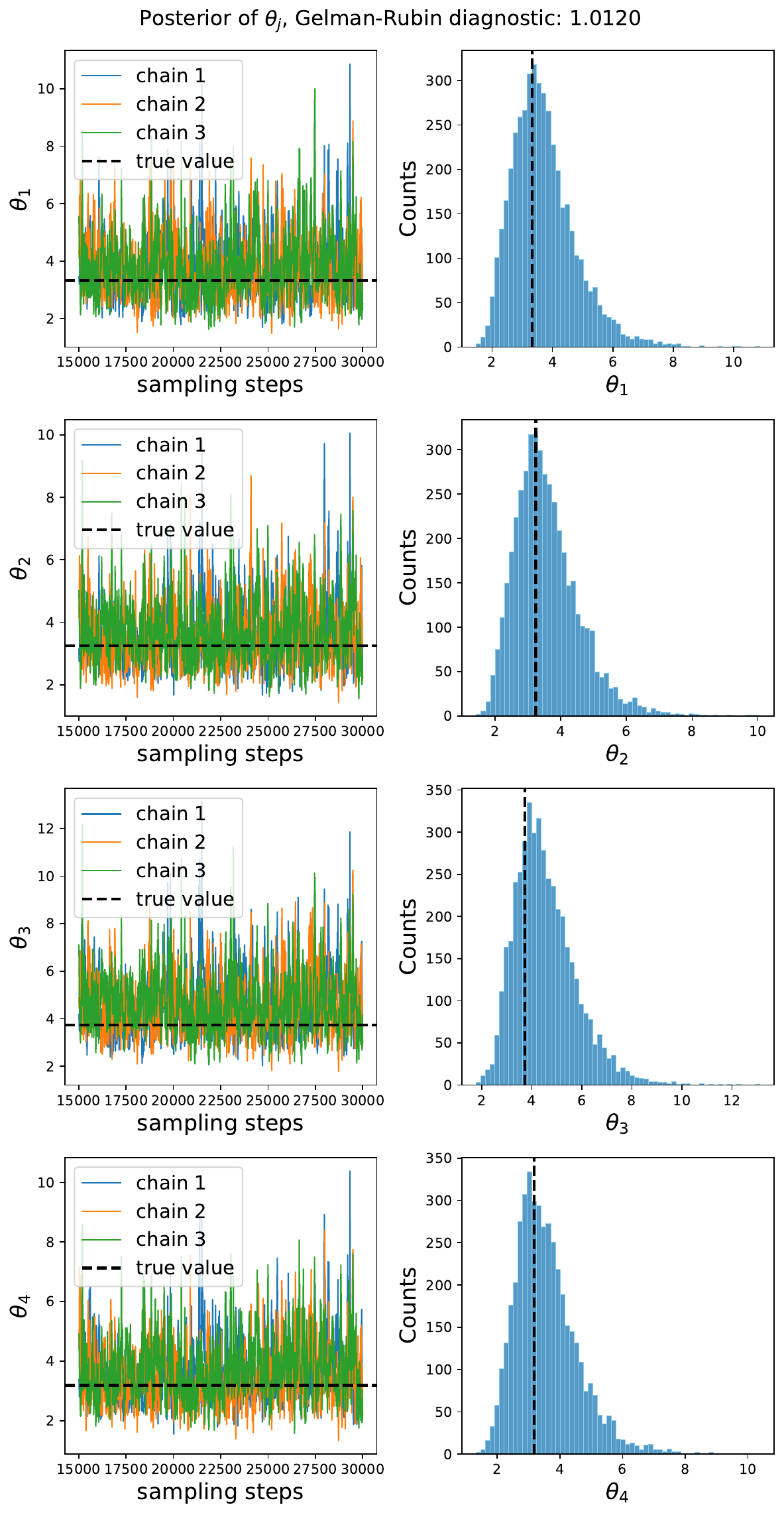}
\caption{Posterior samples of $\theta_j$ for the data generated with $\alpha=0.2$.}
\end{figure}

\begin{figure}
\centering
\includegraphics[width=0.7\linewidth]{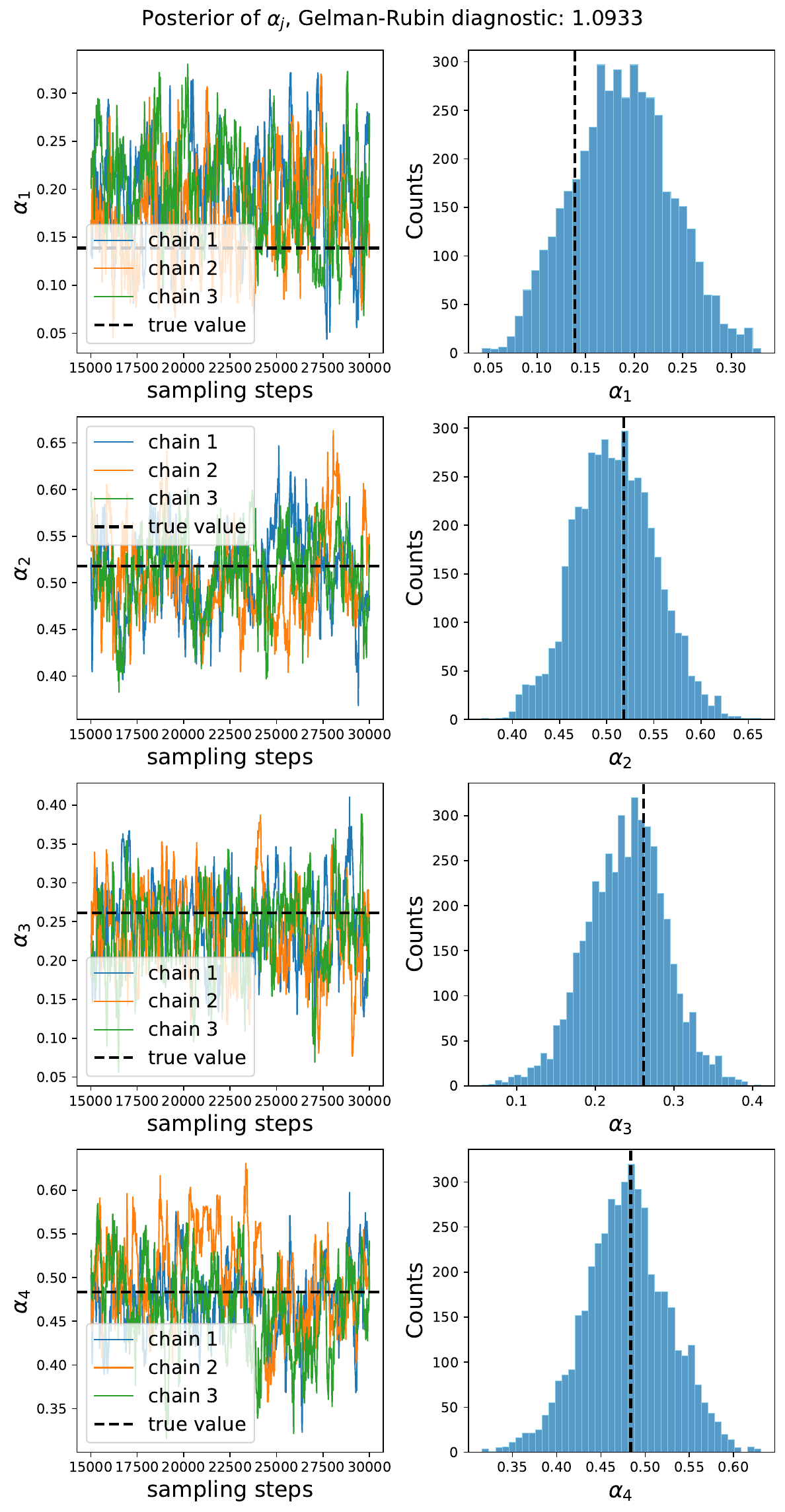}
\caption{Posterior samples of $\alpha_j$ for the data generated with $\alpha=0.2$.}
\end{figure}

\begin{figure}
\centering
\includegraphics[width=0.7\linewidth]{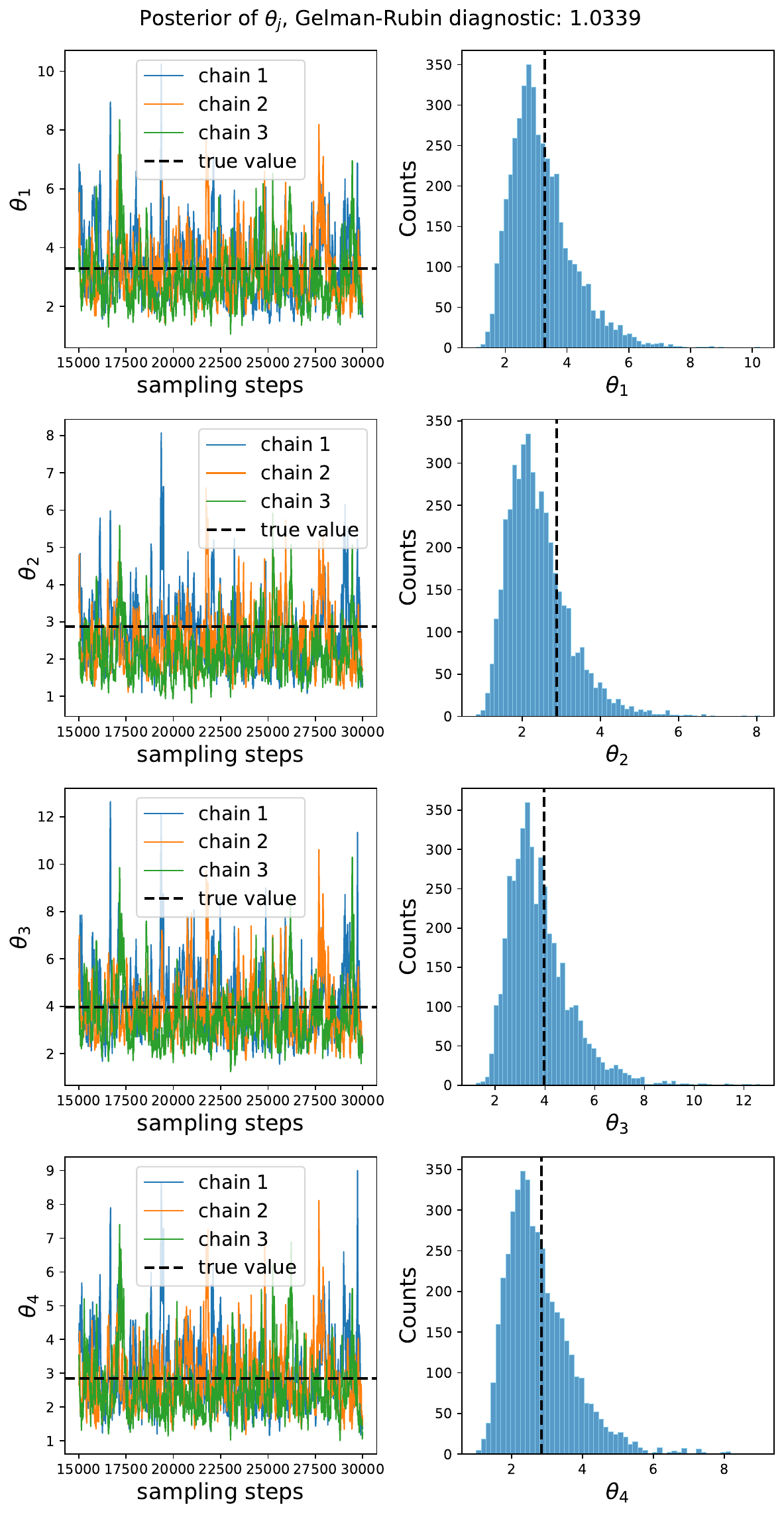}
\caption{Posterior samples of $\theta_j$ for the data generated with $\alpha=0.6$.}
\end{figure}

\begin{figure}
\centering
\includegraphics[width=0.7\linewidth]{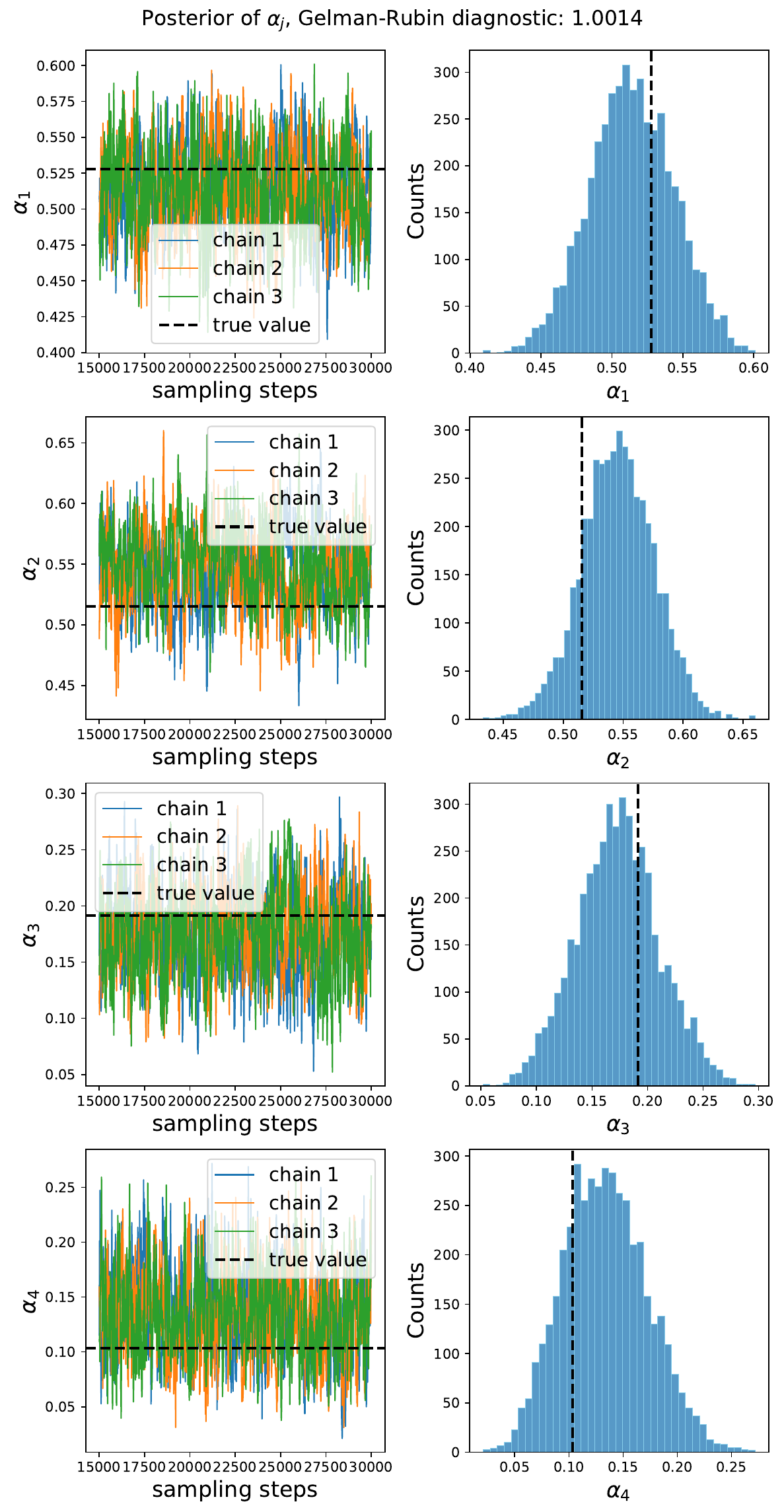}
\caption{Posterior samples of $\alpha_j$ for the data generated with $\alpha=0.6$.}
\end{figure}
\newpage
\bibliographystyle{agsm}
\bibliography{bibliography}

\end{document}